\documentclass[leqno,12pt]{article}

\usepackage[dvips]{graphics}
\usepackage{color}
\usepackage{ifthen}
\usepackage{soul}
\usepackage{latexsym}
\usepackage{amsthm,amssymb}
\usepackage{amsbsy,amsfonts,amsmath}

\newcommand{\bysame}{\leavevmode\hbox
to 3em{\hrulefill}\,}

\usepackage{a4}

\numberwithin{equation}{section}

\newtheorem{theorem}{Theorem}[section]
\newtheorem{proposition}[theorem]{Proposition}

\newtheorem{lemma}[theorem]{Lemma}
\newtheorem{example}[theorem]{Example}

\theoremstyle{definition}

\theoremstyle{remark}
\newtheorem{remark}[theorem]{Remark}

\newcommand{\al}{\alpha}

\newcommand{\trtvt}{\vartheta}
\newcommand{\trtfkAPi}{\fkA}

\newcommand{\chktrUp}{\chktrU^+}
\newcommand{\chktrUm}{\chktrU^-}

\newcommand{\chktrtvt}{{\check{\trtvt}}}

\newcommand{\chktrUgeq}{\chktrU^\geq}
\newcommand{\chktrUleq}{\chktrU^\leq}
\newcommand{\bC}{{\mathbb{C}}}
\newcommand{\bCt}{{\bC^\times}}
\newcommand{\mfkg}{{\mathfrak{g}}}
\newcommand{\mfks}{{\mathfrak{s}}}

\newcommand{\mcD}{{\mathbb{P}}}

\newcommand{\krdel}{\delta}

\newcommand{\pinfty}{\infty}

\newcommand{\Mat}{{\mathrm{Mat}}}

\newcommand{\bR}{\mathbb{R}}
\newcommand{\fkJ}{J}
\newcommand{\bZ}{\mathbb{Z}}
\newcommand{\bZgeqo}{\bZ_{\geq 0}}
\newcommand{\bN}{\mathbb{N}}

\newcommand{\bK}{\mathbb{K}}
\newcommand{\bKt}{\bK^\times}

\newcommand{\kpch}{{\hat{o}}}

\newcommand{\fka}{\mathfrak{a}}

\newcommand{\fkA}{\mathfrak{A}}
\newcommand{\bhm}{\chi}

\newcommand{\quoQ}{{\mathcal{Q}}}
\newcommand{\rmCh}{{\mathrm{Ch}}}

\newcommand{\rmrank}{{\mathrm{rank}}}

\newcommand{\rmGL}{{\mathrm{GL}}}

\newcommand{\tfkA}{\fkA}
\newcommand{\tal}{\al}
\newcommand{\tPi}{\Pi}

\newcommand{\tfkAPi}{\tfkA}
\newcommand{\tfkAPip}{\tfkA_\checktal^+} 
\newcommand{\fkI}{I}
\newcommand{\CardfkI}{{|\fkI|}}

\newcommand{\tq}{q}

\newcommand{\tbX}{{\mathbb{X}}}
\newcommand{\tbUo}{{\mathbb{U}}^0}
       
\newcommand{\tbhm}{\bhm}

\newcommand{\newpara}{\omega}

\newcommand{\newparatlmbtea}{\newpara^\tbhm_{\tlambda,\tmu;\beta}}

\newcommand{\tlambda}{\lambda}
\newcommand{\tmu}{\mu}
\newcommand{\tnu}{\nu}

\newcommand{\tC}{C}
\newcommand{\tc}{c}
\newcommand{\mcA}{{\mathcal{A}}}
\newcommand{\tvsigma}{\tau}

\newcommand{\mcC}{{\mathcal{C}}}
\newcommand{\mcV}{{\mathcal{V}}}

\newcommand{\tils}{s}

\newcommand{\mcH}{{\mathcal{H}}}

\newcommand{\rmMax}{{\mathrm{Max}}}

\newcommand{\tR}{R}
\newcommand{\tRp}{\tR^+}

\newcommand{\mcR}{{\mathcal{R}}}

\newcommand{\tK}{\trK} 
\newcommand{\tL}{\trL} 

\newcommand{\HopfD}{\Delta}
\newcommand{\HopfS}{S}
\newcommand{\Hopfe}{\varepsilon}

\newcommand{\rmSpan}{{\mathrm{Span}}}

\newcommand{\rmid}{{\rm {id}}}

\newcommand{\tbeta}{\beta}

\newcommand{\trU}{U}
\newcommand{\trUp}{U^+}
\newcommand{\trUm}{U^-}
\newcommand{\trUo}{U^0}

\newcommand{\trUtbhmtPi}{\trU(\tbhm,\checktal)}

\newcommand{\trUotbhmtPi}{\trUo(\tbhm,\checktal)}

\newcommand{\trK}{K}
\newcommand{\trL}{L}
\newcommand{\trE}{E}
\newcommand{\trF}{F}

\newcommand{\trzeta}{\zeta}

\newcommand{\trSh}{\mathfrak{Sh}}

\newcommand{\trBp}{{\trU^{+,\flat}}}
\newcommand{\trBm}{{\trU^{-,\flat}}}

\newcommand{\prtrZ}{\mathfrak{Z}}
\newcommand{\prtrZnewpara}{\prtrZ_\newpara}
\newcommand{\prtrZnewparatbhmtPi}{\prtrZnewparatbhmtPi(\tbhm,\tPi)}

\newcommand{\trHCmap}{{\mathfrak{HC}}}

\newcommand{\trpilq}{\wp_\leq}

\newcommand{\rmIm}{{\mathrm{Im}}}

\newcommand{\trLambda}{\Lambda}
\newcommand{\tv}{{\tilde v}}
\newcommand{\mclM}{{\mathcal{M}}}
\newcommand{\mclN}{{\mathcal{N}}}
\newcommand{\mclL}{{\mathcal{L}}}

\newcommand{\MatRank}{{\mathrm{Rank}}}

\newcommand{\trR}{R}
\newcommand{\trRtbhmtPi}{\trR(\tbhm,\checktal)}

\newcommand{\trRp}{\trR^+}
\newcommand{\trRptbhmtPi}{\trRp(\tbhm,\checktal)}

\newcommand{\trvrp}{\varphi}

\newcommand{\originalCardtrRp}{|\trRptbhmtPi|}
\newcommand{\CardtrRp}{\kappa}
\newcommand{\CardtrRptbtP}{|\trRptbhmtPi|}

\newcommand{\trN}{N}

\newcommand{\tPiangi}{{\gravetau_i(\checktal)}}

\newcommand{\trT}{T}

\newcommand{\nblIso}{\nabla}

\newcommand{\trRpmap}{\mathfrak{P}}

\newcommand{\dotbeta}{{\dot \beta}}

\newcommand{\dottrE}{{\dot \trE}}
\newcommand{\dottrF}{{\dot \trF}}

\newcommand{\trjmth}{\gamma}
\newcommand{\trho}{\hat \rho}

\newcommand{\tfkB}{{\mathfrak{B}}}

\newcommand{\trLambdaangi}{{\trLambda^{\tbhm,\checktal,i}}}

\newcommand{\hT}{{\mathcal{T}}}

\newcommand{\dkpch}{{\hat c}}
\newcommand{\ddottrE}{{\ddot{\trE}}}
\newcommand{\ddottrF}{{\ddot{\trF}}}

\newcommand{\mclS}{{\mathbb{S}}}

\newcommand{\trUoSetmsf}{{\trUo\setminus\trUo f}}

\newcommand{\trGamma}{\Gamma}

\newcommand{\mclG}{{\mathcal{G}}}
\newcommand{\mclZ}{{\mathcal{Z}}}

\newcommand{\rmgcd}{{\mathrm{gcd}}}

\newcommand{\bKk}{\bK^k}
\newcommand{\bKtk}{(\bKt)^k}
\newcommand{\KbKk}{\bK[\bKk]}
\newcommand{\KbKtk}{\bK[\bKtk]}

\newcommand{\rmX}{{\mathrm{X}}}
\newcommand{\mclV}{{\mathcal V}}
\newcommand{\mclD}{{\mathcal D}}
\newcommand{\rmx}{{\mathrm{x}}}
\newcommand{\mclF}{{\mathcal F}}

\newcommand{\bKtprim}{\bK_{<\infty}}

\newcommand{\bKkgqi}{{\bK^k_{<\infty}}}

\newcommand{\trxi}{\xi}

\newcommand{\bKtll}{{(\bKt)^{\CardfkIN^2}}}

\newcommand{\bV}{{\mathbb{V}}}
\newcommand{\bD}{{\mathbb{D}}}

\newcommand{\tbXtbhmone}{{\tbX^\tbhm_1}}
\newcommand{\tbXtbhmtwo}{{\tbX^\tbhm_2}}

\newcommand{\tbXtprim}{\tbX_{<\infty}}

\newcommand{\tbXtbhmtlambda}{{\tbX^{\tbhm,\tlambda}_3}}

\newcommand{\ChtbUo}{{\rmCh(\tbUo)}}

\newcommand{\treta}{\eta}
\newcommand{\bKttwl}{{(\bKt)^{2\CardfkIN}}}
\newcommand{\bKtlltwl}{{(\bKt)^{\CardfkIN^2+2\CardfkIN}}}

\newcommand{\bY}{{\mathbb{Y}}}

\newcommand{\bari}{{\bar i}}
\newcommand{\barj}{{\bar j}}
\newcommand{\tromega}{\epsilon}
\newcommand{\bKkkgqi}{{\bK_{<\infty}^{2k}}}

\newcommand{\bG}{{\mathbb{G}}}

\newcommand{\trUodag}{\trU^{0,\dagger}}
\newcommand{\trUdag}{\trU^\dagger}

\newcommand{\chkI}{{\check{\mathcal{I}}}}
\newcommand{\chktrU}{{\check{\trU}}}
\newcommand{\chkpi}{{\check{p}}}
\newcommand{\chktrK}{{\check{\trK}}}

\newcommand{\chktrE}{{\check{\trE}}}
\newcommand{\chktrF}{{\check{\trF}}}
\newcommand{\chktrUo}{{{\check{\trU}}^0}}

\newcommand{\chktrXp}{{{\check{X}}^+}}
\newcommand{\chktrXm}{{{\check{X}}^-}}
\newcommand{\chktrY}{{\check{Y}}}

\newcommand{\chkprtrZ}{{\check{\prtrZ}}}

\newcommand{\chktrSh}{{\check{\trSh}}}
\newcommand{\chktrHCmap}{{\check{\trHCmap}}}

\newcommand{\chktfkB}{{\check{\tfkB}}}

\newcommand{\eparity}{\partial}
\newcommand{\rmAutbZ}{{\mathrm{Aut}}_\bZ}

\newcommand{\phiprime}{\Upsilon}

\newcommand{\trRiso}{\trR_{{\rm{iso}}}}
\newcommand{\trRniso}{\trR_{{\rm{niso}}}}
\newcommand{\espri}{s^\prime}
\newcommand{\eWprime}{W^\prime}
\newcommand{\trhoprime}{\trho^\prime}

\newcommand{\CardfkIN}{\theta}

\newcommand{\newmcV}{\mcV}
\newcommand{\bZleqo}{\bZ_{\leq 0}}
\newcommand{\tRptPi}{\tRp_\tPi}
\newcommand{\tRm}{\tR^-}
\newcommand{\tRmtPi}{\tRm_\tPi}
\newcommand{\tB}{{\mathcal{B}}}
\newcommand{\tiltB}{{\tilde {\tB}}}

\newcommand{\tN}{N}
\newcommand{\tgamma}{\gamma}
\newcommand{\checktB}{{\check{\tB}}}
\newcommand{\checktal}{\pi}
\newcommand{\checktau}{{\check{\tvsigma}}}
\newcommand{\checkV}{{\bar{\newmcV}}}
\newcommand{\checkv}{{\bar{\tal}}}
\newcommand{\bbM}{{\mathbb{B}}}
\newcommand{\bbH}{{\mathcal{H}}}

\newcommand{\checkmcC}{{\check{\mcC}}}
\newcommand{\checkmcR}{{\check{\mcR}}}

\newcommand{\checkmcCctPi}{\checkmcC[\checktal]}
\newcommand{\checkmcRctPi}{\checkmcR[\checktal]}

\newcommand{\checkmcCctPia}{\checkmcC[\checktal^a]}
\newcommand{\checkmcRctPia}{\checkmcR[\checktal^a]}

\newcommand{\bartPi}{{\bar{\checktal}}}

\newcommand{\checkmapu}{{\check{\mapu}}}
\newcommand{\trell}{\ell}

\newcommand{\mapu}{\phi}

\newcommand{\mfkf}{{\mathfrak{F}}}

\newcommand{\gravetau}{\grave{\tvsigma}}

\newcommand{\trNtbhmtPi}{\trN^{\tbhm,\checktal}}

\newcommand{\tBtbhmchPi}{{\tiltB(\tbhm,\checktal)}}
\newcommand{\checktBtbhmchPi}{{\tiltB^\vee(\tbhm,\checktal)}}

\newcommand{\hatf}{{\dot f}}

\newcommand{\trUtbhmtPiith}{\trU(\tbhm,\checktal;i)}

\newcommand{\trUptbhmtPi}{\trUp(\tbhm,\checktal)}
\newcommand{\trUmtbhmtPi}{\trUm(\tbhm,\checktal)}

\newcommand{\trtvttbhmtPi}{\trtvt^{\tbhm,\checktal}}
\newcommand{\trvrptbhmtchecktal}{\trvrp^{\tbhm,\checktal}}

\newcommand{\trvrpptbhmtchecktal}{\trvrp_+^{\tbhm,\checktal}}

\newcommand{\Omegatbhmtchecktal}{\Omega^{\tbhm,\checktal}}

\newcommand{\rmChtbhmchecktal}{\rmCh(\trUotbhmtPi)}

\newcommand{\Mtbctal}{\mclM^{\tbhm,\checktal}}
\newcommand{\MtbctaltrLam}{\Mtbctal(\trLambda)}
\newcommand{\Ltbctal}{\mclL^{\tbhm,\checktal}}
\newcommand{\LtbctaltrLam}{\Ltbctal(\trLambda)}
\newcommand{\NtbctaltrLam}{\mclN^{\tbhm,\checktal}(\trLambda)}

\newcommand{\ckpchi}{{c_i}}

\newcommand{\trhobctal}{\trho^{\tbhm,\checktal}}

\newcommand{\trShbctal}{\trSh^{\tbhm,\checktal}}

\newcommand{\mclSbctaltl}{\mclS^{\tbhm,\checktal}_\tlambda}

\newcommand{\hatX}{{\hat{X}}}
\newcommand{\hatY}{{\hat{Y}}}

\newcommand{\trRpmapbctal}{\trRpmap^{\tbhm,\checktal}}

\newcommand{\bKtkC}{(\bKt)^\CardfkIN}

\newcommand{\checktBtbhmchPictal}{\checktBtbhmchPi[\checktal]}

\newcommand{\trpilqtbtal}{\trpilq^{\tbhm,\checktal}}

\newcommand{\prtrZnewparabctal}{\prtrZnewpara(\tbhm,\checktal)}
\newcommand{\prtrZnewparabctalnangi}{\prtrZnewpara(\tbhm,\tPiangi)}

\newcommand{\trHCmapbctaln}{\trHCmap^{\tbhm,\checktal}_\newpara}
\newcommand{\trHCmapbctalnangi}{\trHCmap^{\tbhm,\tPiangi}_\newpara}

\newcommand{\trTtbtali}{\trT^{\tbhm,\tPiangi}_i}
\newcommand{\smtrTtbtali}{\trT^{(\tbhm,\checktal,\checktal^\prime;i)}}

\newcommand{\hTtbtali}{\hT^{\tbhm,\tPiangi}_i}

\newcommand{\trRtbtal}{\trR(\tbhm,\tPiangi)}

\newcommand{\trjmthtbtalni}{\trjmth^{\tbhm,\tPiangi}_{\newpara;i}}
\newcommand{\trUtbtalni}{\trU(\tbhm,\tPiangi)}
\newcommand{\trUotbtalni}{\trUo(\tbhm,\tPiangi)}
\newcommand{\trUptbtalni}{\trUp(\tbhm,\tPiangi)}
\newcommand{\trUmtbtalni}{\trUm(\tbhm,\tPiangi)}

\newcommand{\trShbctalni}{\trSh^{\tbhm,\tPiangi}}

\newcommand{\acuUtlambda}{{\acute{U}}[\tlambda]}

\newcommand{\tfkBbctanp}{\tfkB^{\tbhm,\checktal}_\newpara}
\newcommand{\tfkBbctanptb}{\tfkBbctanp(\tbeta)}

\newcommand{\tfkBbctanpbctal}{\tfkB^{\tbhm,\tPiangi}_\newpara}
\newcommand{\tfkBbctanpbctaltb}{\tfkBbctanpbctal(\tbeta)}

\newcommand{\trhobctaltb}{\trhobctal(\tbeta)}
\newcommand{\hckpchita}{c_\tal}
\newcommand{\hckpchitb}{c_\tbeta}

\newcommand{\onebctaltrT}{1^{\tbhm,\checktal}\trT}

\newcommand{\tfkAPipni}{\tfkAPi^+_\tPiangi}

\newcommand{\tfkAPiptlambdast}{{\tfkA_\checktal^{+,<\tlambda}}}

\newcommand{\mtlambda}{{m_\tlambda}}
\newcommand{\mtmu}{{m_\tmu}}

\newcommand{\dhXhY}{d_{\hatX,\hatY}}

\newcommand{\remclS}{{\mathcal{S}}}
\newcommand{\retmclS}{{\tilde{\remclS}}}

\newcommand{\talpii}{{\checktal(i)}}
\newcommand{\talpij}{{\checktal(j)}}
\newcommand{\talpik}{{\checktal(k)}}
\newcommand{\talpir}{{\checktal(r)}}
\newcommand{\talpit}{{\checktal(t)}}
\newcommand{\talpidfone}{{\checktal(\hatf(1))}}

\newcommand{\addsig}{\sigma}

\begin{document}

\begin{center}
{\Huge{Centers of Generalized Quantum Groups}}
\end{center}

\vspace{1cm}

\begin{center}
Punita Batra\,\,and\,\,Hiroyuki Yamane$^*$
\end{center}

\footnote[0]{$^*$ The corresponding author is Hiroyuki Yamane. {\it{2000 Mathematics Subject Classification.}}
Primary 17B10, Secondary 17B37}

\begin{abstract}
This paper treats the generalized quantum group $\trU=\trUtbhmtPi$
with a bi-homomorphism $\tbhm$ for which the corresponding generalized root
system is a finite set.
We establish a Harish-Chandra type theorem describing the
(skew) centers of $\trU$.
\end{abstract}

Keywords: Quantum groups, Lie superalgebras, Nichols algebras

\section{Introduction}\label{section:Intr}
In this paper, we study the center and skew-centers of 
the generalized quantum groups.
To do so, we use some facts of Zariski topology.
We introduce a sub-topology of it,
and call it {\it{Nichols topology}}.
Study of the topology has virtually been initiated by \cite{HY10}.
The main result of this paper is Theorem~\ref{theorem:mainTh}. As a bi-product, we also 
give its `symmetric' version Theorem~\ref{theorem:SYMMAIN}. Theorem~\ref{theorem:Shapo} is just \cite[Theorem~7.3]{HY10}.
Theorem~\ref{theorem:SVth} and Theorem~\ref{theorem:basicc} are similar to \cite[Lemma~6.7]{HY10} and
\cite[Lemma~9.4]{HY10} respectively.
We follow the Kac argument~\cite{Kac84} for the center. However 
the proof of Theorem~\ref{theorem:ndgthCor} is very delicate
since we treat `all roots of unit cases.' 
The proofs of Proposition~\ref{proposition:spsibsc} and Lemmas~\ref{lemma:HCinj} and \ref{lemma:estlm} are original.
We adopt the new and simpler definition \cite{Y15} for the generalized root systems (see \eqref{eqn:defstGRS}),
and in Section~\ref{section:setGRSandGQG}, we clarify relation between 
the new and conventional ones.
Argument in Subsection~\ref{subsection:LCSbh} has originally been given in \cite[Section~7]{HY10}.
We need delicate arguments in Section~\ref{section:ShFac} to obtain Theorem~\ref{theorem:ndgthCor}.

Let $\bK$ be an algebraically closed field. Let $\bKt:=\bK\setminus\{0\}$.
Let $\fkI$ be a finite set.
Let $\trtfkAPi$ be a free $\bZ$-module with $\CardfkI=\rmrank_\bZ\trtfkAPi$.
Let $\tbhm:\trtfkAPi\times\trtfkAPi\to\bKt$ be a bi-homomorphism, that is,
\begin{equation}\label{eqn:Intrtbhm}
\tbhm(\tlambda,\tmu+\tnu)=\bhm(\tlambda,\tmu)\tbhm(\tlambda,\tnu)\quad
\mbox{and}\quad\bhm(\tlambda+\tmu,\tnu)=\bhm(\tlambda,\tnu)\bhm(\tmu,\tnu)
\end{equation}
for all $\tlambda$, $\tmu$, $\tnu\in\trtfkAPi$.
In Introduction, for simplicity, we also assume
$\tbhm(\tlambda,\tmu)=\tbhm(\tmu,\tlambda)$ for all 
$\tlambda$, $\tmu\in\trtfkAPi$.
Let $\checktal:\fkI\to\trtfkAPi$ be a map such that
$\checktal(\fkI)$ is a $\bZ$-base of $\trtfkAPi$.
In the same way as in the Lusztig's definition \cite{b-Lusztig93}
of the quantum groups, to the pair $(\tbhm,\checktal)$, we can associate  
a unique $\bK$-algebra $\chktrU(\tbhm,\checktal)$
characterized by the following properties.
\newline\par
${\rm{(i)}}$ 
As a $\bK$-algebra, $\chktrU$ has generators
$\chktrK_\tlambda$ ($\tlambda\in\trtfkAPi$),
$\chktrE_i$, $\chktrF_i$ ($i\in \fkI$)
satisfying the equations
$\chktrK_0=1$, 
$\chktrK_\tlambda \chktrK_\tmu=\chktrK_{\tlambda+\tmu}$,
$\chktrK_\tlambda \chktrE_i =\tbhm(\tlambda,\talpii)\chktrE_i \chktrK_\tlambda$,
$\chktrK_\tlambda \chktrF_i  =\tbhm(\tlambda,-\talpii)\chktrF_i \chktrK_\tlambda$,
$\chktrE_i\chktrF_j-\chktrF_j\chktrE_i=\delta_{ij}(-\chktrK_\talpii+\chktrK_{-\talpii})$.
As a $\bK$-linear space, the elements $\chktrK_\tlambda$ ($\tlambda\in\trtfkAPi$)
are linearly independent.
\par
${\rm{(ii)}}$ 
Let $\chktrUo$ be the subalgebra of 
$\chktrU$ generated by $\chktrK_\tlambda$ ($\tlambda\in\trtfkAPi$).
Let $\chktrUp$ (resp. $\chktrUm$) be the subalgebra of 
$\chktrU$ generated by
$\chktrE_i$ (resp. $\chktrF_i$)
($i\in \fkI$) and $1$.
One has a $\bK$-linear isomorphism
$\chktrUm\otimes_\bK\chktrUo\otimes_\bK\chktrUp
\to\chktrU$ with $Y\otimes Z\otimes X\mapsto YZX$.
One also has $\bK$-linear subspaces $\chktrU_\tlambda$
($\tlambda\in\trtfkAPi$)
such that $\chktrU=\oplus_{\tlambda\in\trtfkAPi}\chktrU_\tlambda$,
$\chktrU_\tlambda\chktrU_\tmu\subset\chktrU_{\tlambda+\tmu}$,
$\chktrUo\subset\chktrU_0$, and 
$\chktrE_i\in\chktrU_\talpii$, $\chktrF_i\in\chktrU_{-\talpii}$
($i\in \fkI$).
\par
${\rm{(iii)}}$ 
Let $\chktrUgeq$ (resp. $\chktrUleq$) be the subalgebra of 
$\chktrU$ generated by $\chktrUp$ (resp. $\chktrUm$)
and $\chktrUo$.
One has a Drinfeld bilinear map
$\chktrtvt=\chktrtvt^{\tbhm,\checktal}:
\chktrUgeq\times\chktrUleq\to\bK$ such that
$\chktrtvt_{|\chktrUp\times\chktrUm}:\chktrUp\times\chktrUm\to\bK$ is non-degenerate,
and $\chktrtvt(\chktrK_\tlambda,\chktrK_\tmu)=\tbhm(\tlambda,-\tmu)$,
$\chktrtvt(\chktrE_i,\chktrF_j)=\delta_{ij}$, 
$\chktrtvt(\chktrE_i,\chktrK_\tmu)=\chktrtvt(\chktrK_\tlambda,\chktrF_j)=0$.
\newline\newline
We note that for each pair $(\tbhm,\checktal)$,
one define
the Karchenko's root system $\trRtbhmtPi$,
which is a subset of
$\trtfkAPi$ including $\checktal(\fkI)$;
$\trRtbhmtPi$ is used to give a
PBW-type theorem of $\chktrU(\tbhm,\checktal)$,
see Theorem~\ref{theorem:KhaPBW}. 

\begin{example} \label{example:ExaInto}
{\rm{
Assume that $\bK$ is the complex field $\bC$.
Let $\tq\in\bCt\setminus\{1\}$.
Let
$\langle\,,\,\rangle:\tfkAPi\times\tfkAPi\to\bZ$
be a $\bZ$-module bihomomorphism.

(1) Let $\mfkg$ be a complex simple Lie algebra.
Assume that $\langle\talpii,\talpii\rangle\ne 0$ ($i\in\fkI$) and
$\Bigl[{\frac {2\langle\talpii,\talpij\rangle} {\langle\talpii,\talpii\rangle}}\Bigr]_{i,j\in\fkI}$ coincides with
the Cartan matrix of
$\mfkg$.
Assume that $\tq$ is not a root of unity, and that
$\tbhm(\talpii,\talpij)=\tq^{\langle\talpii,\talpij\rangle}$
($i$, $j\in \fkI$). Then
$\chktrU$ can be identified  with the quantum group
$U_\tq(\mfkg)$,
and $\trRtbhmtPi$ can be identified with
the root system of $\mfkg$.

(2) Assume that $\langle\,,\,\rangle$ is as in (1).
Assume that $\tq$ is an $n$-th primitive root of unity for some $n\in\bN$
with $n\geq 2$.
Assume that 
$\tbhm(\talpii,\talpij)=\tq^{\langle\talpii,\talpij\rangle}$
($i$, $j\in \fkI$), $\tbhm(\talpii,\talpii)\ne 1$ ($i\in \fkI$)
and $\tbhm(\talpii,\talpii)\ne \tbhm(\talpij,\talpij)$
($i$, $j\in \fkI$ with 
$\langle\talpii,\talpii\rangle\ne\langle\talpij,\talpij\rangle$).
Let $X$ be the two-sided ideal of $\chktrU$ generated by
$\chktrK_\talpii^n-1$ ($i\in \fkI$).
Then the quotient algebra $\chktrU/X$ can be identified with
the Lusztig's small quantum group $u_\tq(\mfkg)$.
It follows that $\trRtbhmtPi$ is the same as in (1).

(3) (See also Subsection~\ref{subsection:SymExamples}.) Let $\mfks$ be a basic classical Lie superalgebra
(which is also called a complex simple Lie supserlagebra of type A-G).
Identify $\trtfkAPi$ with the root lattice of $\mfks$
for which $\checktal(\fkI)$ is also a set of simple roots.
Identify $\langle\,,\,\rangle$ with the Killing form of $\mfks$.
Let $\fkI^\prime:=\{i\in\fkI| \mbox{ $\talpii $ is an odd simple root} \}$. 
Define the $\bZ$-module homomorphism
$\eparity:\trtfkAPi\to\bZ$ by $\eparity(\talpii):=0$
($i\in\fkI\setminus\fkI^\prime$) and $\eparity(\talpij):=1$
($j\in\fkI^\prime$).
Assume that $\tq$ is not a root of unity, and that
$\tbhm(\talpii,\talpij)=(-1)^{\eparity(\talpii)\eparity(\talpij)}
\tq^{\langle\talpii,\talpij\rangle}$
($i$, $j\in \fkI$). Let
$\chktrU^\addsig=\chktrU\oplus\chktrU\addsig$
be the $\bC$-algebra obtained from $\chktrU$ by
adding an element $\addsig$ such that $\addsig^2=1$,
$\addsig\chktrK_\tlambda\addsig=\chktrK_\tlambda$, 
$\addsig\chktrE_i\addsig=(-1)^{\eparity(\talpii)}\chktrE_i$, and 
$\addsig\chktrF_i\addsig=(-1)^{\eparity(\talpii)}\chktrF_i$ ($i\in\fkI$).
Then the quantum superalgebra 
$U_\tq(\mfks)$ can be identified with 
the subagebra of $\chktrU^y$
generated by $\addsig^{\eparity(\talpii)}\chktrK_\talpii$,
$\chktrE_i$, and $\chktrF_i\addsig^{\eparity(\talpii)}$,
so $\chktrU^\addsig=U_\tq(\mfks)\oplus U_\tq(\mfks)\addsig$. 
It follows that 
the root system of $\mfks$ is
given by $\trRtbhmtPi
\cup\{2\tbeta|\tbeta\in\trRtbhmtPi,(-1)^{\eparity(\tbeta)}=-1,
\langle\tbeta,\tbeta\rangle\ne 0\}$.
}}
\end{example}
Let $\newpara:\trtfkAPi\to\bKt$ be a $\bZ$-module homomorphism. 
We call the $\bK$-linear space
\begin{equation*}
\chkprtrZ_\newpara=\chkprtrZ_\newpara(\tbhm,\checktal)
:=\{\,Z\in\chktrU_0\,|\,\forall \tlambda \in\trtfkAPi,\,\forall X\in\chktrU_\tlambda,\,
ZX
=\newpara(\tlambda)XZ\,\}
\end{equation*} the {\it{$\newpara$-skew graded center}} of $\chktrU$.
Define the $\bK$-linear map $\chktrSh=\chktrSh^{\tbhm,\checktal}:\chktrU\to\chktrUo$
by $\chktrSh_{|\chktrUo}=\rmid_{\chktrUo}$ and $\chktrSh(\chktrU\chktrE_i)=\chktrSh(\chktrF_i\chktrU)=\{0\}$
($i\in\fkI$).
Let $\chktrHCmap_\newpara=\chktrHCmap^{\tbhm,\checktal}_\newpara:=\chktrSh_{|\chkprtrZ_\newpara}$.
It is easy to see that 
$\chktrHCmap_\newpara$ is injective.

The $\chktrU$ is also appeared in the context of the Schneider-Andruskiewitsch program of classification
of pointed Hopf algebras (see \cite{AS10}).
Let $\trR:=\trRtbhmtPi$.
Heckenberger classified $\trR$ with ${\mathrm{Char}}(\bK)=0$ and $|\trR|<\infty$
(see \cite{Hec09}).
In 2010, Heckenberger and the second author \cite{HY10}
gave a factorization 
formula of the Shapovalov determinants of $\trU$
defined for $\tbhm$ with $|\trR|<\infty$ and $\tbhm(\tal,\tal)\ne 1$ for
all $\tal\in\trR$,
see Theorem~\ref{theorem:Shapo}.
Using it, 
we prove Theorem~\ref{theorem:mainTh} along the same argument as that given by  
Kac~\cite{Kac84}.
As in \cite{HY10}, we also use a `density argument', that is,
we first discuss `at roots of unity' and then pass to the general cases. 
In \cite{AY2015}, Angiono and the second author gave the explicit formula
of the universal ${\mathrm{R}}$-matrix of $\trU$
defined for $\tbhm$ with $|\trR|<\infty$
(${\mathrm{R}}$ is different from $\trR$). 
In \cite{AYY15}, Azam, Yousofzadeh and the second author classified the finite-dimensional
irreducible representations of $\trU$
defined for $\tbhm$ with ${\mathrm{Char}}(\bK)=0$ and $|\trR|<\infty$.
In \cite{BY15}, the authors of this paper studied skew centers for rank-one cases.
In \cite{Y15}, the second author introduced a new definition
of generalized root systems, which is by far easier than that of the original one. 

Let $\trR:=\trRtbhmtPi$. For $\tal\in\trR$, let $\tq_\tal:=\tbhm(\tal,\tal)$.
\begin{center}
Assume that $\trR$ is a finite set. Assume that $\tq_\tal\ne 1$ for
all $\tal\in\trR$.
\end{center}
Let $\chktfkB^{\tbhm,\checktal}_\newpara$
be the $\bK$-subspace of $\chktrUo$ formed by
the elements 
$\sum_{\tlambda\in\tfkAPi}a_\tlambda
\chktrK_\tlambda$ with $a_\tlambda\in\bK$
satisfying the following (${\check{e}}1$)-(${\check{e}}4$).
For $x\in\bKt\setminus\{1\}$, let $\mcD(x):=\{x^n|n\in\bZ\}$,
and let $\kpch(x):=\min\{r\in\bN|x^r=1\}$
(resp $:=0$) if $\mcD(x)$ is a finite
(resp. infinite) set.
Define the $\bZ$-module homomorphism $\trho=\trho^{\tbhm,\checktal}:\tfkAPi\to\bKt$ 
by $\trho(\talpij):=\tq_\talpij$ ($j\in\fkI$).
Let $\newpara^\tbhm_{\tlambda;\tbeta}:={\frac {\newpara(\tbeta)}
{\tbhm(\tlambda,\tbeta)}}$
for $\tlambda\in\tfkAPi$ and $\tbeta\in\trR$.
\newline\par
(${\check{e}}1$) $a_{\tlambda+2t\tbeta}
=\trho(\tbeta)^t\cdot a_\tlambda$ for $\tlambda\in\tfkAPi$, $\tbeta\in\trR$ and $t\in\bZ\setminus\{0\}$
with 
$\kpch(\tq_\tbeta)=0$ 
and $\newpara^\tbhm_{\tlambda;\tbeta}=\tq_\tbeta^t$.
\par
(${\check{e}}2$) $a_\tlambda=0$ for $\tlambda\in\tfkAPi$
and $\tbeta\in\trR$
with 
$\kpch(\tq_\tbeta)=0$ 
and $\newpara^\tbhm_{\tlambda;\tbeta}\notin\mcD(\tq_\tbeta)$.
\par
(${\check{e}}3$) 
$\sum_{x=-\infty}^\pinfty (a_{\tlambda+2(\kpch(\tq_\tbeta)x+t)\tbeta}
-\trho(\tbeta)^ta_{\tlambda+2\kpch(\tq_\tbeta)x\tbeta})
\trho(\tbeta)^{-\kpch(\tq_\tbeta)x}
=0$ 
for $\tlambda\in\tfkAPi$, $\tbeta\in\trRp$
and $t\in\{1,2,\ldots,\kpch(\tq_\tbeta)-1\}$
with 
$\kpch(\tq_\tbeta)\geq 2$ and
$\newpara^\tbhm_{\tlambda;\tbeta}=\tq_\tbeta^t$.
\par
(${\check{e}}4$) 
$\sum_{y=-\infty}^\pinfty (a_{\tlambda+2(\kpch(\tq_\tbeta)y+k)\tbeta}
-\trho(\tbeta)^ka_{\tlambda+2\kpch(\tq_\tbeta)y\tbeta})
\trho(\tbeta)^{-\kpch(\tq_\tbeta)y}
=0$ 
for $\tlambda\in\tfkAPi$, $\tbeta\in\trRp$
and $k\in\{1,2,\ldots,\kpch(\tq_\tbeta)-1\}$
with 
$\kpch(\tq_\tbeta)\geq 2$ and
$\newpara^\tbhm_{\tlambda;\tbeta}\notin\mcD(\tq_\tbeta)$.

\begin{theorem}\label{theorem:SYMMAIN}
We have $\rmIm\chktrHCmap^{\tbhm,\checktal}_\newpara
=\chktfkB^{\tbhm,\checktal}_\newpara$.
{\rm{(}}Note that $\chktrHCmap^{\tbhm,\checktal}_\newpara$ is injective
by Lemmas~{\rm{\ref{lemma:HCinj}}} and {\rm{\ref{lemma:preSymMain}}}.{\rm{)}}
\end{theorem}
The proof will be given in Subsection~\ref{subsection:CmSgc}.
Our main result Theorem~\ref{theorem:mainTh} for $\trU$ is similar to the above theorem
(see Theorem~\ref{theorem:DefofGQG} and \eqref{eqn:assNlcs} for $\trU$),
but we do not assume $\tbhm(\tlambda,\tmu)
=\tbhm(\tmu,\tlambda)$
($\tlambda$, $\tmu\in\tfkAPi$).

This paper is organized as follows.
Section~\ref{section:Intr} is Introduction.
In Section~\ref{section:notation}, we collect necessary notations in the paper.
In Section~\ref{section:App},
we collect results obtained by
well-known argument using
the structure theorem for finitely generated modules over a principal ideal domain.
In Section~\ref{section:AppII}, we treat the `density argument'.
In Section~\ref{section:setGRSandGQG}, we collect necessary facts of the generalized root systems.
In Section~\ref{section:GQGr}, we introduce the generalized quantum groups $\trU=\trUtbhmtPi$
and the Kharchenko's PBW theorem of $\trU$,
whose statement is described using its generalized root system.
In Section~\ref{section:SingVec}, we give explicit formulas of singular vectors of Verma modules of $\trU$
defined `at roots of unity'.
In Section~\ref{section:ShFac}, we give a key theorem Theorem~\ref{theorem:nodgthm}
using the results of Sections~\ref{section:App} and \ref{section:AppII}.
In Section~\ref{section:SkewCenters}, we introduce the equations $(e1)_\tbeta$-$(e4)_\tbeta$ ($\tbeta\in\trRtbhmtPi$),
which formulate our Harich-Chandra-type theorem Theorem~\ref{theorem:mainTh}.
In Section~\ref{section:KKarg}, we complete a proof of our main result Theorem~\ref{theorem:mainTh}.
In Section~\ref{section:symcase}, we assume $\tbhm$ to be symmetric,
treat the quotient $\bK$-algebra $\chktrU=\trU/(\trL_\tlambda-\trK_{-\tlambda}
|\tlambda\in\trtfkAPi)$, and explain a Harish-Chandra-type theorem
of the quantum superalgebra $U_q(\mfks)$ of a
basic classical Lie superalgebra $\mfks$.

As to previous study of centers of quantum groups, we
mention \cite{C1Lach}, \cite{C2Tnage},
and \cite{C3Tanisaki}.

\tableofcontents 

\section{Notation} \label{section:notation}
Let $\krdel_{a,b}$ means the Kronecker's delta, i.e.,
$\krdel_{a,a}:=1$,
and $\krdel_{a,b}:=0$ if $a\ne b$.

For $a$, $b\in\bR$, let $\fkJ_{a,b}:=\{\,n\in\bZ\,|
\,a\leq n\leq b\,\}$,
and $\fkJ_{a,\pinfty}:=\{\,n\in\bZ\,|\,a\leq n\,\}$.
Let $\bZgeqo:=\fkJ_{0,\pinfty}$. Note
$\bN=\fkJ_{1,\pinfty}$.

Let $\bK$ be an algebraically closed field.
Let $\bKt:=\bK\setminus\{0\}$.
For $n\in\bZgeqo$ and $x\in\bK$, let $(n)_x:=\sum_{r=1}^n x^{r-1}$,
and $(n)_x!:=\prod_{r=1}^n(r)_x$.
For $n\in\bZgeqo$, $m\in\fkJ_{0,n}$ and $x\in\bK$,
define ${n\choose m}_x\in\bK$ by
${n\choose 0}_x:={n\choose n}_x:=1$,
and ${n\choose m}_x:={n-1\choose m}_x+x^{n-m}{n-1\choose m-1}_x
=x^m{n-1\choose m}_x+{n-1\choose m-1}_x$
(if $m\in\fkJ_{1,n-1}$).
If $(m)_x!(n-m)_x!\ne 0$, then ${n\choose m}_x=
{\frac {(n)_x!} {(m)_x!(n-m)_x!}}$.
For $x$, $y$, $z\in\bK$, and $n\in\bN$, we have
$\prod_{t=0}^{n-1}(y+x^tz)=\sum_{m=0}^n
x^{{\frac {m(m-1)} 2}}{n\choose m}_xy^{n-m}z^m$.

For $n\in\bZgeqo$, and $x$, $y\in\bK$,
let $(n;x,y):=1-x^{n-1}y$
and $(n;x,y)!:=\prod_{m=1}^n(m;x,y)$.

For $x\in\bKt$, define $\kpch(x)\in\bZgeqo\setminus\{1\}$ by
\begin{equation}\label{eqn:chd}
\kpch(x):= 
\left\{\begin{array}{l}\min\{\,r^\prime\in\fkJ_{2,\pinfty}\,|\,(r^\prime)_x!=0\,\} 
\,\,\mbox{if $(r^{\prime\prime})_x!=0$ for some $r^{\prime\prime}\in\fkJ_{2,\pinfty}$}, \\
0 \quad\mbox{otherwise.}
\end{array}\right.
\end{equation}

For an associative $\bK$-algebra $\fka$ and $X$, $Y\in \fka$,
let $[X,Y]:=XY-YX$. For a set $Z$, let $|Z|$ denote the cardinality of $Z$.

For $m\in\bN$, $x$, $y\in\fkJ_{1,m}$, and a commutative ring $X$ having the unit $1$, 
let $\Mat(m,X)$ be the set of $m\times m$ matrices over $X$, 
let $\rmGL(m,X)$ be the set of invertible $m\times m$ matrices over $X$, and 
let $E^{X,m}_{x,y}$ denote the element of $\Mat(m,X)$
whose $(x,y)$-component is $1$ and whose $(x^\prime,y^\prime)$-components are $0$
for $(x^\prime,y^\prime)\in\fkJ_{1,m}^2\setminus\{(x,y)\}$,
that is, $E^{X,m}_{x,y}$ is the $(x,y)$-matrix unit.
However, as usual, we shall always abbreviate $E^{X,m}_{x,y}$ to $E_{x,y}$; we shall be able to deduce 
$X$ and $m$ from context.
For a non-empty subset $Y$ of $X$, 
let $\Mat(m,Y):=\sum_{x=1}^m\sum_{y=1}^mYE^{X,m}_{x,y}$.

For a $\bK$-algebra $T$ (with $1$),
let $\rmCh(T)$ be the 
set of
$\bK$-algebra 
homomorphisms from $T$ to $\bK$.

For an associative $\bK$-algebra $\fka$, if $\fka$ is an integral domain, 
we let $\quoQ(\fka)$ denote the field of fractions of $\fka$, and see that
$\fka$ is regarded as a $\bK$-subalgebra of $\quoQ(\fka)$ in a natural way.

\section{Removal argument}\label{section:App}

\subsection{Factorization property of matrices over a polynomial algebra}
\label{subsection:sApp}

Recall that for a topological space $Z$, a subset $Y$ of $Z$
is called {\it{locally closed}}
if $Y$ is an intersection of an open subset and a closed subset.
As for coordinate ring of algebraic varieties,
see \cite[Observation~4.23,\,Exercise~6.26]{SR05}, for example. 

Let $k\in\bN$. In Sections~\ref{section:App} and \ref{section:AppII},
we study a locally closed subset of the affine variety 
$\bKtk$,
whose coordinate $\bK$-algebra
$\KbKtk$ can be identified with the Laurent polynomial
algebra $\bK[\,\rmX_i^{\pm 1}\,|\,i\in\fkJ_{1,k}\,]$ in $k$-variables $\rmX_i$
($i\in\fkJ_{1,k}$).
Recall that the coordinate $\bK$-algebra $\KbKk$ of the affine variety $\bKk$
can be identified with the polynomial $\bK$-algebra 
$\bK[\,\rmX_i\,|\,i\in\fkJ_{1,k}\,]$ 
in $k$-variables $\rmX_i$
($i\in\fkJ_{1,k}$). Under the identification, we can regard
that $\bKtk\subset\bKk$, and 
\begin{equation*}
\KbKk=\bK[\,\rmX_i\,|\,i\in\fkJ_{1,k}\,]\subset\KbKtk=\bK[\,\rmX_i^{\pm 1}\,|\,i\in\fkJ_{1,k}\,].
\end{equation*}

In Section~\ref{section:App}, we use the notation as follows.
Fix $f=f(\rmX)\in\KbKtk$
and assume that $f$ is a non-zero non-invertible
irreducible element of $\KbKtk$.
Let  
$\mclF:=\{h\in \KbKtk|\forall p\in\KbKtk, h\ne pf \}$.
Let $\mclF^{-1}\KbKtk:=\{{\frac {p_1} {p_2}}\in\quoQ(\KbKtk)|p_1\in\KbKtk,
p_2\in\mclF\}$, i.e., it is identified with the localization 
of $\KbKtk$ by $\mclF$.
Note the following.
\begin{equation}\label{eqn:elfPID}
\begin{array}{l}
\mbox{(1) $\KbKk$, $\KbKtk$ and $\mclF^{-1}\KbKtk$
are unique factorization domains.} \\
\mbox{(2) $\mclF^{-1}\KbKtk$
is a principal ideal domain whose ideals are} \\
\mbox{
$(\mclF^{-1}\KbKtk)f^t$ ($t\in\bZgeqo$).}
\end{array}
\end{equation}
Fix $g=g(\rmX)\in\mclF$.
For $h=h(\rmX)\in\KbKtk$, let
$\mclV(h):=\{\,\rmx\in\bKtk\,|\,h(\rmx)=0\,\}$,
and $\mclD(h):=\bKtk\setminus\mclV(h)$.

We can easily show the following two lemmas using Hilbert's Nullstellensatz
and  the Zariski topology.

\begin{lemma}\label{lemma:hilzr} 
Let $h\in\KbKtk$ be such that $\mclV(h)\supset\mclV(f)\cap\mclD(g)$.
Then $h\in\KbKtk f$. In particular, $\mclV(h)\supset\mclV(f)$.
\end{lemma}

\begin{lemma}\label{lemma:hilzrd} We have $\mclV(f)\cap\mclD(g)\ne\emptyset$.
\end{lemma}
By (2) of \eqref{eqn:elfPID} and \cite[XV~\S2~Theorem~5]{Lang65}, we have
\begin{lemma}\label{lemma:pidd}
Let $G\in\Mat(m,\KbKtk)$.
Then there exist $m^\prime\in\fkJ_{0,m}$, $c_x\in\bZgeqo$
$(x\in\fkJ_{1,m^\prime})$, 
$h\in\mclF$ and 
$P_t\in\Mat(m,\KbKtk)$ $(t\in\fkJ_{1,2})$
with $\det(P_t)\in\mclF$ such that
\begin{equation}\label{eqn:eqlmd}
P_1GP_2=h\cdot\sum_{x=1}^{m^\prime}f^{c_x}E_{x,x}.
\end{equation}
In particular, if $\det(G)\ne 0$, then
$m^\prime=m$ and  $\det(G)=p\cdot f^{c_1+\cdots+ c_m}$ 
for some $p\in\mclF$.
Moreover, letting $m^{\prime\prime}:=|\{\,x\in\fkJ_{1,m}\,|\,c_x= 0\}|$,
we have $\MatRank(G(\rmx))=m^{\prime\prime}$
{\rm{(}}$\rmx\in\mclV(f)\cap\mclD(p^\prime)${\rm{)}}
for some $p^\prime\in\mclF$.
\end{lemma}

\begin{lemma}\label{lemma:aplmt} 
Keep the notation as in Lemma~{\rm{\ref{lemma:pidd}}}.
Let $r\in\fkJ_{0,m}$.
Let $r\in\fkJ_{0,m}$.
Assume that
\begin{equation}\label{eqn:aplmt}
\MatRank(G(\rmx))\leq m-r \quad (\rmx\in\mclV(f)\cap\mclD(g)).
\end{equation} Then 
\begin{equation}\label{eqn:rtaplmt}
\det G\in\KbKtk f^r.
\end{equation} 
Moreover we have $m^{\prime\prime}\leq m-r$.
\end{lemma}
{\it{Proof.}} 
Use the same notation as in Lemma~\ref{lemma:pidd}
for $G$.
We may assume $m^\prime=m$.
By Lemma~\ref{lemma:hilzrd}, $\emptyset\ne\mclV(f)\cap\mclD(g\cdot\det(P_1P_2)h)
\subset\mclV(f)\cap\mclD(g)$. By \eqref{eqn:eqlmd} and \eqref{eqn:aplmt},
we have $m^{\prime\prime}\leq m-r$.
By \eqref{eqn:eqlmd} and (1) of \eqref{eqn:elfPID},
we have \eqref{eqn:rtaplmt}. \hfill $\Box$
\newline\par
We only use Lemma~\ref{lemma:aplmt} to incidentally give a proof of Theorem~\ref{theorem:Shapo}.
It has originally been proved in \cite{HY10}, and we essentially use it in this paper. After 
Theorem~\ref{theorem:keythma}, we give a short proof of it by using Lemma~\ref{lemma:aplmt}
and other results in this paper although it is essentially the same as the proof given in
\cite{HY10}.

\subsection{Key Lemma to removal of a denominator} \label{subsection:Van} 

The same result as Lemma~\ref{lemma:maintool} below
with $\KbKk$ in place of $\KbKtk$ has been used to study 
the center of the universal enveloping algebra of a basic classical Lie superalgebras,
see \cite[Lemma~2]{Kac84}, \cite[Lemma~13.2.6]{Musson12}.
See also the paragraph after Lemma~\ref{lemma:maintool}.

\begin{lemma}\label{lemma:maintool} Keep the notation as in Lemma~{\rm{\ref{lemma:pidd}}}.
Let $r\in\fkJ_{0,m}$. Let $p\in\mclF$.
Let $Y=Y(\rmX)\in\Mat(m,\KbKtk)$. Assume {\rm{(i)}}, {\rm{(ii)}} and {\rm{(iii)}} below. 
\newline\par
{\rm{(i)}} $\det(G(\rmX))=f(\rmX)^rg(\rmX)$. \par 
{\rm{(ii)}} $\MatRank(G(\rmx))\leq m-r$ $(\rmx\in\mclV(f)\cap\mclD(p))$. \par
{\rm{(iii)}} $\ker G(\rmx)\subset\ker Y(\rmx)$ $(\rmx\in\mclV(f)\cap\mclD(p))$.
\newline\newline
Then there exists $Z=Z(\rmX)\in\Mat(m,\KbKtk)$
such that
\begin{equation}\label{eqn:vn}
Y(\rmX)G(\rmX)^{-1}={\frac 1 {g(\rmX)}}Z(\rmX)
\end{equation} as an equation of $\Mat(m,\quoQ(\KbKtk))$.
Moreover $\MatRank(G(\rmx))= m-r$ $(\rmx\in\mclV(f)\cap\mclD(p^\prime))$
for some $p^\prime\in\mclF\cap\KbKtk p$.
Furthermore, we have $m=m^\prime=m^{\prime\prime}+r$, and assuming $c_x\leq c_{x+1}$ {\rm{(}}$x\in\fkJ_{1,m^\prime-1}${\rm{)}},
we have $c_{x^\prime}=1$ {\rm{(}}$x^\prime\in\fkJ_{m^{\prime\prime}+1,m}${\rm{)}}.
\end{lemma}
{\it{Proof.}}
By (i), since $g\ne 0$, we have
\begin{equation}\label{eqn:mpem}
m^\prime=m\quad\mbox{and}\quad\sum_{x=m^{\prime\prime}+1}^m c_x=r,\,\,
\mbox{which implies}\,\,m-m^{\prime\prime}\leq r.
\end{equation}
Let $p^\prime=p^\prime(\rmX):=p(\rmX)h(\rmX)\det(P_1(\rmX)P_2(\rmX))\in\mclF$.
Note that $\mclV(f)\cap\mclD(p^\prime)\ne\emptyset$ holds by Lemma~\ref{lemma:hilzrd}.
Note that $m^{\prime\prime}=\MatRank(G(\rmx))$
($\rmx\in\mclV(f)\cap\mclD(p^\prime)$).
By (ii),
$m^{\prime\prime}\leq m-r$.
By \eqref{eqn:mpem}, we conclude that $m^{\prime\prime}= m-r$ and $c_x=1$ for all $x\in\fkJ_{m-r+1,m}$. 

Define $T=T(\rmX)\in\Mat(m,\KbKtk)$ by
\begin{equation*}
T(\rmX):=\sum_{x=1}^{m-r} E_{x,x}
+f(\rmX)\sum_{y=m-r+1}^m E_{y,y}.
\end{equation*}
By the above argument, we have
\begin{equation}\label{eqn:mpemd}
P_1(\rmX)G(\rmX)P_2(\rmX)=h(\rmX)T(\rmX).
\end{equation} By (iii) and \eqref{eqn:mpemd}, we see that
for any $\rmx\in\mclV(f)\cap\mclD(p^\prime)$, any $x\in\fkJ_{1,m}$
and any $y\in\fkJ_{m-r+1,m}$, the $(x,y)$-component of 
$P_1(\rmx)Y(\rmx)P_2(\rmx)$ is zero.
Hence by Lemma~\ref{lemma:hilzr}, there exists 
$D(\rmX)\in\Mat(m,\KbKtk)$
such that 
\begin{equation}\label{eqn:mpemdd}
P_1(\rmX)Y(\rmX)P_2(\rmX)=D(\rmX)T(\rmX).
\end{equation} 
By \eqref{eqn:mpemd} and \eqref{eqn:mpemdd},
we have $Y(\rmX)={\frac 1 {h^{\prime\prime}(\rmX)}}D^\prime(\rmX)G(\rmX)$
for some $h^{\prime\prime}(X)\in\mclF$
and some $D^\prime(\rmX)\in\Mat(m,\KbKtk)$, i.e.,
${\frac 1 {h^{\prime\prime}(\rmX)}}D^\prime(\rmX)={\frac 1 {h(\rmX)}}P_1(\rmX)^{-1}D(\rmX)P_1(\rmX)$.
By (i), we see that this lemma holds. \hfill $\Box$
\newline\par
The reason we give Lemma~\ref{lemma:maintool} is as follows. Let $\trU$ be as in Theorem~\ref{theorem:DefofGQG}
below. Assume the same condition as in \eqref{eqn:assNlcs}. 
By \cite[Theorem~7.3]{HY10}, we have already obtained Theorem~\ref{theorem:Shapo},
which is a Shapovalov determinant factorization formula for $\trU$.
In Theorem~\ref{theorem:ndgthCor}, we see an explicit property
of singular vectors of Verma module of $\trU$.
Let ${\widehat{\trU}}$ be a completion of $\trU$ defined as in \cite[(1.3.3)]{C3Tanisaki}.
In the proof of Lemma~\ref{lemma:EXTlem}, using Lemma~\ref{lemma:maintool} we develop an argument
similar to a Kac's well-known argument of \cite{Kac84},
see also \cite[The proof of Theorem~13.1.1]{Musson12}. 
By Theorems \ref{theorem:Shapo}, \ref{theorem:ndgthCor} and Lemma~\ref{lemma:EXTlem},
we get an element of the (skew) center of ${\widehat{\trU}}$. 
By Lemma~\ref{lemma:estlm}, we show that it can be expressed as a finite sum, so
it belongs to the (skew) center $\prtrZnewparabctal$ of $\trU$;
in fact it is $V$ of Lemma~\ref{lemma:premainTh}. Thus we see that the Harish-Chandra map 
$\trHCmapbctaln:\prtrZnewparabctal\to\tfkBbctanp$ is surjective, see our main result Theorem~\ref{theorem:mainTh}.
(See \eqref{eqn:defprtrZomega}, \eqref{eqn:deftfkBbctanp}, \eqref{eqn:deftrHComega}
for the definitions of $\prtrZnewparabctal$, $\tfkBbctanp$, $\trHCmapbctaln$ respectively.
See also \eqref{eqn:spsibsc-1}.)

\section{Density argument}\label{section:AppII}

The main results in this Section have already been obtained 
in \cite[Section~9]{HY10}. 
In \cite[Section~9]{HY10}, those
were written in terms of coordinate rings of algebraic varieties.
We reformulate them in more elemental way in order to make the reader efficiently 
understand delicate arguments
in this paper.

\subsection{Easy lemma}

Recall that for a topological space $X$,
a subset $D$ of $X$ is said to be {\it{dense}}
if
\begin{equation}\label{eqn:FXDF}
\cap_{F:closed\,in X,D\subset F}F=X.
\end{equation}

We need the following easy lemma,
which can be proved in a standard way.

\begin{lemma} \label{lemma:basica}
Let $X$ be a topological space.
Let $D$ be a dense subset of $X$.
Let $V$ be an open subset of $X$.
Then
$V\subset\cap_{F:closed\,in X,D\cap V\subset F}F$,
that is, $D\cap V$ is a dense subset of $V$
under the relative topology.
\end{lemma}

\subsection{Nichols topology}

Here we express elements of $\bK^k$ for some $k\in\bN$ by column vectors.

Let $A=[a_{ij}]_{1\leq i,j\leq k}\in\rmGL(k,\bZ)$.
Define the map
$f_A:\bKtk\to\bKtk$ by 
\begin{equation*}
f_A(\left[\begin{array}{c} z_1 \\ \vdots \\ z_k
\end{array}\right] ):=\left[\begin{array}{c} z_1^{a_{11}}\cdots z_k^{a_{1k}}
\\ \vdots \\ z_1^{a_{k1}}\cdots z_k^{a_{kk}}
\end{array}\right].
\end{equation*} Then $f_A$ is bijective, and 
\begin{equation}\label{eqn:prefAmV}
f_A^{-1}=f_{A^{-1}},
\end{equation} since $f_{A_1}f_{A_2}=f_{A_1A_2}$ ($A_1$, $A_2\in\bK^k$).

Let
\begin{equation}\label{eqn:dfKprime}
\bKtprim:=\{\,x\in\bKt\,|\,\exists r\in\bN,\,x^r=1\,\}.
\end{equation}

For $k\in \bN$, $d\in\bKtprim$ and $m_x\in\bZ$ ($x\in\fkJ_{1,k}$),
let 
\begin{equation}\label{eqn:dfprcrs}
V(d;m_1,\dots ,m_k):=\{\,\left[\begin{array}{c} z_1 \\ \vdots \\ z_k
\end{array}\right] \in\bKtk\,|\,z_1^{m_1}\cdots z_k^{m_k}=d\,\},
\end{equation} call it
a {\it{$1$-pre-fundamental Nichols closed subset}} of $\bKtk$,
and call its complement 
an {\it{$1$-pre-fundamental Nichols open subset}} of $\bKtk$.
In addition to the Zariski topology on $\bKtk$, we also treat the
topology on $\bKtk$ such that the family of all $1$-pre-fundamental Nichols 
open subsets forms its open subbase; we
call it the {\it{Nichols topology}} on $\bKtk$.
We call an open (resp. closed) subset in the sense of the Nichols topology
a {\it{Nichols open {\rm{(}}resp. closed{\rm{)}} subset}}.  

As for \eqref{eqn:dfprcrs},
\begin{equation}\label{eqn:TrivNch}
\begin{array}{l}
\mbox{$V(d;m_1,\dots ,m_k)$ is $\bKtk$
(resp. $\emptyset$)
if and only if} \\
\mbox{$m_x=0$ ($x\in\fkJ_{1,k}$)
and $d=1$ (resp. $d\ne 1$).}
\end{array}
\end{equation}

For $A\in\rmGL(k,\bZ)$, and 
$^t[m_1,\ldots,m_k]$, $^t[{\check m}_1,\ldots,{\check m}_k]\in\bZ^k$ with
\begin{equation*}
[{\check m}_1,\ldots,{\check m}_k]=[m_1,\ldots,m_k]A,
\end{equation*} by \eqref{eqn:prefAmV}, we have
\begin{equation}\label{eqn:fAmV}
(f_A)^{-1}(V(d;m_1,\dots ,m_k))=V(d;{\check m}_1,\dots ,{\check m}_k).
\end{equation}

Let $r\in\fkJ_{1,k}$, $d_x\in\bKtprim$ ($x\in\fkJ_{1,k}$),
$M=[m_{ij}]_{1\leq i,j\leq k}\in\rmGL(k,\bZ)$, and
\begin{equation*}
V^\prime:=\cap_{x=1}^r V(d_x;m_{x,1},\ldots,m_{x,k}).
\end{equation*}
By \eqref{eqn:fAmV}, we have
\begin{equation}\label{eqn:rprecrs}
V^\prime= f_M^{-1}(\{\,
\left[\begin{array}{c} 
d_1 \\ \vdots \\ d_r \\ z_{r+1} \\ \vdots \\ z_k
\end{array}\right]\in\bKtk\,|\,
\left[\begin{array}{c} 
z_{r+1} \\ \vdots \\ z_k
\end{array}\right]\in(\bKt)^{k-r}\,\}).
\end{equation} 
We call such a Nichols closed subset of $\bKtk$ as $V^\prime$
an {\it{$r$-fundamental Nichols closed subset}} of $\bKtk$,
or a {\it{fundamental Nichols closed subset}} for simplicity.

\begin{lemma} \label{lemma:onecrs} 
Any non-empty proper $1$-pre-fundamental Nichols
closed subset of $\bKtk$ is a finite union of
some $1$-fundamental Nichols
closed subsets.
\end{lemma}
{\it{Proof.}} Let $V$ be as in \eqref{eqn:dfprcrs}, and
assume that $\emptyset\subsetneq V\subsetneq\bKtk$.
By \eqref{eqn:TrivNch}, we have $[m_1,\ldots,m_k]\ne [0,\ldots,0]$.
Let $m^\prime:=\rmgcd\{\,m_x\,|\,x\in\fkJ_{1,k}\,\}$.
Let $M^{\prime\prime}=[m^{\prime\prime}_{xy}]_{1\leq x,y\leq k}\in\rmGL(k,\bZ)$
be such that $m^{\prime\prime}_{1y}={\frac {m_y} {m^\prime}}$
($y\in\fkJ_{1,k}$).
Let $X:=\{\,z\in\bKt\,|\,z^{m^\prime}=1\,\}$.
Then $X$ is a finite set.
By \eqref{eqn:fAmV}, we have
\begin{equation*}
\begin{array}{l}
V=V(d;m_1,\ldots,m_k)=f_{M^{\prime\prime}}^{-1}(V(d;m^\prime,0,\ldots,0)) \\
\quad = \cup_{z\in X}f_{M^{\prime\prime}}^{-1}(V(z;1,0,\ldots,0))
=\cup_{z\in X}V(z;m^{\prime\prime}_{11},\ldots,m^{\prime\prime}_{1k}),
\end{array}
\end{equation*} as desired. This completes the proof. \hfill $\Box$

\begin{lemma}\label{lemma:fncvr}
The complement of any finite union of proper pre-$1$-fundamental Nichols closed subsets
of $\bKtk$
is a dense subset of $\bKtk$ under the Zariski topology.
\end{lemma}
{\it{Proof.}} We use the notation as in
Subsection~\ref{subsection:sApp}.
Let $g:=\prod_{y=1}^p(d_y
-\prod_{x=1}^k\rmX_x^{m_{y,x}})\in\KbKtk$
for some $p\in\bN$, some $d_y\in\bKt$
($y\in\fkJ_{1,p}$)  and some $m_{y,x}\in\bZ$ 
($x\in\fkJ_{1,k}$) with
$[m_{y,1},\ldots,m_{y,x}]\ne [0,\ldots,0]$,
see also \eqref{eqn:TrivNch}.
Let $h\in\KbKtk$ be such that $\mclD(g)\subset\mclV(h)$,
so  $\mclV(g)\cup\mclV(h)=\mclV(gh)=\bKtk$.
By an argument similar to that for the proof of Lemma~\ref{lemma:hilzr},
we see that there exists $t_p\in\bN$ ($p\in\fkJ_{1,k}$) such that
$gh\prod_{p=1}^k X_p^{t_p}=0$
(as an element of $\KbKk$).
Hence $h=0$.
This completes the proof. \hfill $\Box$
\newline\par
If a subset of $\bKtk$ is a finite intersection of some
$1$-pre-fundamental Nichols closed subsets, we call it a
{\it{pre-fundamental Nichols closed subset}}.

\begin{lemma} \label{lemma:intcrs} 
Let $r\in\fkJ_{1,k-1}$, and let $V$ be an $r$-fundamental Nichols closed subset of $\bKtk$.
If $V^\prime$ is a $1$-pre-fundamental Nichols closed subset of $\bKtk$
with $\emptyset\subsetneq V^\prime\cap V\subsetneq V$,
then $V^\prime\cap V$ is a finite union of some
$(r+1)$-fundamental Nichols closed subsets. Moreover,
if $V^{\prime\prime}$ is a pre-fundamental Nichols closed subset of $\bKtk$
with $\emptyset\subsetneq V^{\prime\prime}\cap V\subsetneq V$,
then $V^{\prime\prime}\cap V=\cup_{y=1}^xV^{\prime\prime\prime}_y$
for some $x\in\bN$ and some $r_y$-fundamental Nichols closed subsets
$V^{\prime\prime\prime}_y$ $(y\in\fkJ_{1,x})$ with
$r_y\in\fkJ_{r+1,k}$.
\end{lemma}
{\it{Proof.}} This follows easily from Lemma~\ref{lemma:onecrs}
and \eqref{eqn:rprecrs}. \hfill $\Box$

\begin{lemma} {\rm{(}}See also {\rm{\cite[Lemma~9.3]{HY10}}}.{\rm{)}} \label{lemma:maincrs} 
Any non-empty proper Nichols closed subset of $\bKtk$ is a finite union of some
fundamental Nichols closed subsets.
\end{lemma}
{\it{Proof.}} By Lemma~\ref{lemma:onecrs}, it suffices to prove the fact $(\sharp)$ below.
\newline\par
$(\sharp)$ Let $X$ be an non-empty set, and let $f:X\to\bN$ be a map.
Assume that for any $x\in X$ and $y\in\fkJ_{1,f(x)}$, a fundamental Nichols closed subset
$V_{x,y}$ is given. Then 
\begin{equation}\label{eqn:cpcuVxy}
\bigcap_{x\in X}(\bigcup_{y=1}^{f(x)}V_{x,y})
\end{equation} is a finite union of some fundamental closed Nichols
subsets.
\newline\par 
We prove $(\sharp)$. Let $V$ be the Nichols closed set in \eqref{eqn:cpcuVxy}.
We may assume that $|X|\geq 2$.
For $x\in X$ and $y\in\fkJ_{1,f(x)}$,
let $r_{x,y}\in\fkJ_{1,k}$ be such that
$V_{x,y}$ in \eqref{eqn:cpcuVxy} is an $r_{x,y}$-fundamental Nichols closed subset. 
We use induction on $m:=\rmMax\{\,k-r_{x,y}\,|\,x\in X,\,y\in\fkJ_{1,f(x)}\,\}$.
If $m=0$, $(\sharp)$ is clear from \eqref{eqn:rprecrs}.
Assume $m\geq 1$.
Let $x^\prime\in X$ and $y^\prime\in\fkJ_{1,f(x^\prime)}$
be such that $m=k-r_{x^\prime,y^\prime}$. 
Note that
\begin{equation*}
V=\bigcup_{y^{\prime\prime}=1}^{f(x^\prime)}
(V_{x^\prime,y^{\prime\prime}}\cap(\bigcap_{x\in X\setminus\{x^\prime\}}
(\bigcup_{y=1}^{f(x)}V_{x,y}))).
\end{equation*} 
Hence we may assume that $f(x^\prime)=1$ and $y^\prime=1$.
Moreover we may assume that $\emptyset\subsetneq V_{x^\prime,1}\cap(
\cup_{y=1}^{f(x)}V_{x,y})\subsetneq V_{x^\prime,1}$
for all $x\in X\setminus\{x^\prime\}$.
Furthermore we may assume $\emptyset\subsetneq V_{x^\prime,1}\cap V_{x,y}\subsetneq V_{x^\prime,1}$
for all $x\in X\setminus\{x^\prime\}$
and all $y\in\fkJ_{1,f(x)}$.
Note that
\begin{equation*}
V=\bigcap_{x\in X\setminus\{x^\prime\}}
(\bigcup_{y=1}^{f(x)}(V_{x^\prime,1}\cap V_{x,y})).
\end{equation*}
Hence, using Lemma~\ref{lemma:intcrs} and the induction, we see that $(\sharp)$ holds.
This completes the proof.
\hfill $\Box$

\subsection{Density for Nichols closed subsets}

Let $\bKkgqi:=(\bKtprim)^k$.
Note that
\begin{equation}\label{eqn:prsX}
f_A(\bKkgqi)=\bKkgqi\quad\quad(A\in\rmGL(k,\bZ)).
\end{equation}

By \eqref{eqn:rprecrs}, \eqref{eqn:prsX} and Lemmas~\ref{lemma:basica} and \ref{lemma:maincrs},
we have
\begin{theorem}
{\rm{(}}See also {\rm{\cite[Lemma~9.4]{HY10}}}.{\rm{)}} \label{theorem:basicc} 
Let $V$ be a Nichols closed subset of $\bKtk$.
Let $X$ be a Zariski open subset of $\bKtk$.
Then $\bKkgqi\cap X\cap V$ is
a dense subset of $X\cap V$
under the {\rm{(}}relative{\rm{)}} Zariski topology on $\bKtk$.
\end{theorem}

\section{Generalized root systems}\label{section:setGRSandGQG}

As a definition of the generalized root system,
we adopt the new and simpler one introduced by \cite{Y15}
in order to make this paper more readable.

\subsection{Set-theoretical generalized root system
(SGRS)}\label{subsection:setGRS}
Let $\CardfkIN\in\bN$.
Let $\newmcV$ be a $\CardfkIN$-dimensional $\bR$-linear space.
($\newmcV$ in 
Section~\ref{section:setGRSandGQG} does not relate to $\mclV(\cdot)$ in 
Sections~\ref{section:App} and \ref{section:AppII}.) 
Let $\fkI$ be a set with $\CardfkI=\CardfkIN$.
Let $\tR$ be a subset of $\newmcV$
such that $\rmSpan_\bR(\tR)=\newmcV$
and $0\notin\tR$.
For a non-empty subset $\tPi$ of $\tR$,
let $\tRptPi:=\tR\cap\rmSpan_{\bZgeqo}(\tPi)$ and
$\tRmtPi:=\tR\cap\rmSpan_{\bZleqo}(\tPi)$, and call
$\tPi$ {\it{a base of $\tR$}}
if $|\tPi|=\CardfkIN$ and $\tR=\tRptPi\cup\tRmtPi$.
Let $\tiltB$ be the set of all bases of $\tR$.
Let $\tB$ be a nonempty subset of $\tiltB$.
We call the pair $(\tR,\tB)$ {\it{a set-theoretical generalized root system}}
({\it{an SGRS}} for short) if 
\begin{equation}\label{eqn:defstGRS}
\forall\tPi\in\tB, \forall\tal\in\tPi, 
\tR\cap\bR\tal=\{\tal,-\tal\},
\exists\tPi^{(\tal)}\in\tB, 
\tRp_{\tPi^{(\tal)}}\cap\tRmtPi=\{-\tal\},
\end{equation} where see also \cite{Y15}.

From now on until end of Subsection~\ref{subsection:setGRS},
we fix an SGRS $(\tR,\tB)$.
Note that $(\tR,-\tB)$ are also an SGRS, where
$-\tB:=\{-\tPi|\tPi\in\tB\}$
(where $-\tPi:=\{-\tal|\tal\in\tPi\}$).
Note that $(\tR\cap(-\tR),\tB)$ and $(\tR\cup(-\tR),\tB)$
are also SGRSs.
Let $\tPi\in\tB$ and 
$\tal\in\tPi$. Let $\tN^\tPi_{\tal,\tal}:=-2$ $(\in\bZ)$.
We see that there exist $\tN^\tPi_{\tal,\tbeta}\in\bZgeqo$
($\tbeta\in\tPi\setminus\{\tal\}$)
such that $\tPi^{(\tal)}=\{\tgamma+\tN^\tPi_{\tal,\tgamma}\tal|\tgamma\in\tPi\}$.
Since $(\tPi^{(\tal)})^{(-\tal)}=\tPi$, we have
\begin{equation}\label{eqn:invtN}
\tN^{\tPi^{(\tal)}}_{-\tal,\tbeta+\tN^\tPi_{\tal,\tbeta}\tal}=\tN^\tPi_{\tal,\tbeta}
\quad(\tbeta\in\tPi).
\end{equation} 
Let $\tB[\tPi]:=\{\tPi^\prime\in\tB||\tRm_\tPi\cap\tRp_{\tPi^\prime}|<\infty\}$.
Then for $\tPi^\prime\in\tB$, we see that
$\tPi^\prime\in\tB[\tPi]$ if and only if there exist $k\in\bZgeqo$, $\tPi_t\in\tB$, $\tbeta_t\in\tPi_t$
($t\in\fkJ_{0,k-1})$ such that $\tPi_0=\tPi$, $\tPi_k=\tPi^\prime$
and $\tPi_t=\tPi_{t-1}^{(\tbeta_{t-1})}$ ($t\in\fkJ_{1,k}$).
Note that $(\tR,\tB[\tPi])$ is also an SGRS.
We say that $(\tR,\tB)$ is {\it{connected}} if
$\tB[\tPi]=\tB$ for some $\tPi\in\tB$.
\begin{lemma} \label{lemma:BtB}
{\rm{(See also \cite[Lemmas~2.2(2) and 2.3(1)]{Y15}.)}}
Assume $|\tR|<\infty$.
Then $\tiltB=\tB$, $\cup_{\tPi\in\tB}\tPi=\tR$,
and  $(\tR,\tB)$ is connected.
In particular $\tRmtPi=-\tRptPi$
for $\tPi\in\tB$. so $|\tR|=2|\tRp|$.
\end{lemma}
{\it{Proof.}} The claim easily follows 
from the definition of $\tB[\tPi]$.
\hfill $\Box$ 
\newline\par
Let $\checktB$ be the set of all maps $\checktal:\fkI\to\newmcV$
with $\checktal(\fkI)\in\tB$.
For $i\in\fkI$, define the map $\checktau_i:\checktB\to\checktB$ by
$\checktau_i(\checktal)(j):=\talpij
+\tN^{\checktal(\fkI)}_{\talpii,\talpij}\talpii$
($j\in\fkI$). By \eqref{eqn:invtN}
\begin{equation}\label{eqn:checkinvtN}
\tN^{\checktau_i(\checktal)(\fkI)}_{\checktau_i(\checktal)(i),
\checktau_i(\checktal)(j)}=\tN^{\checktal(\fkI)}_{\talpii,\talpij}
\quad(\checktal\in\checktB,\,i,j\in\fkI).
\end{equation}
Hence $\checktau_i^2=\rmid_\checktB$.
Let $\checkV$ be an $\CardfkIN$-dimensional $\bR$-linear space.
Let $\{\checkv_i|i\in\fkI\}$ be an $\bR$-basis of $\checkV$.
For $\checktal\in\checktB$ and $i\in\fkI$, define
$\tils^\checktal_i\in\rmGL(\checkV)$ by
$\tils^\checktal_i(\checkv_j):=\checkv_j+\tN^{\checktal(\fkI)}_{\talpii,\talpij}\checkv_i$.
Then $\tils^\checktal_i=(\tils^\checktal_i)^{-1}=\tils^{\checktau_i(\checktal)}_i$.
Let $\bbM$ be the set of all maps from $\bN$ to $\fkI$.
For $\checktal\in\checktB$, $f\in\bbM$ and $t\in\bZgeqo$, define 
$\checktal_{f,t}\in\checktB$ by
\begin{equation}\label{eqn:checktalft}
\checktal_{f,t}:=
\left\{
\begin{array}{ll}
\checktal & \quad\mbox{if $t=0$}, \\ 
\checktau_{f(t)}(\checktal_{f,t-1}) & \quad\mbox{if $t\in\bN$}, \\
\end{array}\right.
\end{equation}
and define $1^\checktal\tils_{f,t}\in\rmGL(\checkV)$ by
\begin{equation*}
1^\checktal\tils_{f,t}:=
\left\{
\begin{array}{ll}
\rmid_\checkV & \quad\mbox{if $t=0$}, \\ 
1^\checktal\tils_{f,t-t}\circ 
\tils^{\checktal_{f,t}}_{f(t)} & \quad\mbox{if $t\in\bN$}. \\
\end{array}\right.
\end{equation*} 
For $\checktal\in\checktB$,
let $\checktB[\checktal]:=\{\checktal_{f,t}
|f\in\bbM, t\in\bZgeqo\}$
and
$\bbH[\checktal]:=
\{1^\checktal\tils_{f,t}\in\rmGL(\checkV)|f\in\bbM, t\in\bZgeqo\}$.
Note that $\{\checktal^\prime(\fkI)|\checktal^\prime\in\checktB[\checktal]\}
=\tB[\checktal(\fkI)]$.
Note that $\checktB[\checktal]=\checktB[\checktal^\prime]$
for $\checktal^\prime\in\checktB[\checktal]$.

\begin{lemma} 
{\rm{(See also \cite[Lemma~2.4]{Y15}.)}} 
\label{lemma:basech}
Let $\checktal\in\checktB$, $f\in\bbM$, $t\in\bZgeqo$, 
$\checktal^\prime:=\checktal_{f,t}$,
and $w:=1^\checktal\tils_{f,t}\in\bbH[\checktal]$.
Let $d_{ij}\in\bZ$ $(i,j\in\fkI)$ be such that 
$w(\checkv_j)=\sum_{i\in\fkI}d_{ij}\checkv_i$.
Then for all $j\in\fkI$,
we have $\checktal^\prime(j)=\sum_{i\in\fkI}d_{ij}\talpii$.
In particular, for $f^\prime\in\bbM$ and $t^\prime\in\bZgeqo$,
if $1^\checktal\tils_{f^\prime,t^\prime}=w$ {\rm{(}}that is, $1^\checktal\tils_{f^\prime,t^\prime}$
and $w$ are the same element of $\rmGL(\checkV)${\rm{)}}, 
we have $\checktal^\prime=\checktal_{f^\prime,t^\prime}$.
\end{lemma}
{\it{Proof.}} We use induction on $t$.
It $t=0$, the claim is clear.
Assume $t\geq 1$.
Let $\checktal^{\prime\prime}:=\checktal_{f,t-1}$
and $w^\prime:=1^\checktal\tils_{f,t-1}$.
Let $d^\prime_{ij}\in\bZ$ $(i,j\in\fkI)$ be such that 
$w^\prime(\checkv_j)=\sum_{i\in\fkI}d^\prime_{ij}\checkv_i$.
Let $j\in\fkI$.
Then $\checktal^{\prime\prime}(j)=\sum_{i\in\fkI}d^\prime_{ij}\talpii$.
Let $k:=\tN^{\checktal^\prime(\fkI)}_{\checktal^\prime(f(t)),\checktal^\prime(j)}$.
By \eqref{eqn:checkinvtN}, $k=\tN^{\checktal^{\prime\prime}(\fkI)}_{\checktal^{\prime\prime}(f(t)),\checktal^{\prime\prime}(j)}$.
We have  $w(\checkv_j)=w^\prime\tils^{\checktal^\prime}_{f(t)}(\checkv_j)
=w^\prime(\checkv_j+k\checkv_{f(t)})
=\sum_{i\in\fkI}(d^\prime_{ij}
+kd^\prime_{if(t)})\checkv_i$,
and
$\checktal^\prime(j)=\checktau_{f(t)}(\checktal^{\prime\prime})(j)
=\checktal^{\prime\prime}(j)
+k
\checktal^{\prime\prime}(f(t))
=\sum_{i\in\fkI}(d^\prime_{ij}
+kd^\prime_{if(t)})\talpii$.
Hence the claim holds. 
\hfill $\Box$

\subsection{Categorical generalized root system (CGRS)}\label{subsection:CatGRS}
Let $\CardfkIN$ and $\fkI$ be as in Subsection~\ref{subsection:setGRS}.
We say that
$\tC=[\tc_{ij}]_{i,j\in\fkI}\in\Mat(\CardfkIN,\bZ)$
is a {\it{generalized Cartan matrix}}
if the conditions ${\rm{(M1)}}$, ${\rm{(M2)}}$ below
are satisfied.
\newline\par
${\rm{(M1)}}$ $\tc_{ii}=2$ ($i\in I$).  \par 
${\rm{(M2)}}$ $\tc_{jk}\leq 0$, $\delta_{\tc_{jk},0}=\delta_{\tc_{kj},0}$ 
($j$, $k\in I$, $j\ne k$).
\newline\par
Let $\mcA$ be a non-empty set.
Assume that for each $i\in\fkI$,
a map $\tvsigma_i:\mcA\to\mcA$ is given.
Assume that for each $a\in\mcA$, and each $i$, $j\in\fkI$,
a generalized Cartan matrix
$\tC^a=[\tc^a_{ij}]_{i,j\in\fkI}$ is given.
We say that \begin{equation*}
\mcC=\mcC(\fkI,\mcA,(\tvsigma_i)_{i\in\fkI},(\tC^a)_{a\in\mcA})
\end{equation*} is a {\it{Cartan scheme}} if
the conditions ${\rm{(C1)}}$, ${\rm{(C2)}}$ below
are satisfied.
\newline\par
${\rm{(C1)}}$ \quad $\tvsigma_i^2=\rmid_\mcA$ ($i\in\fkI$). \par
${\rm{(C2)}}$ \quad $\tc^{\tvsigma_i(a)}_{ij}=\tc^a_{ij}$ ($a\in\mcA$, $i,j\in\fkI$).
\newline\par
Let 
$\mcC=\mcC(\fkI,\mcA,(\tvsigma_i)_{i\in\fkI},(\tC^a)_{a\in\mcA})$
be a Cartan scheme. Assume that for each $a\in\mcA$,
there exist
a $\CardfkIN$-dimensional $\bR$-linear space
$\newmcV^a$ and a map
$\checktal^a:\fkI\to\newmcV^a$ such that 
$\checktal^a(\fkI)$ is an $\bR$-basis of $\newmcV^a$.
Define the $\bR$-linear isomorphism
$
\tils^a_i:\newmcV^a\to\newmcV^{\tvsigma_i(a)}\quad(a\in\mcA,\,i\in\fkI)
$
by
$\tils^a_i(\checktal^a(j)):=\checktal^{\tvsigma_i(a)}(j)-
\tc^a_{ij}\checktal^{\tvsigma_i(a)}(i)$ ($j\in\fkI$),
so $\tils^{\tvsigma_i(a)}_i\tils^a_i=\rmid_{\newmcV^a}$ by ${\rm{(C2)}}$.
Let $\bbM$ be as in Subsection~\ref{subsection:setGRS}.
For $a\in\mcA$ and $f\in\bbM$, let
$a_{f,0}:=a$ and $1^a\tils_{f,0}:=\rmid_{\newmcV^a}$
and let $a_{f,t}:=\tvsigma_{f(t)}(a_{f,t-1})$
and $1^a\tils_{f,t}:=1^a\tils_{f,t-t}\circ 
\tils^{a_{f,t}}_{f(t)}$
for $t\in\bN$.
Assume that for each $a\in\mcA$, a subset
$\tR^a$ of $\newmcV^a$ is given.
Let $\tR^{a,+}:=\tR^a\cap(\oplus_{i\in\fkI}\bZgeqo\checktal^a(i))$
and $\tR^{a,-}:=\tR^a\cap(\oplus_{i\in\fkI}\bZleqo\checktal^a(i))$
($a\in\mcA$). We say that the data 
\begin{equation} \label{eqn:cnnRsystem}
\mcR=\mcR(\mcC,(\newmcV^a,\checktal^a,\tR^a)_{a\in\mcA})
\end{equation} is 
{\it{a categorical generalized root system of type $\mcC$}}
({\it{a CGRS of type $\mcC$}} for short)
if the following conditions
${\rm{(R1)}}$-${\rm{(R4)}}$ are satisfied.
\newline\par
${\rm{(R1)}}$ \quad $\tR^a=\tR^{a,+}\cup\tR^{a,-}$ \quad ($a\in\mcA$).\par
${\rm{(R2)}}$ \quad $\tR^a\cap\bZ\checktal^a(i)=\{\checktal^a(i),-\checktal^a(i)\}$ \quad 
($a\in\mcA$, $i\in\fkI$).\par
${\rm{(R3)}}$ \quad $\tils^a_i(\tR^a)=\tR^{\tvsigma_i(a)}$ \quad 
($a\in\mcA$, $i\in\fkI$). \par
${\rm{(R4)}}$ \quad
For $a\in\mcA$, $f\in\bbM$ and $t\in\bZgeqo$, 
if $1^a\tils_{f,t}(\checktal^{a_{f,t}}(i))=\checktal^a(i)$ for all $i\in\fkI$,
then $a_{f,t}=a$.
\newline\newline
(See \cite[Definition~1.2]{AYY15}, \cite{CH09}.)
\newline\par
From now on until end of Subsection~\ref{subsection:CatGRS},
we fix a Cartan scheme $\mcC=\mcC(\fkI,\mcA,(\tvsigma_i)_{i\in\fkI},(\tC^a)_{a\in\mcA})$
and a CGRS $\mcR=\mcR(\mcC,(\newmcV^a,\checktal^a,\tR^a)_{a\in\mcA})$ of type $\mcC$.
It is clear that
\begin{equation}\label{eqn:minusCGRS}
\begin{array}{l}
\mbox{$\mcR(\mcC,(\newmcV^a,\checktal^a,\tR^a\cup(-\tR^a))_{a\in\mcA})$, 
$\mcR(\mcC,(\newmcV^a,\checktal^a,\tR^a\cap(-\tR^a))_{a\in\mcA})$ and} \\
\mbox{$\mcR(\mcC,(\newmcV^a,-\checktal^a,\tR^a)_{a\in\mcA})$ are also CGRSs of type $\mcC$.}
\end{array}
\end{equation}
Note that
\begin{equation}\label{eqn:minuscaij}
-c^a_{ij}=\max\{k\in\bZgeqo|\checktal^a(j)+k\checktal^a(i)\in\tR^a\}
\quad(a\in\mcA,\,i,j\in\fkI,i\ne j).
\end{equation}
Let $a\in\mcA$.
Let $\mcA[a]:=\{a_{f,t}
|f\in\bbM, t\in\bZgeqo\}$,
$\mcC[a]:=\mcC(\fkI,\mcA[a],(\tvsigma_i)_{i\in\fkI},(\tC^a)_{a\in\mcA[a]})$, and 
$\mcR[a]:=\mcR(\mcC[a],(\newmcV^{a^\prime},\checktal^{a^\prime},\tR^{a^\prime})_{a^\prime\in\mcA[a]})$.
Then $\mcC[a]$ is a Cartan scheme, and
\begin{equation} \label{eqn:mcCamcRa}
\mbox{$\mcR[a]$ is a CGRS of type $\mcC[a]$.}
\end{equation} If $\mcA[a]=\mcA$, then we see that $\mcC[a]=\mcC$ 
and $\mcR[a]=\mcR$ and
we say that $\mcC$ is {\it{connected}}
and that $\mcR$ is {\it{connected}}.

By \cite[Lemma~1.5]{AYY15} and \eqref{eqn:minusCGRS}, we have
\begin{lemma}
\label{lemma:fij}
Let $i$, $j\in\fkI$ be such that $i\ne j$.
Let $m:=|\tR^{a,+}\cap(\bZgeqo\checktal^a(i)\oplus \bZgeqo\checktal^a(j))|$,
and assume $m<\infty$.
Let $f\in\bbM$ be such that $f(2x-1)=i$,
$f(x)=j$ $(x\in\bN)$.
Then $a_{f,2m}=a$ and $1^a\tils_{f,2m}=\rmid_{\newmcV^a}$.
\end{lemma}
 
Let $a\in\mcA$.
Let $b$, $b^\prime\in\mcA[a]$.
Let $\mcH(b,b^\prime):=\{
1^b\tils_{f,t}:\newmcV^{b^\prime}\to\newmcV^b|f\in\bbM,t\in\bZgeqo,b_{f,t}=b^\prime\}$.
For $w\in\mcH(b,b^\prime)$, let 
\begin{equation*}
\trell_{b,b^\prime}(w)
:=\min\{t\in\bZgeqo|\exists f\in\bbM,b_{f,t}=b^\prime,1^b\tils_{f,t}=w\}.
\end{equation*} 
The following lemma is well-known. 
\begin{lemma}\label{lemma:prptell}
Let $b$, $b^\prime\in\mcA[a]$\ and $w\in\mcH(b,b^\prime)$. Then we have the followings.
\par
{\rm{(1)}} We have $\trell_{b,b^\prime}(w)=|\tR^{b,-}\cap w(\tR^{b^\prime,+})|$
and $\tR^{b,+}\cap w(\tR^{b^\prime,-})=-(\tR^{b,-}\cap w(\tR^{b^\prime,+}))$.
In particular, for $w^\prime\in\mcH(a,b)$
and $w^{\prime\prime}\in\mcH(a,b^\prime)$,
if $w^\prime(\checktal^b(\fkI))=w^{\prime\prime}(\checktal^{b^\prime}(\fkI))$
as a subset of $\newmcV^a$, then $b=b^\prime$
and $w^\prime=w^{\prime\prime}$. \par
{\rm{(2)}} Let $l\in\bN$. 
Let $f\in\bbM$ be such that $1^b\tils_{f,l}=w$. Then
$l=\trell_{b,b^\prime}(w)$
if and only if
$\tR^{b,+}\cap w(\tR^{b^\prime,-})=\{1^b\tils_{f,t-1}(\checktal^{b_{f,t-1}}(f(t)))|t\in\fkJ_{1,l}\}$.
\end{lemma}
{\it{Proof.}}
(1) This follows from \cite[Lemma~8(iii) (see also Proposition~1)]{HY08},
Lemma~\ref{lemma:fij} and \eqref{eqn:minusCGRS}. \par
(2) Let $l^\prime:=\trell_{b,b^\prime}(w)$.
Let $X^\prime$ be a subset of $\tR^{b^\prime,-}$ such that
$w(X^\prime)\subset\tR^{b,+}$. Assume that there exists $i\in\fkI$ such that
$w(\checktal^{b^\prime}(i))\in\tR^{b,+}$. Then $-\checktal^{b^\prime}(i)\notin X^\prime$.
Let $b^{\prime\prime}:=\tvsigma_i(b^\prime)$.
Let $X^{\prime\prime}:=\tils^{b^{\prime}}_i(X^\prime)$.
Note that $\tils^{b^{\prime\prime}}_i(-\checktal^{b^{\prime\prime}}(i))=\checktal^{b^\prime}(i)$
and $\tils^{b^{\prime\prime}}_i(\tR^{b^{\prime\prime},-}\setminus\{-\checktal^{b^{\prime\prime}}(i)\})
=\tR^{b^\prime,-}\setminus\{-\checktal^{b^\prime}(i)\}$.
Then $X^{\prime\prime}\subset\tR^{b^{\prime\prime},-}$
and $-\checktal^{b^{\prime\prime}}(i)\notin X^{\prime\prime}$. 
Then the claim~(2) follows form the claim~(1).
\hfill $\Box$

\begin{lemma} {\rm{(See also \cite[Proposition~2.12]{CH09}.)}}
\label{lemma:tRainfty}
Assume $|\tR^a|<\infty$.
Let $\CardtrRp:=|\tR^{a,+}|$. Then $\tR^{a,-}=-\tR^{a,+}$,
and there exist a unique $b\in\mcA[a]$ and a unique $1^aw_0\in\mcH(a,b)$
such that $\tR^{a,+}=1^aw_0(\tR^{b,-})$.
Moreover $\trell_{a,b}(1^aw_0)=\CardtrRp$.
Furthermore, for $f\in\bbM$ with $1^a\tils_{f,\CardtrRp}=1^aw_0$,
we have $\tR^{a,+}=\{1^a\tils_{f,t-1}(\checktal^{a_{f,t-1}}(f(t)))|t\in\fkJ_{1,\CardtrRp}\}$.
\end{lemma}
{\it{Proof.}} Let $b^\prime\in\mcA[a]$, $w\in\mcH(a,b^\prime)$
and $r:=\trell_{a,b^\prime}(w)$. Assume $r\in\fkJ_{0,\CardtrRp-1}$.
By Lemma~\ref{lemma:prptell}(1), there exists $i\in\fkI$ such that
$w(\checktal^{b^\prime}(i))\in\tR^{a,+}$.
Let $b^{\prime\prime}:=\tvsigma_i(b^\prime)$.
By an argument similar to that in the proof of Lemma~\ref{lemma:prptell},
we have $\trell_{a,b^{\prime\prime}}(w\tils^{b^{\prime\prime}}_i)=r+1$.
Then the claim follows from Lemma~\ref{lemma:prptell}.
\hfill $\Box$

\subsection{Finite SGRS} \label{subsection:FiniteSGRS}
Let $(\tR,\tB)$ be an SGRS, i.e., it is as in \eqref{eqn:defstGRS}.
Let $\checktal\in\checktB$.
For $\checktal^\prime\in\checktB[\checktal]$,
let $\tc^{\checktal^\prime}_{ij}:=-
\tN^{\checktal^\prime(\fkI)}_{\checktal^\prime(i),\checktal^\prime(j)}$
($i$, $j\in\fkI$) and $\tC^{\checktal^\prime}:=[\tc^{\checktal^\prime}_{ij}]_{i,j\in\fkI}$,
and define the $\bK$-linear isomorphism $\checkmapu_{\checktal^\prime}:
\newmcV\to\checkV$ by $\checkmapu_{\checktal^\prime}(\checktal^\prime(i))
:=\checkv_i$ ($i\in\fkI$).
Let $\checkmcCctPi:=\mcC(\fkI,(\checktau_i)_{i\in\fkI},\checktB[\checktal],
(\tC^{\checktal^\prime})_{\checktal^\prime\in\checktB[\checktal]})$.
Then $\checkmcCctPi$ is a Cartan scheme.
Define the map $\bartPi:\fkI\to\checkV$ by $\bartPi(i):=\checkv_i$ ($i\in\fkI$).
Let $\checkmcRctPi:=\mcR(\checkmcCctPi,(\checkV,\bartPi,
\checkmapu_{\checktal^\prime}(\tR))_{\checktal^\prime\in\checktB[\checktal]})$.
Then $\checkmcRctPi$ is a CGRS of type $\checkmcCctPi$ 
(where for this $\checkmcRctPi$, we have $\newmcV^{\checktal^\prime}=\checkV$
and $\checktal^{\checktal^\prime}=\bartPi$ for all 
$\checktal^\prime\in\checktB[\checktal]$),
since $1^{\checktal^\prime}s_{f,t}\in\rmGL(\checkV)$
($\checktal^\prime\in\checktB[\checktal]$
$f\in\bbM$, $t\in\bZgeqo$) 
for $\checkmcRctPi$ is the same as in Lemma~\ref{lemma:basech}
(with $\checktal^\prime$ in place of $\checktal$).
By Lemma~\ref{lemma:basech}, we have
\begin{equation}\label{eqn:ucaaft}
(\checkmapu_\checktal\circ\checktal_{f,t})(i)=1^\checktal\tils_{f,t}(\checkv_i)
\quad(f\in\bbM, t\in\bZgeqo, i\in\fkI).
\end{equation}
By Lemma~\ref{lemma:prptell} and \eqref{eqn:ucaaft}, we see that
\begin{equation*}
\mbox{$\checktal^\prime=\checktal^{\prime\prime}$
for $\checktal^\prime$, $\checktal^{\prime\prime}\in\checktB[\checktal]$
with $\checktal^\prime(\fkI)=\checktal^{\prime\prime}(\fkI)$.
}
\end{equation*}

For Lemma~\ref{lemma:PsSGRS} below, recall Lemma~\ref{lemma:BtB}.
\begin{lemma} \label{lemma:PsSGRS}
Assume $|\tR|<\infty$. Let $\tPi\in\tB$.
Let $\CardtrRp:=|\tRptPi|$.
Let $\checktal\in\checktB$ be such that $\checktal(\fkI)=\tPi$.
Then there exists $f\in\bbM$ such that $\checktal_{f,\CardtrRp}(\fkI)=-\tPi$.
Moreover for $f^\prime\in\bbM$,
$\checktal_{f^\prime,\CardtrRp}(\fkI)=-\tPi$
if and only if $\tRptPi=\{\checktal_{f^\prime,t-1}(f^\prime(t))|t\in\fkJ_{1,\CardtrRp}\}$.
\end{lemma}
{\it{Proof.}} This follows from Lemmas~\ref{lemma:prptell},~\ref{lemma:tRainfty}
and \eqref{eqn:ucaaft}.
\hfill $\Box$

\subsection{Relation between an SGRS and a CGRS} \label{subsection:RelSGRSCGRS}
In Subsection~\ref{subsection:RelSGRSCGRS}, fix a Cartan scheme $\mcC=\mcC(\fkI,\mcA,(\tvsigma_i)_{i\in\fkI},(\tC^a)_{a\in\mcA})$, 
and fix a CGRS
$\mcR=\mcR(\mcC,(\newmcV^a,\checktal^a,\tR^a)_{a\in\mcA})$
of type $\mcC$.
Fix $a\in\mcA$.
Let $\tB[a]$ be the set of the subsets $1^a\tils_{f,t}(\checktal^{a_{f,t}}(\fkI))$ of $\tR^a$
with $f\in\bbM$ and $t\in\bZgeqo$.
Then $(\tR^a,\tB[a])$ is an SGRS. 
Let $\checkmcRctPia=\mcR(\checkmcCctPia,(\checkV,\bartPi,
\checkmapu_{\checktal^\prime}(\tR))_{\checktal^\prime\in\checktB[\checktal^a]})$ be those defined 
for the SGRS $(\tR^a,\tB[a])$.
For $a^\prime\in\mcA[a]$, define the $\bK$-linear isomorphism
$\mapu^{a^\prime}:\newmcV^{a^\prime}\to\checkV$ 
by $\mapu^{a^\prime}(\checktal^{a^\prime}(i))
:=\checkv_i$ ($i\in\fkI$).
We have
\begin{lemma}\label{lemma:uaftRaft}
We have
$\mapu^{a_{f,t}}(\tR^{a_{f,t}})=\checkmapu_{\checktal^a_{f,t}}(\tR^a)$
$(f\in\bbM, t\in\bZgeqo)$.
In particular, 
\begin{equation} \label{eqn:uaftRafteq}
c^{a_{f,t}}_{ij}
=c^{\checktal^a_{f,t}}_{ij}\quad(i, j\in\fkI).
\end{equation}
\end{lemma}
{\it{Proof.}} We use an induction on $t$.
The case where $t=0$ is clear.
Assume that $t\geq 1$ and $\mapu^{a_{f,t-1}}(\tR^{a_{f,t-1}})
=\checkmapu_{\checktal^a_{f,t-1}}(\tR^a)$.
By \eqref{eqn:minuscaij}, we have $c^{a_{f,t-1}}_{ij}
=c^{\checktal^a_{f,t-1}}_{ij}$ ($i$, $j\in\fkI$). 
Note $c^{a_{f,t}}_{f(t),j}=c^{a_{f,t-1}}_{f(t),j}$ and
$c^{\checktal^a_{f,t-1}}_{f(t),j}=c^{\checktal^a_{f,t}}_{f(t),j}$
hold for $j\in\fkI$.
Hence $\mapu^{a_{f,t}}=\tils^{\checktal^a_{f,t-1}}_{f(t)}
\circ\mapu^{a_{f,t-1}}\circ\tils^{a_{f,t}}_{f(t)}$
and $\tils^{\checktal^a_{f,t-1}}_{f(t)}
\circ\checkmapu_{\checktal^a_{f,t-1}}=\checkmapu_{\checktal^a_{f,t}}$.
Hence $\mapu^{a_{f,t}}(\tR^{a_{f,t}})
=(\tils^{\checktal^a_{f,t-1}}_{f(t)}
\circ\mapu^{a_{f,t-1}}\circ\tils^{a_{f,t}}_{f(t)})(\tR^{a_{f,t}})
=(\tils^{\checktal^a_{f,t-1}}_{f(t)}
\circ\mapu^{a_{f,t-1}})(\tR^{a_{f,t-1}})
=(\tils^{\checktal^a_{f,t-1}}_{f(t)}
\circ\checkmapu_{\checktal^a_{f,t-1}})(\tR^a)
=\checkmapu_{\checktal^a_{f,t}}(\tR^a)$.
The claim \eqref{eqn:uaftRafteq} follows from \eqref{eqn:minuscaij}.
\hfill $\Box$

\begin{lemma}\label{lemma:aftafptp}
Let $f$, $f^\prime\in\bbM$
and $t$, $t^\prime\in\bZgeqo$ be such that
$\checktal^a_{f,t}=\checktal^a_{f^\prime,t^\prime}$.
Then $a_{f,t}=a_{f^\prime,t^\prime}$. 
\end{lemma}
{\it{Proof.}} It follows from \eqref{eqn:uaftRafteq}
that $(\mapu^a)^{-1}\circ 1^{\checktal^a}\tils_{f^{\prime\prime},t^{\prime\prime}}
\circ\mapu^{a_{f^{\prime\prime},t^{\prime\prime}}} 
=1^a\tils_{f^{\prime\prime},t^{\prime\prime}}$
($f^{\prime\prime}\in\bbM$, $t^{\prime\prime}\in\bZgeqo$).
Since 
$1^{\checktal^a}\tils_{f,t}
=1^{\checktal^a}\tils_{f^\prime,t^\prime}$ holds by \eqref{eqn:ucaaft},
we have $1^a\tils_{f,t}
=1^a\tils_{f^\prime,t^\prime}$.
By Lemma~\ref{lemma:prptell}~(2), $a_{f,t}=a_{f^\prime,t^\prime}$.
\hfill $\Box$
\newline\par
Let $\checktB$ be the one for $\tB[a]$.
By Lemma~\ref{lemma:aftafptp},
we can define the surjection $\mfkf:\checktB[\checktal^a]\to\mcA[a]$
by $\mfkf(\checktal^a_{f,t}):=a_{f,t}$
($f\in\bbM$, $t\in\bZgeqo$).
Then $\mfkf(\checktau_i(\checktal^\prime)):=\tvsigma_i(\mfkf(\checktal^\prime))$
($\checktal^\prime\in\checktB[\checktal^a]$, $i\in\fkI$).
By \eqref{eqn:uaftRafteq}, we have 
$c^{\mfkf(\checktal^\prime)}_{ij}=c^{\checktal^\prime}_{ij}$
($\checktal^\prime\in\checktB[\checktal^a]$, $i$, $j\in\fkI$).

\section{Generalized quantum groups} \label{section:GQGr}

\subsection{Definition of generalized quantum groups}\label{subsection:DefofGQGr}
Let $\CardfkIN$, $\newmcV$ 
and $\fkI$ be as in Subsection~\ref{subsection:setGRS}.
Let $\tPi$ be an $\bR$-basis of $\newmcV$.
Let $\checktal:\fkI\to\newmcV$ be a map with $\checktal(\fkI)=\tPi$.
Define the $\bZ$-submodule $\trtfkAPi$ of $\newmcV$ by
$\trtfkAPi:=\oplus_{i\in\fkI}\bZ\talpii$, so 
$\rmrank_\bZ\trtfkAPi=\CardfkIN(=|\fkI|)$. 
Let $\tbhm:\trtfkAPi\times\trtfkAPi\to\bKt$ be a map 
such that
\begin{equation}\label{eqn:Intrtbhm}
\tbhm(\tlambda,\tmu+\tnu)=\bhm(\tlambda,\tmu)\tbhm(\tlambda,\tnu)\quad
\mbox{and}\quad\bhm(\tlambda+\tnu,\mu)=\bhm(\tlambda,\tmu)\bhm(\tnu,\tmu)
\end{equation}
($\tlambda$, $\tmu$, $\tnu\in\trtfkAPi$).
We call such a map as $\tbhm$ {\it{a bi-homomorphism on $\trtfkAPi$}}.
\begin{equation}\label{eqn:fixPix}
\mbox{From now on until the end of this paper, fix $\checktal$, $\trtfkAPi$ and $\tbhm$ as above.} 
\end{equation}
It is well-known that 
\begin{theorem}\label{theorem:DefofGQG} Note \eqref{eqn:fixPix}. \par
{\rm{(1)}} There exists a unique associative $\bK$-algebra {\rm{(}}with $1${\rm{)}} $\trU$ 
satisfying the following
conditions {\rm{(GQG1)}}-{\rm{(GQG4)}}.
\newline\newline
{\rm{(GQG1)}} As a $\bK$-algebra, $\trU$ is generated by the elements:
\begin{equation}
\trK_\tlambda,\,\trL_\tlambda\,\,(\tlambda\in\tfkAPi),\quad
\trE_i, \trF_i\,\,(i\in \fkI).
\end{equation} 
{\rm{(GQG2)}} The following equations hold:
\begin{equation}\label{eqn:relone}
\begin{array}{l}
\trK_0=\trL_0=1,\,
\trK_\tlambda \trK_\tmu=\tK_{\tlambda+\tmu},\,
\trL_\tlambda \trL_\tmu=\tL_{\tlambda+\tmu},\,
\trK_\tlambda \trL_\tmu=\tL_\tmu \trK_\tlambda, \\
\trK_\tlambda \trE_i =\tbhm(\tlambda,\talpii)\trE_i \trK_\tlambda,\,
\trL_\tlambda \trE_i  =\tbhm(-\talpii,\tlambda)\trE_i \trL_\tlambda,\\
\trK_\tlambda \trF_i  =\tbhm(\tlambda,-\talpii)\trF_i \trK_\tlambda,\,
\trL_\tlambda \trF_i  =\tbhm(\talpii,\tlambda)\trF_i\trL_\tlambda,\\
\mbox{$[\trE_i,\trF_j]=\delta_{ij}(-\trK_\talpii+\trL_\talpii)$}.
\end{array}
\end{equation}
{\rm{(GQG3)}} Regard $\trU\otimes\trU$ as a $\bK$-algebra {\rm{(}}with $1${\rm{)}}
by $(X_1\otimes Y_1)(X_2\otimes Y_2):=X_1X_2\otimes Y_1Y_2$
{\rm{(}}$X_t$, $Y_t\in\trU$ $(t\in\fkJ_{1,2})${\rm{)}}.
Then there exists a $\bK$-algebra homomorphism
$\HopfD:\trU\to\trU\otimes\trU$ such that
$\HopfD(\trK_\tlambda)=\trK_\tlambda\otimes  \trK_\tlambda$,
$\HopfD(\trL_\tlambda)=\tL_\tlambda\otimes  \trL_\tlambda$,
$\HopfD(\trE_i)=\trE_i\otimes  1+\trK_\talpii\otimes\trE_i$,
$\HopfD(\trF_i)=\trF_i\otimes
\trL_\talpii+1\otimes  \trF_i$.
\newline\newline
{\rm{(GQG4)}} 
Let $\trUo$ be the associative $\bK$-algebra {\rm{(}}with $1${\rm{)}} of $\trU$
generated by the elements $\trK_\tlambda\trL_\tmu$ with 
$\tlambda$, $\tmu\in\tfkAPi$. 
Let $\trUp$ {\rm{(}}resp. $\trUm${\rm{)}} be the $\bK$-subalgebra of 
$\trU$
generated by the elements $1$ and $\trE_i$ {\rm{(}}resp. $1$ and $\trF_i${\rm{)}} with all $i\in\fkI$.
Let $\trBp:=\sum_{\tlambda\in\tfkAPi}\trUp\tK_\tlambda$,
and $\trBm:=\sum_{\tlambda\in\tfkAPi}\trUm\tL_\tlambda$.
Then the following conditions {\rm{(GQG4-1)}}-{\rm{(GQG4-3)}} hold. \newline\newline
{\rm{(GQG4-1)}} 
The elements $\trK_\tlambda\trL_\tmu$ with 
$(\tlambda,\tmu)\in\tfkAPi^2$ form a $\bK$-basis of $\trUo$. \newline
{\rm{(GQG4-2)}} The $\bK$-linear map 
$m:\trUm\otimes\trUo\otimes\trUp\to\trU$
defined by $m(Y\otimes Z\otimes X):=YZX$
is a $\bK$-linear isomorphism. \newline\newline
{\rm{(GQG4-3)}} There exists a $\bK$-bilinear map
$\trtvttbhmtPi:\trBp\times\trBm\to\bK$ 
satisfying the following conditions {\rm{(GQG4-3-1)}}-{\rm{(GQG4-3-3)}},
where $\trtvt$ means $\trtvttbhmtPi$.
\newline\newline
{\rm{(GQG4-3-1)}} $\trtvt_{|\trUp\times\trUm}$
is non-degenerate.
\newline\newline
{\rm{(GQG4-3-2)}} The equations 
$\trtvt(\trK_\tlambda,\trL_\tmu)=\tbhm(\tlambda,\tmu)$,
$\trtvt(\trE_i,\trF_j)=\delta_{ij}$,
$\trtvt(\trK_\tlambda,\trF_j)=\trtvt(\trE_i,\trL_\tlambda)=0$ hold.
\newline\newline
{\rm{(GQG4-3-3)}}
For $X_t\in\trBp$ $(t\in\fkJ_{1,2})$
and $Y\in\trBm$ with $\HopfD(Y)=\sum_{y=1}^kY^{(1)}_y\otimes Y^{(2)}_y$,
the equation $\trtvt(X_1X_2,Y)=\sum_{y=1}^k\trtvt(X_1,Y^{(2)}_y)\trtvt(X_2,Y^{(1)}_y)$ holds.
For $X\in\trBp$ with $\HopfD(X)=\sum_{x=1}^rX^{(1)}_x\otimes X^{(2)}_x$
and $Y_t\in\trBm$ $(t\in\fkJ_{1,2})$,
the equation $\trtvt(X,Y_1Y_2)=\sum_{x=1}^r\trtvt(X^{(2)}_x,Y_1)\trtvt(X^{(1)}_x,Y_2)$
holds.
\newline\par
{\rm{(2)}} 
There exists a $\bK$-algebra epimorpism
$\Hopfe:\trU\to\bK$ such that
$\Hopfe(\trK_\tlambda)=\Hopfe(\trL_\tlambda)=1$,
$\Hopfe(\trE_i)=\Hopfe(\trF_i)=0$, and
there exists a $\bK$-algebra anti-automorphism
$\HopfS:\trU\to\trU$ such that
$\HopfS(\trK_\tlambda)=\trK_{-\tlambda}$,
$\HopfS(\trL_\tlambda)=\trL_{-\tlambda}$,
$\HopfS(\trE_i)=-\trK_{-\talpii}\trE_i$,
$\HopfS(\trF_i)=-\trF_i\trL_{-\talpii}$. 
Moreover $\trU$ can be regarded as a Hopf algebra
$(\trU,\HopfD,\HopfS,\Hopfe)$. Furthermore,
letting $\trtvt$ mean $\trtvttbhmtPi$, for
$X\in\trBp$ and $Y\in\trBm$, the following equations hold{\rm{:}}
\begin{equation}\label{eqn:bprtoftB}
\begin{array}{l}
\trtvt(\HopfS(X),Y)=\trtvt(X,\HopfS^{-1}(Y)),
\trtvt(X,1)=\Hopfe(X),
\trtvt(1,Y)=\Hopfe(Y), \\
YX=\sum_{x=1}^k\sum_{y=1}^r
\trtvt(X^{(1)}_x, \HopfS(Y^{(1)}_y))
\trtvt(X^{(3)}_x, Y^{(3)}_y)X^{(2)}_xY^{(2)}_y, \\
XY =\sum_{x=1}^k\sum_{y=1}^r
\trtvt(X^{(3)}_x, \HopfS(Y^{(3)}_y))
\trtvt(X^{(1)}_x, Y^{(1)}_y) Y^{(2)}_yX^{(2)}_x,
\end{array}
\end{equation} where $((\rmid_\trU\otimes\HopfD)\circ\HopfD)(X)=
\sum_{x=1}^k X^{(1)}_x\otimes X^{(2)}_x\otimes X^{(3)}_x$ and 
$((\rmid_\trU\otimes\HopfD)\circ\HopfD)(Y)=
\sum_{y=1}^r Y^{(1)}_y\otimes Y^{(2)}_y\otimes Y^{(3)}_y$.
\end{theorem}
{\it{Proof.}} This theorem can be proved in a well-known argument introduced by Drinfeld~\cite{Dr86}.
For a detailed and essentially the same argument, see \cite[CHAPTER~3]{b-Lusztig93},
\cite[Proof of Theorem~2.9.4]{Y94},
\cite[Subsections~6.3, 6.4]{Y99}. \hfill $\Box$
\newline\par
We have fixed $\tbhm$ and $\checktal$, see \eqref{eqn:fixPix}.
For a bi-homomorphism $\tbhm^\prime$ on $\tfkAPi$ and 
a map $\checktal^\prime:\fkI\to\newmcV$ such that $\checktal^\prime(\fkI)$
is an $\bR$-basis of $\newmcV$, let  
$\trU(\tbhm^\prime,\checktal^\prime)$ 
(resp. $\trUo(\tbhm^\prime,\checktal^\prime)$, resp. $\trUp(\tbhm^\prime,\checktal^\prime)$,
resp. $\trUm(\tbhm^\prime,\checktal^\prime)$)
to mean the one defined in the same way as that for $\trU$ (resp. $\trUo$, resp. $\trUp$,
resp.  $\trUm$) of Theorem~\ref{theorem:DefofGQG} 
with $\tbhm^\prime$ and $\checktal^\prime$ in place of $\tbhm$ and $\checktal$.
\begin{equation}\label{eqn:abbrU}
\begin{array}{l}
\mbox{We also often use
the symbols $\trU$, $\trUo$, $\trUp$
and $\trUm$,} \\
\mbox{which mean $\trUtbhmtPi$, $\trUotbhmtPi$, $\trUptbhmtPi$
and $\trUmtbhmtPi$ respectively.}
\end{array}
\end{equation} 

Let $\tfkAPip:=\oplus_{i\in\fkI}\bZgeqo\talpii(\subset\tfkAPi)$.
Recall \eqref{eqn:abbrU}.
We can also regard $\trU$ as a $\tfkAPi$-graded algebra
\begin{equation}\label{eqn:abbrUgr}
\trU=\oplus_{\tlambda\in\tfkAPi}\trU_\tlambda
\end{equation}
such that
$\trU_\tlambda\trU_\tmu\subset \trU_{\tlambda+\tmu}$
($\tlambda$, $\mu\in\tfkAPi$)
$\trE_i\in\trU_\talpii$, $\trF_i\in\trU_{-\talpii}$
($i\in\fkI$)
and $\trUo\subset\trU_0$.
Let $\trUp_\tlambda:=\trUp\cap\trU_\tlambda$ and $\trUm_\tlambda:=\trUm\cap\trU_\tlambda$
($\tlambda\in\tfkAPi$).
Then we have $\trUp=\oplus_{\tlambda\in\tfkAPip}\trUp_\tlambda$ and 
$\trUm=\oplus_{\tlambda\in\tfkAPip}\trUm_{-\tlambda}$.
We can easily see that 
\begin{equation} \label{eqn:protvt}
\begin{array}{l}
\trtvttbhmtPi(X\tK_\tlambda,Y\tL_\tmu)=\delta_{\tnu,\tnu^\prime}\cdot\tbhm(\tlambda,\tmu)\trtvttbhmtPi(X,Y) \\
\quad (\tlambda,\tmu\in\tfkAPi,\,
\tnu,\tnu^\prime\in\tfkAPip,\,X\in\trUp_\tnu,
Y\in\trUm_{-\tnu^\prime}),
\end{array}
\end{equation}
\begin{equation} \label{eqn:dimUmlUpl}
\dim\trUm_{-\tlambda}=\dim\trUp_\tlambda\quad (\tlambda\in\tfkAPip),
\end{equation} and that
\begin{equation} \label{eqn:dimUmlUpld}
\begin{array}{l}
\mbox{${\trtvttbhmtPi}_{|\trUp_\tlambda\times\trUm_{-\tlambda}}$ is non-degenerate for $\tlambda\in\tfkAPip$} \\
\mbox{and ${\trtvttbhmtPi}(\trUp_\tmu,\trUm_{-\tnu})=\{0\}$
for $\tmu$, $\tnu\in\tfkAPip$ with $\tmu\ne\tnu$.}
\end{array}
\end{equation} Hence using \eqref{eqn:bprtoftB}, we can see Remark~\ref{remark:GQGthreed} below.
\begin{remark}\label{remark:GQGthreed}
As a $\bK$-algebra, 
$\trU=\trUtbhmtPi$ can also be characterized by {\rm{(GQG1)}}, {\rm{(GQG2)}}, {\rm{(GQG4-1)}}, {\rm{(GQG4-2)}} and 
\begin{equation*}
\mbox{{\rm{(GQG4-$3^\prime$)}}  
$\{Y\in\trUm|\forall i\in\fkI, [\trE_i,Y]=0\}
=\{X\in\trUp|\forall i\in\fkI, [X,\trF_i]=0\}=\bK\cdot 1$.}
\end{equation*} See also \cite[(4.8)]{AYY15}.
\end{remark}

Let $i\in\fkI$.
Let $\trUtbhmtPiith$ be the $\bK$-subalgebra of $\trUtbhmtPi$
generated by $\trK_\tlambda$, $\trL_\tlambda$ ($\tlambda\in\tfkAPi$)
and $\trE_i$, $\trF_i$. Then we can easily see: 
\begin{lemma}\label{lemma:Uith}
Let $\ckpchi:=\kpch(\tbhm(\talpii,\talpii))$. Then we have the following.
\par {\rm{(1)}}
If $\ckpchi=0$ {\rm{(}}resp. $\ckpchi\ne 0${\rm{)}}, then as a $\bK$-algebra, $\trUtbhmtPiith$ can also be defined by the generators
$\trK_\tlambda$, $\trL_\tlambda$ $(\tlambda\in\tfkAPi)$
and $\trE_i$, $\trF_i$ and the relations composed of those of \eqref{eqn:relone}
{\rm{(}}resp. those of \eqref{eqn:relone} and 
$\trE_i^\ckpchi=\trF_i^\ckpchi=0${\rm{)}}.
\par {\rm{(2)}} If $\ckpchi=0$ {\rm{(}}resp. $\ckpchi\ne 0${\rm{)}}, then as a $\bK$-linear space,
the elements $\trF_i^y\trK_\tlambda\trL_\tmu\trE_i^x$ with $\tlambda$, $\tmu\in\tfkAPi$
and $x$, $y\in\bZgeqo$ {\rm{(}}resp. $x$, $y\in\fkJ_{0,\ckpchi-1}${\rm{)}}
form a $\bK$-basis of $\trUtbhmtPiith$.

\par {\rm{(3)}} Let $\checktal^\prime:\fkI\to\newmcV$ be a map such that $\checktal^\prime(\fkI)$
is a $\bZ$-base of $\tfkAPi$. Let $j\in\fkI$.
Assume $\checktal^\prime(j)=-\talpii$.
Then there exists
a unique $\bK$-algebra isomorphism
$\smtrTtbtali:\trU(\tbhm,\checktal^\prime;j)\to\trUtbhmtPiith$
such that $\smtrTtbtali(\trK_\tlambda)=\trK_\tlambda$, 
$\smtrTtbtali(\trL_\tlambda)=\trL_\tlambda$ $(\tlambda\in\tfkAPi)$
and $\smtrTtbtali(\trE_j)=\trF_i\trL_{-\talpii}$,
$\smtrTtbtali(\trF_j)=\trK_{-\talpii}\trE_i$.
\end{lemma}

\subsection{Kharchenko's PBW theorem}\label{subsection:sbsecKhaPBW}

For the statement of Theorem~\ref{theorem:KhaPBW},
we need the notations as follows.
Let $Y$ be a non-empty subset of $\tfkAPip$.
Let $z:Y\to\bN$ be a map.
Define the subset $Y^{\langle z\rangle}$ of $Y\times\bN$
by $Y^{\langle z\rangle}:=\{(\lambda,k)|\lambda\in Y, k\in\fkJ_{1,z(\lambda)}\}$.
Define the surjection $p^{\langle Y, z\rangle}:Y^{\langle z\rangle}\to Y$ by
$p^{\langle Y, z\rangle}(y)=\lambda$ for $y:=(\lambda,k)\in Y^{\langle z\rangle}$.
Let $\trRpmap^{\langle Y, z\rangle}$ be the set of 
maps $f:Y^{\langle z\rangle}\to\bZgeqo$
such that 
\begin{equation*}
\mbox{$(f(y))_{\tbhm(p^{\langle Y, z\rangle}(y),p^{\langle Y, z\rangle}(y))}!\ne 0$
for all $y\in Y^{\langle z\rangle}$.}
\end{equation*} For $\tlambda\in\tfkAPip$, let
\begin{equation*}
\trRpmap^{\langle Y, z\rangle}_\tlambda
:=\{f\in\trRpmap^{\langle Y, z\rangle}|\sum_{y\in Y^{\langle z\rangle}}f(y)p^{\langle Y, z\rangle}(y)=\tlambda\}.
\end{equation*}

The following is well-known.

\begin{theorem}\label{theorem:KhaPBW}
{\rm{(}}Kharchenko's PBW theorem {\rm{\cite{Kha99}}}{\rm{)}}
Note \eqref{eqn:abbrU}.

{\rm{(1)}} There exists a unique pair
\begin{equation*}
(\trRptbhmtPi,\trvrpptbhmtchecktal) 
\end{equation*}
of a subset $\trRptbhmtPi$ of $\tfkAPip\setminus\{0\}$ and
a map $\trvrpptbhmtchecktal:\trRptbhmtPi\to\bN$ such that
\begin{equation*}
\dim\trUp_\tlambda=|\trRpmap^{\langle\trRptbhmtPi,\trvrpptbhmtchecktal\rangle}_\tlambda|
\quad\mbox{for all $\tlambda\in\tfkAPip$.}
\end{equation*} 

{\rm{(2)}} 
Let
$X:=\trRptbhmtPi^{\langle\trvrpptbhmtchecktal\rangle}$ and
$p:=p^{\langle\trRptbhmtPi,\trvrpptbhmtchecktal\rangle}$.
Let $Z:=\trRpmap^{\langle\trRptbhmtPi,\trvrpptbhmtchecktal\rangle}$.
Let $Z^\prime$ be the subset of $Z$ of elements
$f$ with $|\{x\in X|f(x)\ne 0\}|<\infty$.
Then there exist a total order $\preceq$ on $X$ and
elements ${\grave E}_x$ of $\trUp_{p(x)}$ for $x\in X$ such that
\begin{equation*}
\mbox{the elements $\prod^{\preceq}_{x\in X}{\grave E}_x^{f(x)}$
with $f\in Z^\prime$
form a $\bK$-basis of $\trUp$,}
\end{equation*} where
$\prod^{\preceq}_{x\in X}{\grave E}_x^{f(x)}:=
{\grave E}_{x_1}^{f(x_1)}{\grave E}_{x_2}^{f(x_2)}\cdots{\grave E}_{x_k}^{f(x_k)}$
for $f$ and $x_1,\ldots,x_k$ with $x_1\prec\cdots\prec x_k$ and
$f(x^\prime)=0$ $(x^\prime\in X\setminus\{x_1,\ldots,x_k\})$.

{\rm{(3)}} The same claim as that of {\rm{(2)}}
with $\trUm$ and $\trUm_{-p(x)}$ in place of
$\trUp$ and $\trUp_{p(x)}$ respectively holds,
where $\preceq$ is the same as the one of {\rm{(2)}}.

\end{theorem}

More precisely, as for Theorem~\ref{theorem:KhaPBW}~(3), 
we have the $\bK$-algebra automorphism
$\Omegatbhmtchecktal$ of $\trU$ with 
$\Omegatbhmtchecktal(\trK_\tlambda\trL_\tmu):=\trK_{-\tlambda}\trL_{-\tmu}$
($\tlambda$, $\tmu\in\tfkAPi$),
$\Omegatbhmtchecktal(\trE_i):=\trF_i\trL_{-\talpii}$
$\Omegatbhmtchecktal(\trF_i):=\trK_{-\talpii}\trE_i$
($i\in\fkI$),
see \cite[Subsection~4.2]{AYY15} for example.

\begin{equation*}
\begin{array}{l}
\mbox{Let $\trRtbhmtPi:=\trRptbhmtPi\cup(-\trRptbhmtPi)$. Define the map $\trvrptbhmtchecktal:\trRtbhmtPi\to\bN$} \\
\mbox{by
$\trvrptbhmtchecktal(-\tal):=\trvrptbhmtchecktal(\tal):=\trvrpptbhmtchecktal(\tal)$
($\tal\in\trRptbhmtPi$).}
\end{array}
\end{equation*}
It is clear that
\begin{equation}\label{eqn:simplpr}
\bZ\talpii\cap\trRtbhmtPi=\{-\talpii,\talpii\},\,\trvrptbhmtchecktal(\talpii)=\trvrptbhmtchecktal(-\talpii)=1
\quad(i\in\fkI).
\end{equation}

For $i$, $j\in\fkI$ with $i\ne j$, there exists
$\trNtbhmtPi_{i,j}\in\bZgeqo\cup\{\pinfty\}$ such that
\begin{equation}\label{eqn:defNij}
\begin{array}{lcl}
\fkJ_{0,\trNtbhmtPi_{i,j}}&=&\{\,m\in\bZgeqo\,|\,\talpij+m\talpii\in\trRtbhmtPi\,\} \\
&=& \{\,m\in\bZgeqo\,|\,(m)_{\tq_{ii}}!(m;\tq_{ii},\tq_{ij}\tq_{ji})!\ne
0\,\},
\end{array}
\end{equation} where we let $\tq_{ik}:=\tbhm(\talpii,\talpik)$
for $k\in\fkI$ (see \cite[(2.17)]{HY10} for the second equation for
example).
Let $\trNtbhmtPi_{i,i}:=-2$ for $i\in\fkI$.

Let $\tBtbhmchPi$ be the set of all bases of $\trRtbhmtPi$.
Let $\checktBtbhmchPi$ be the set of all maps $\checktal^\prime:\fkI\to\newmcV$ 
with $\checktal^\prime(\fkI)\in\tBtbhmchPi$.
Let $i\in\fkI$. Define the map $\gravetau_i:\checktBtbhmchPi\to\checktBtbhmchPi$ by
\begin{equation}\label{eqn:Defgrvt}
\left\{\begin{array}{ll}
\gravetau_i(\checktal^\prime)(j):=
\checktal^\prime(j)+\trN^{\tbhm,\checktal^\prime}_{i,j}\checktal^\prime(i) \,\,(j\in\fkI)
& \mbox{if $\trN^{\tbhm,\checktal^\prime}_{i,j^\prime}\ne \pinfty$ for all $j^\prime\in\fkI$}, \\
\gravetau_i(\checktal^\prime):=\checktal^\prime
& \mbox{otherwise}
\end{array}\right.
\end{equation} ($\checktal^\prime\in\checktBtbhmchPi$).

The following proposition is well-known.
\begin{proposition} \label{proposition:HecProp}
{\rm{({\it{See}} \cite[Proposition~1]{Hec06}.)}}
Let $\checktal^\prime\in\checktBtbhmchPi$. Then we have 
the following equations \eqref{eqn:pprprvprt}-\eqref{eqn:prvprtaf}.
\begin{equation}\label{eqn:pprprvprt}
\trN^{\tbhm,\gravetau_i(\checktal^\prime)}_{i,j}=\trN^{\tbhm,\checktal^\prime}_{i,j}
\quad\quad(i,j\in\fkI).
\end{equation}
\begin{equation}\label{eqn:prprvprt}
\gravetau_i^2=\rmid_\checktBtbhmchPi.
\end{equation}
\begin{equation}\label{eqn:prvprt}
\trR(\tbhm,\gravetau_i(\checktal^\prime))=\trR(\tbhm,\checktal^\prime),
\end{equation} that is, 
\begin{equation}\label{eqn:prvprtd}
\trRp(\tbhm,\gravetau_i(\checktal^\prime))=
\left\{\begin{array}{ll}
\{-\talpii\}\cup(\trRp(\tbhm,\checktal^\prime)\setminus\{\talpii\})
& \mbox{if $\trN^{\tbhm,\checktal^\prime}_{i,j}\ne \pinfty$ for all $j\in\fkI$}, \\
\trRp(\tbhm,\checktal^\prime) & \mbox{otherwise}.
\end{array}\right.
\end{equation}
\begin{equation}\label{eqn:prvprtaf}
\trvrp^{\tbhm,\gravetau_i(\checktal^\prime)}(\tal)=\trvrp^{\tbhm,\checktal^\prime}(\tal)
\quad(\tal\in\trRtbhmtPi).
\end{equation}
\end{proposition}

\begin{lemma}\label{lemma:trRinfty}
Assume $|\trRtbhmtPi|<\infty$. 
Then $(\trRtbhmtPi,\tBtbhmchPi)$ is a connected SGRS
with $\checktal^\prime(\fkI)^{(\checktal^\prime(i))}=\gravetau_i(\checktal^\prime)(\fkI)$
$(\checktal^\prime\in\checktBtbhmchPi,i\in\fkI)$. 
In particular, for every $\tal\in\trRtbhmtPi$, there exists a base $B$ of $\trRtbhmtPi$
such that $\tal\in B$. {\rm{(}}$B$ is a $\bZ$-base of $\tfkAPi$.{\rm{)}}
Moreover $\trvrptbhmtchecktal(\tal)=1$ for all $\tal\in\trRtbhmtPi$.
\end{lemma}
{\it{Proof.}} The claim follows from Lemma~\ref{lemma:BtB}
and \eqref{eqn:simplpr}, \eqref{eqn:prvprt}, \eqref{eqn:prvprtaf}.
\hfill $\Box$

\begin{theorem}{\rm{({\it{See}} \cite[{\it{Theorem}}~6.11]{Hec10}.)}}\label{theorem:Liso}
Let $i\in\fkI$ be such that
$\trNtbhmtPi_{ij}\ne \pinfty$ for all $j\in\fkI$.
Then there exists
a unique $\bK$-algebra isomorphism
\begin{equation}\label{eqn:LuIso}
\trTtbtali:\trUtbtalni\to\trUtbhmtPi
\end{equation}
satisfying the following equations \eqref{eqn:LuIsoeqn},
where we mean
$\trT_i:=\trTtbtali$, $\trN_{ij}:=\trNtbhmtPi_{ij}$
and $\tq_{tr}:=\tbhm(\talpit,\talpir)$  $(t,r\in\fkI)$.
\begin{equation}\label{eqn:LuIsoeqn}
\begin{array}{l}
\trT_i(X)=\trT^{(\tbhm,\checktal,\tPiangi;i)}(X)\quad(X\in\trU(\tbhm,\tPiangi;i)), \\
\trT_i(\trE_j)= \sum_{k=0}^{N_{ij}}(-\tq_{ij})^k\tq_{ii}^{\frac
{k(k-1)} {2}}
{{N_{ij}}\choose{k}}_{\tq_{ii}}\trE_i^{\trN_{ij}-k}E_jE_i^k, \\
\trT_i(\trF_j)= {\frac 1
{(\trN_{ij})_{\tq_{ii}}!(\trN_{ij};\tq_{ii},\tq_{ij}\tq_{ji})!}}
\sum_{k=0}^{N_{ij}}(-\tq_{ji})^k\tq_{ii}^{\frac {k(k-1)} {2}}
{{N_{ij}}\choose{k}}_{\tq_{ii}}\trF_i^{\trN_{ij}-k}\trF_j\trF_i^k,
\end{array}
\end{equation} {\rm{(}}$\tal\in\tfkAPi$,
$j\in\fkI\setminus\{i\}${\rm{)}},
where recall Lemma~{\rm{\ref{lemma:Uith}~(3)}} for
$\trT^{(\tbhm,\checktal,\tPiangi;i)}$.
{\rm{(}}Note that $\trE_j$ {\rm{(}}resp. $\trF_j${\rm{)}} of LRSs of \eqref{eqn:LuIsoeqn} 
is an element of 
$\trUtbtalni_{\talpij+\trN_{ij}\talpii}$ {\rm{(}}resp. 
$\trUtbtalni_{-\talpij-\trN_{ij}\talpii}${\rm{)}}
since $\tPiangi(j)=\talpij+\trN_{ij}\talpii$.{\rm{)}}
In particular,
\begin{equation*}
\trTtbtali(\trUtbtalni_\tlambda)=\trUtbhmtPi_\tlambda\quad
(\tlambda\in\tfkAPi).
\end{equation*}
\end{theorem}

Let $\CardtrRp:=\originalCardtrRp$. Assume $\CardtrRp<\infty$.
By Lemma~\ref{lemma:trRinfty}, $(\trRtbhmtPi,\tBtbhmchPi)$ is an SGRS. 
For $f\in\bbM$ and $t\in\bZgeqo$, 
recall $\checktal_{f,t}\in\checktBtbhmchPi$ from \eqref{eqn:checktalft}, and 
define the $\bK$-linear isomorphism
$\onebctaltrT_{f,t}:\trU(\tbhm,\checktal_{f,t})\to\trUtbhmtPi$ by
\begin{equation*}
\onebctaltrT_{f,0}:=\rmid_{\trUtbhmtPi}, 
\quad\mbox{and}\quad
\onebctaltrT_{f,t}:=\onebctaltrT_{f,t-1}\circ\trT^{\tbhm,\checktal_{f,t}}_{f(t)}
\,\,(t\in\bN).
\end{equation*}
\begin{equation}\label{eqn:defhatf}
\begin{array}{l}
\mbox{Let $\hatf\in\bbM$ be such that $\checktal_{\hatf,\CardtrRp}(\fkI)=-\tPi$.} \\ 
\mbox{Note that $\hatf(1)$ can be any element of $\fkI$.} \\
\mbox{Let $\dotbeta_t:=\checktal_{\hatf,t-1}(\hatf(t))$ ($t\in\fkJ_{1,\CardtrRp}$).}
\end{array}
\end{equation}
By Lemmas~\ref{lemma:PsSGRS} and \ref{lemma:trRinfty}, we have 
\begin{equation}\label{eqn:trRptbhmtPi}
\trRptbhmtPi=\{\,\dotbeta_t\,|\,t\in\fkJ_{1,\CardtrRp}\,\}.
\end{equation}
Let $t\in\fkJ_{1,\CardtrRp}$. Let
\begin{equation}\label{eqn:DEFdtEF}
\dottrE_t:=\onebctaltrT_{\hatf,t-1}
(\trE_{\hatf(t)}),
\quad
\dottrF_t:=\onebctaltrT_{\hatf,t-1}
(\trF_{\hatf(t)}).
\end{equation}
It is clear that $\dottrE_t\in\trUtbhmtPi_{\dotbeta_t}$
and $\dottrF_t\in\trUtbhmtPi_{-\dotbeta_t}$.
By Lemma~\ref{lemma:Uith}~(1), 
for $k\in\fkJ_{2,\infty}$,
\begin{equation}\label{eqn:abcequi}
(k)_{\tbhm(\dotbeta_t,\dotbeta_t)}!=0\,\,\Leftrightarrow\,\, 
\dottrE_t^k=0\,\,\Leftrightarrow\,\,\dottrF_t^k=0.
\end{equation}

\begin{theorem}\label{theorem:LusPBW}
{\rm{(}}See {\rm{\cite[Theorems~4.8,\,4.9]{HY10}}}{\rm{)}} Note \eqref{eqn:abbrU}.
Let $\CardtrRp:=\originalCardtrRp$. Assume $\CardtrRp<\infty$.
Let $\dottrE_t$, $\dottrF_t$ $(t\in\fkJ_{1,\CardtrRp})$ be as above.
Recall the convention~\eqref{eqn:abbrU}.
Then for any bijection $x:\fkJ_{1,\CardtrRp}\to\fkJ_{1,\CardtrRp}$,
the elements 
\begin{equation*}
\dottrE_{x(1)}^{k_{x(1)}}
\cdots\dottrE_{x(\CardtrRp)}^{k_{x(\CardtrRp)}}
\,\,(\mbox{resp.}\,\dottrF_{x(1)}^{k_{x(1)}}\cdots\dottrF_{x(\CardtrRp)}^{k_{x(\CardtrRp)}}))\quad
(k_t\in\bZgeqo,\,(k_t)_{\tbhm(\dotbeta_t,\dotbeta_t)}!\ne 0\,(t\in\fkJ_{1,\CardtrRp}))
\end{equation*}
form a $\bK$-basis of $\trUp$
{\rm{(}}resp. $\trUm${\rm{)}}.
In particular, $\dottrE_t\in\trUp_{\dotbeta_t}$
and $\dottrF_t\in\trUm_{-\dotbeta_t}$ for $t\in\fkJ_{1,\CardtrRp}$.
\end{theorem}
{\it{Proof.}} 
Let $X$ be the set of all bi-homomorphisms
$\tbhm^\prime$ on $\trtfkAPi$ with $|\trR(\tbhm^\prime,\checktal)|<\infty$.
Then $\tbhm\in X$.
Let $c^{\tbhm^\prime}_{ij}:=-\trN^{\tbhm^\prime,\checktal}_{i,j}$
($\tbhm^\prime\in X$, $i,j\in\fkI$).
Let $\tC^{\tbhm^\prime}=[\tc^{\tbhm^\prime}_{ij}]_{i,j\in\fkI}\in\Mat(\CardfkIN,\bZ)$.
For $i\in\fkI$, define the map
${\dot\tau}_i:X\to X$ by ${\dot\tau}_i(\tbhm^\prime)(\talpij,\talpik):=
\tbhm^\prime(\talpij-c^{\tbhm^\prime}_{ij}\talpii,\talpik-c^{\tbhm^\prime}_{ik}\talpii)$.
($j$, $k\in\fkI$).
Then $\mcC_X:=\mcC(X,\fkI,({\dot\tau}_i)_{i\in\fkI},
(\tC^{\tbhm^\prime})_{\tbhm^\prime\in X})$ is a Cartan scheme,
and $\mcR_X:=\mcR(\mcC_X,(\newmcV,\checktal,\trR(\tbhm^\prime,\checktal))_{\tbhm^\prime\in X})$
is a CGRS of type $\mcC_X$.
Recall from \eqref{eqn:mcCamcRa} that
$\mcR_X[\tbhm]$ is a CGRS of type $\mcC_X[\tbhm]$.
Note that for $f\in\bbM$ and $t\in\bZgeqo$, there exists a $\bK$-algebra isomorphism
$\nblIso_{f,t}:\trU(\tbhm_{f,t},\checktal)\to\trU(\tbhm,\checktal_{f,t})$
such that $\nblIso_{f,t}(\trK_\talpii)=\trK_{\checktal_{f,t}(i)}$,
$\nblIso_{f,t}(\trL_\talpii)=\trL_{\checktal_{f,t}(i)}$,
$\nblIso_{f,t}(\trE_i)=\trE_i$ and $\nblIso_{f,t}(\trF_i)=\trF_i$
($i\in\fkI$). Note $1^\checktal\tils_{\hatf,l}(\tPi)=-\tPi$.
Then the claim follows from \cite[Theorems~4.8,\,4.9]{HY10}.
\hfill $\Box$
\newline\par Let $\CardtrRp:=\originalCardtrRp$. Assume $\CardtrRp<\infty$.
Recall \eqref{eqn:defhatf}. Then we have
\begin{equation}\label{eqn:extlgdo}
\begin{array}{lcl}
\trRp(\tbhm,\checktal_{\hatf,1})&=&\{-\talpidfone\}\cup
(\trRptbhmtPi\setminus\{\talpidfone\}) \\
&=&\{\,\dotbeta_t\,|\,t\in\fkJ_{2,\CardtrRp}\,\}
\cup\{-\talpidfone\}.
\end{array}
\end{equation} Moreover we have:
\begin{lemma}\label{lemma:extlgd} Let $i\in\fkI$.
Let $\hatf$ of \eqref{eqn:defhatf} be such that
$i=\hatf(1)$.
Let $j\in\fkI$ be such that
$\checktal_{\hatf,\CardtrRp}(j)=-\talpii$.
Let $\hatf^\prime\in\bbM$ be such that $\hatf^\prime(t)=\hatf(t+1)$
$(t\in\fkJ_{1,\CardtrRp -1})$ and $\hatf^\prime(\CardtrRp)=j$. 
{\rm{(}}Recall $\checktal_{\hatf,1}=\tPiangi$.{\rm{)}} Then we have
\begin{equation}\label{eqn:extlgda}
\tPiangi_{\hatf^\prime,t-1}(\hatf^\prime(t))=\dotbeta_{t+1}\,(t\in\fkJ_{1,\CardtrRp-1}),\,\,
\tPiangi_{\hatf^\prime,\CardtrRp-1}(\hatf^\prime(\CardtrRp))=-\talpii,
\end{equation} 
\begin{equation}\label{eqn:extlgdb}
\tPiangi_{\hatf^\prime,\CardtrRp}(\fkI)=-\tPiangi(\fkI),
\end{equation} and
\begin{equation}\label{eqn:extlgdc}
\begin{array}{l}
1^{\tbhm,\tPiangi}\trT_{\hatf^\prime,\CardtrRp-1}(\trE_{\hatf^\prime(\CardtrRp)})
=x\trE_i\,
(\in\trUp(\tbhm,\tPiangi)_{-\talpii}), \\
1^{\tbhm,\tPiangi}\trT_{\hatf^\prime,\CardtrRp-1}(\trF_{\hatf^\prime(\CardtrRp)})
=x^{-1}\trF_i\,
(\in\trUp(\tbhm,\tPiangi)_\talpii)
\end{array}
\end{equation} for some $x\in\bKt$. Moreover
for $\tlambda\in\tfkAPipni$,
\begin{equation}\label{eqn:extlgdd}
\begin{array}{l}
\trTtbtali(\trUptbtalni_\tlambda)\subset
\oplus_{k=0}^\infty\rmSpan_\bK(\trUm_{-k\talpii}\trUo\trUp_{\tlambda+k\talpii}), \\
\trTtbtali(\trUmtbtalni_{-\tlambda})\subset
\oplus_{k=0}^\infty\rmSpan_\bK(\trUm_{-\tlambda-k\talpii}\trUo\trUp_{k\talpii}).
\end{array}
\end{equation}
\end{lemma}
{\it{Proof.}} The equations of \eqref{eqn:extlgda} are clear since
$\tPiangi_{\hatf^\prime,t-1}=\checktal_{\hatf,t}$ for $t\in\fkJ_{1,\CardtrRp}$.
Then \eqref{eqn:extlgdb} follows from Lemma~\ref{lemma:PsSGRS}.
Then \eqref{eqn:extlgdc} and \eqref{eqn:extlgdd} follows from Theorem~\ref{theorem:LusPBW}.
\hfill $\Box$

\section{Singular vectors}\label{section:SingVec}
\subsection{Verma modules}\label{subsection:VermaMd}
Recall the convention \eqref{eqn:fixPix} and \eqref{eqn:abbrU}.
By {\rm{(GQG4-2)}}, for $\trLambda\in\rmChtbhmchecktal$, we have a unique
left $\trU$-module $\MtbctaltrLam$
with a non-zero element $\tv_\trLambda\in\MtbctaltrLam$
satisfying the conditions ($v1$) and ($v2$) below.
\newline\par
($v1$) Equations $Z\tv_\trLambda=\trLambda(Z)\tv_\trLambda$
($Z\in\trUo$) and  $\trE_i\tv_\trLambda=0$ ($i\in\fkI$) hold. \par
($v2$) The $\bK$-linear homomorphism
$\trUm\to\MtbctaltrLam$ defined by $Y\mapsto Y\tv_\trLambda$
is bijective.
\newline\par
Let $\trLambda\in\rmChtbhmchecktal$. For $i\in\fkI$, 
define 
$\trLambdaangi\in\rmCh(\trUo(\tbhm,\tPiangi))$ by 
\begin{equation}\label{eqn:Lmbdai}
\trLambdaangi(\trK_{\tlambda}\trL_{\tmu}):=
{\frac {\tbhm(\talpii,\tmu)^{\ckpchi-1}} 
{\tbhm(\tlambda,\talpii)^{\ckpchi-1}}}
\trLambda(\trK_{\tlambda}\trL_{\tmu})\quad(\tlambda,\,\tmu\in\tfkAPi),
\end{equation} where $\ckpchi:=\kpch(\tbhm(\talpii,\talpii))$. 

We can easily see:
\begin{theorem}\label{theorem:LisoML}
{\rm{(\cite[{\it{Proposition}}~5.11]{HY10})}} Let $\trLambda\in\rmChtbhmchecktal$. 
Let $i\in\fkI$. Let $\ckpchi:=\kpch(\tbhm(\talpii,\talpii))$.
Assume $\ckpchi\ne 0$. 
Assume that $\trLambda(-\trK_\talpii+\tbhm(\talpii,\talpii)^t\trL_\talpii)$ $\ne 0$
for all $t\in\fkJ_{0,\ckpchi-2}$.
Then there exists a $\bK$-linear isomorphism
\begin{equation*}
\hTtbtali:\mclM^{\tbhm,\tPiangi}(\trLambdaangi)\to 
\MtbctaltrLam
\end{equation*} such that
\begin{equation*}
\hTtbtali(X\tv_\trLambdaangi)=\trTtbtali(X)\trF_i^{\ckpchi-1}\tv_\trLambda\quad
(X\in\trUtbtalni).
\end{equation*} 
\end{theorem}

Define the $\bZ$-module homomorphism $\trhobctal:\tfkAPi\to\bKt$ 
by 
\begin{equation*}
\trhobctal(\talpij):=\tbhm(\talpij,\talpij)\quad (j\in\fkI).
\end{equation*}

\begin{lemma}\label{lemma:lmbbrho}
{\rm{(}}See also {\rm{\cite[{\it{Lemma}}~2.17]{HY10}.}}{\rm{)}}
Let $i\in\fkI$, and $\ckpchi:=\kpch(\tbhm(\talpii,\talpii))$. Then we have
\begin{equation}\label{eqn:tbtbrho}
\forall
\tlambda\in\tfkAPi,\,\,\tbhm(\talpii,\tlambda)^{\ckpchi-1}\tbhm(\tlambda,\talpii)^{\ckpchi-1}
={\frac {\trho^{\tbhm,\tPiangi}(\tlambda)} {\trhobctal(\tlambda)}}.
\end{equation}
{\rm{(}}This lemme is a wider version of {\rm{\cite[{\it{Lemma}}~2.17]{HY10}}}
where $\ckpchi\ne 0$ was assumed.{\rm{)}}
\end{lemma}
{\it{Proof.}} For $i$, $j\in\fkI$, let
$\tq_{ij}:=\tbhm(\talpii,\talpij)$ and $\trN_{ij}:=\trNtbhmtPi_{i,j}$.
so $\tPiangi(j)=\talpij+\trN_{ij}\talpii$ by \eqref{eqn:Defgrvt}.
Let
$f$, $f^\prime:\tfkAPi\to\bKt$ be
the two $\bZ$-module homomorphisms
sending $\tlambda\in\tfkAPi$ to
the same values as LHS and RHS of \eqref{eqn:tbtbrho}
respectively.
Let $j\in\fkI$. Then $f(\talpij)=\tq_{ij}^{\ckpchi-1}\tq_{ji}^{\ckpchi-1}$,
and
\begin{equation*}
\begin{array}{l}
f^\prime(\talpij)=\trho^{\tbhm,\tPiangi}(\talpij)\trhobctal(\talpij)^{-1} \\
\quad=\trho^{\tbhm,\tPiangi}(\tPiangi(j)+\trN_{ij}\tPiangi(i))\tq_{jj}^{-1} \\
\quad=\tbhm(\tPiangi(j),\tPiangi(j))\tbhm(\tPiangi(i),\tPiangi(i))^{N_{ij}}\tq_{jj}^{-1} \\
\quad=\tq_{jj}(\tq_{ij}\tq_{ji})^{\trN_{ij}}\tq_{ii}^{\trN_{ij}^2}\cdot\tq_{ii}^{\trN_{ij}}\cdot\tq_{jj}^{-1} \\
\quad=(\tq_{ij}\tq_{ji})^{\trN_{ij}}\tq_{ii}^{\trN_{ij}(\trN_{ij}+1)}.
\end{array}
\end{equation*}
Note $\tq_{ii}^\ckpchi=1$.
If $j=i$, then
$f(\talpii)=\tq_{ii}^{-2}=f^\prime(\talpii)$ since $\trN_{ii}=-2$.
If $j\ne i$ and $\tq_{ii}^{N_{ij}}\tq_{ij}\tq_{ji}=1$,
then $f(\talpij)=\tq_{ii}^{-(\ckpchi-1)N_{ij}}=\tq_{ii}^{N_{ij}}=f^\prime(\talpij)$.
Assume that $j\ne i$ and $\tq_{ii}^{N_{ij}}\tq_{ij}\tq_{ji}\ne 1$.
Then $\ckpchi=N_{ij}+1$. Hence $\tq_{ii}^{N_{ij}+1}=1$,
Hence $f(\talpij)=f^\prime(\talpij)$.
Thus we conclude $f=f^\prime$, as desired. \hfill $\Box$
\newline\par
Since the statement of \cite[Lemma~6.7]{HY10} did not give 
singular vectors, we reformulate it to Theorem~\ref{theorem:SVth} below
although they are similar.
To give an explicit formula of singular vectors seems to become necessary in
future studies of many areas such as representation theory,
conformal field theory and vertex operator algebras.

\begin{theorem}\label{theorem:SVth}
{\rm{(}}See also {\rm{\cite[{\it{Lemma}}~6.7]{HY10}.)}}
Let $\CardtrRp:=\originalCardtrRp$. Assume $\CardtrRp<\infty$. 
Keep the notation of \eqref{eqn:trRptbhmtPi} and \eqref{eqn:DEFdtEF}.
Let $m\in\fkJ_{1,\CardtrRp}$. 
Let $a_{m^\prime}:=\trhobctal(\dotbeta_{m^\prime})$,
$b_{m^\prime}:=\tbhm(\dotbeta_{m^\prime},\dotbeta_{m^\prime})$
and $\dkpch_{m^\prime}:=\kpch(b_{m^\prime})$
for $m^\prime\in\fkJ_{1,m}$.
Assume $\dkpch_{m^\prime}\ne 0$ for all $m^\prime\in\fkJ_{1,m}$. 
Let $t\in\fkJ_{1,\dkpch_m-1}$.
For $m^\prime$, $m^{\prime\prime}\in\fkJ_{1,m}$, let
$p_{m^\prime,m^{\prime\prime}}:=
\tbhm(\dotbeta_{m^\prime},\dotbeta_{m^{\prime\prime}})
\tbhm(\dotbeta_{m^{\prime\prime}},\dotbeta_{m^\prime})$.
Assume that the following conditions {\rm{(i)}}, {\rm{(ii)}} and {\rm{(iii)}} hold.
\newline\par
{\rm{(i)}} $\trLambda(a_m\trK_{\dotbeta_m}-b_m^t\trL_{\dotbeta_m})=0$. \par
{\rm{(ii)}} $\trLambda(a_{m^\prime}\trK_{\dotbeta_{m^\prime}}
-b_{m^\prime}^{t^\prime}\trL_{\dotbeta_{m^\prime}})\ne 0$
for all ${m^\prime}\in\fkJ_{1,m-1}$
and all ${t^\prime}\in\fkJ_{1,\dkpch_{m^\prime}-1}$. 
 \par
{\rm{(iii)}} $\trLambda(\trK_{\dotbeta_{m^\prime}}-b_{m^\prime}^{t^\prime-1}
(\prod_{m^{\prime\prime}=1}^{m^\prime-1}
p_{m^\prime,m^{\prime\prime}}^{\dkpch_{m^{\prime\prime}}-1})
p_{m^\prime,m}^t
\trL_{\dotbeta_{m^\prime}})\ne 0$
for all ${m^\prime}\in\fkJ_{1,m-1}$
and all ${t^\prime}\in\fkJ_{1,\dkpch_{m^\prime}-1}$. 
\newline\newline
Define $\tv^\prime\in\mclM(\trLambda)_{-t\dotbeta_m}$ by
\begin{equation*}
\tv^\prime:=\dottrE_1^{\dkpch_1-1}\cdots \dottrE_{m-1}^{\dkpch_{m-1}-1}
\dottrF_m^t\dottrF_{m-1}^{\dkpch_{m-1}-1}\cdots \dottrF_1^{\dkpch_1-1}\tv_\trLambda.
\end{equation*} Then we have the following.
\newline\par
{\rm{(1)}} $\tv^\prime\ne 0$, and $\trE_j\tv^\prime= 0$ for all $j\in\fkI$.
\par
{\rm{(2)}} The elements
\begin{equation*}
\begin{array}{l}
\dottrE_{m-1}^{r_{m-1}}\cdots \dottrE_1^{r_1}\dottrF_l^{r_l}\cdots \dottrF_{m+1}^{r_{m+1}}
\dottrF_m^p\dottrF_{m-1}^{\dkpch_{m-1}-1}\cdots \dottrF_1^{\dkpch_1-1}\tv^\prime \\
(\in\mclM(\trLambda)_{-((t+p)\dotbeta_m+(\sum_{y=m+1}^lr_y\dotbeta_y)+
(\sum_{z=1}^{m-1}(\dkpch_z-1-r_z)\dotbeta_z))})
\\
\,\,(\,p\in\fkJ_{0,\dkpch_m-1-t},\,r_x\in\fkJ_{0,\dkpch_x-1}\,
(x\in\fkJ_{1,\CardtrRp}\setminus\{m\})\,)
\end{array}
\end{equation*} form the $\bK$-basis of the left $\trU$-submodule 
$\trU\cdot\tv^\prime$ of $\MtbctaltrLam$.
\end{theorem}
{\it{Proof.}} For $i$, $j\in\fkI$, let $\tq_{ij}:=\tbhm(\talpii,\talpij)$.
For $t\in\bN$ and $i$, $j\in\fkI$, we have
\begin{equation}\label{eqn:FEm}
[\trE_i, \trF_j^t]=\delta_{ij}(t)_{\tq_{ii}}(-\trK_\talpii
+\tq_{ii}^{-t+1}\tL_\talpii)\trF_j^{t-1}.
\end{equation}
If $m=1$, the statement can easily be proved by (i), \eqref{eqn:FEm} and
Theorem~\ref{theorem:LusPBW}.

Assume $m\in\fkJ_{2,\CardtrRp}$.
Let $i:=\hatf(1)\in\fkI$, see \eqref{eqn:defhatf} for $\hatf$.
Then $\dotbeta_1=\talpii$,
$\dottrE_1=\trE_i$, and $\dottrF_1=\trF_i$.
For ${m^\prime}\in\fkJ_{2,m}$
and ${t^\prime}\in\fkJ_{1,\dkpch_{m^\prime}-1}$, we have
\begin{equation}\label{eqn:steppr}
\begin{array}{l}
\trLambdaangi(\trho^{\tbhm,\tPiangi}(\dotbeta_{m^\prime})\trK_{\dotbeta_{m^\prime}}
-b_{m^\prime}^{t^\prime}\trL_{\dotbeta_{m^\prime}}) \\
\quad =\trLambda(\trho^{\tbhm,\tPiangi}(\dotbeta_{m^\prime})
\tbhm(\talpii,\dotbeta_{m^\prime})^{-{c_1}+1}\trK_{\dotbeta_{m^\prime}}-b_{m^\prime}^{t^\prime}
\tbhm(\dotbeta_{m^\prime},\talpii)^{{c_1}-1}\trL_{\dotbeta_{m^\prime}})    \\
\quad\quad\mbox{(by \eqref{eqn:Lmbdai})}\\
\quad =\tbhm(\dotbeta_{m^\prime},\talpii)^{{c_1}-1}\trLambda(a_{m^\prime}\trK_{\dotbeta_{m^\prime}}
-b_{m^\prime}^{t^\prime}\trL_{\dotbeta_{m^\prime}}) 
\quad\mbox{(by \eqref{eqn:tbtbrho})}.
\end{array}
\end{equation} For ${m^\prime}\in\fkJ_{2,m-1}$
and ${t^\prime}\in\fkJ_{1,\dkpch_{m^\prime}-1}$, we have
\begin{equation}\label{eqn:stepprr}
\begin{array}{l}
\trLambdaangi(\trK_{\dotbeta_{m^\prime}}-b_{m^\prime}^{t^\prime-1}
(\prod_{m^{\prime\prime}=2}^{m^\prime-1}
p_{m^\prime,m^{\prime\prime}}^{\dkpch_{m^{\prime\prime}}-1})
p_{m^\prime,m}^t
\trL_{\dotbeta_{m^\prime}}) \\
\quad =\tbhm(\dotbeta_{m^\prime},\talpii)^{-(c_1-1)}\trLambda(\trK_{\dotbeta_{m^\prime}}
-b_{m^\prime}^{t^\prime-1}
(\prod_{m^{\prime\prime}=1}^{m^\prime-1}
p_{m^\prime,m^{\prime\prime}}^{\dkpch_{m^{\prime\prime}}-1})
p_{m^\prime,m}^t
\trL_{\dotbeta_{m^\prime}}) \\
\quad\quad \mbox{(by \eqref{eqn:Lmbdai})}.
\end{array}
\end{equation}
Let
\begin{equation}\label{eqn:defvpp}
\begin{array}{lcl}
\tv^{\prime\prime}&:=&(\dottrE_2^{\dkpch_2-1}\cdots \dottrE_{m-1}^{\dkpch_{m-1}-1}
\dottrF_m^t\dottrF_{m-1}^{\dkpch_{m-1}-1}\cdots\dottrF_2^{\dkpch_2-1})
\dottrF_1^{c_1-1}\tv_\trLambda \\
& & \in\MtbctaltrLam_{-t\dotbeta_m-(c_1-1)\dotbeta_1}.
\end{array}
\end{equation} By \eqref{eqn:defvpp}, we have 
\begin{equation}\label{eqn:vpEikmvpp}
\tv^\prime=\dottrE_1^{c_1-1}\tv^{\prime\prime}.
\end{equation}

Let $\trT_i:=\trTtbtali$.
Let $\ddottrE_x:=(\trT_i)^{-1}(\dottrE_{x+1})$,
($x\in\fkJ_{1,m-2}$) and $\ddottrF_y:=(\trT_i)^{-1}(\dottrF_{y+1})$,
($y\in\fkJ_{1,m-1}$).
By (ii) and Theorem~\ref{theorem:LisoML}, we have $\hTtbtali$
and conclude that
\begin{equation}\label{eqn:prevpp}
\begin{array}{l}
(\hTtbtali)^{-1}(\tv^{\prime\prime}) \\
\quad =\ddottrE_1^{\dkpch_2-1}\cdots \ddottrE_{m-2}^{\dkpch_{m-1}-1}
\ddottrF_{m-1}^t\ddottrF_{m-2}^{\dkpch_{m-1}-1}\cdots\ddottrF_1^{\dkpch_2-1}
\tv_\trLambdaangi \\
\quad\quad \in\mclM^{\tbhm,\tPiangi}(\trLambdaangi)_{-t\dotbeta_m}.
\end{array}
\end{equation} 
By induction and by \eqref{eqn:steppr}, \eqref{eqn:stepprr},
\eqref{eqn:prevpp}, Lemma~\ref{lemma:extlgd} and Theorems~\ref{theorem:Liso}
and \ref{theorem:LisoML},
we have 
\newline\par
$(1^\prime)$ $\tv^{\prime\prime}\ne 0$, 
and $\trT_i(\trE_j)\tv^{\prime\prime}= 0$ ($j\in\fkI$). In particular,
$\trF_i\tv^{\prime\prime}= 0$.
\par
$(2^\prime)$ The elements
$
\dottrE_{m-1}^{r_{m-1}}\cdots \dottrE_2^{r_2}\dottrE_1^{r_1}
\dottrF_l^{r_l}\cdots \dottrF_{m+1}^{r_{m+1}}
\dottrF_m^p\dottrF_{m-1}^{\dkpch_{m-1}-1}\cdots \dottrF_2^{\dkpch_2-1}\tv^{\prime\prime}
$
($p\in\fkJ_{0,\dkpch_m-1-t}$, $r_x\in\fkJ_{0,\dkpch_x-1}$ ($x\in\fkJ_{1,l}\setminus\{m\}$))
form the $\bK$-basis of the left $\trU$-submodule 
$\trU\tv^{\prime\prime}$ of $\MtbctaltrLam$.
\newline\newline
By \eqref{eqn:abcequi}, $\trE_i^{\dkpch_1}=0$. 
Hence by \eqref{eqn:vpEikmvpp},
\begin{equation}\label{eqn:Eivpo}
\trE_i\tv^\prime=0.
\end{equation}
We have
\begin{equation*}
\begin{array}{l}
\dottrF_1^{{\dkpch_1}-1}\tv^\prime \\
\quad = \trF_i^{{\dkpch_1}-1}\trE_i^{{\dkpch_1}-1}\tv^{\prime\prime} \\
\quad = (-1)^{{\dkpch_1}-1}({\dkpch_1}-1)_{\tq_{ii}}!(\prod_{x=1}^{{\dkpch_1}-1}
(-\trK_\talpii+\tq_{ii}^{-x+1}\trL_\talpii))\tv^{\prime\prime} \\
\quad\quad\quad\mbox{(by \eqref{eqn:FEm} and $(1^\prime)$)} \\
\quad = ({\dkpch_1}-1)_{\tq_{ii}}!(\prod_{x=1}^{{\dkpch_1}-1}\trLambda(
\tbhm(\talpii,\dotbeta_m)^{-t}\tq_{ii}^{-{\dkpch_1}+1}\trK_\talpii-\tq_{ii}^{-x+{\dkpch_1}}
\tbhm(\dotbeta_m,\talpii)^t\trL_\talpii))\tv^{\prime\prime} \\
\quad\quad\quad\mbox{(by \eqref{eqn:defvpp})} \\
\quad = ({\dkpch_1}-1)_{\tq_{ii}}!\tq_{ii}^{{\dkpch_1}-1}(\prod_{x=1}^{{\dkpch_1}-1}\trLambda(
\tbhm(\talpii,\dotbeta_m)^{-t}\trK_\talpii-\tq_{ii}^{-x-1}
\tbhm(\dotbeta_m,\talpii)^t\trL_\talpii))\tv^{\prime\prime} \\
\quad = ({\dkpch_1}-1)_{\tq_{ii}}!\tq_{ii}^{{\dkpch_1}-1}\tbhm(\talpii,\dotbeta_m)^{-t({\dkpch_1}-1)}(\prod_{x=1}^{{\dkpch_1}-1}(\trLambda(
\trK_\talpii-b_1^{x-1}
p_{1m}^t\trL_\talpii))\tv^{\prime\prime} \\
\quad \in  \bKt\tv^{\prime\prime}\quad\mbox{(by (iii))}.
\end{array}
\end{equation*} Hence $\tv^\prime\ne 0$, and, by $(2^\prime)$, we have $(2)$.

We show $(1)$. 
Let $j\in\fkI\setminus\{i\}$. Let $\trN_{ij}:=\trNtbhmtPi_{i,j}$.
Let 
\begin{equation*}
X_z:=\sum_{t=0}^z(-\tq_{ij})^t\tq_{ii}^{\frac
{t(t-1)} {2}}
{{z}\choose{t}}_{\tq_{ii}}\trE_i^{z-t}\trE_j\trE_i^t
\quad\quad(z\in\fkJ_{-1,{\dkpch_1}}),
\end{equation*} where $X_{-1}=0$.
Then $\trE_iX_z-\tq_{ii}^z\tq_{ij}X_z\trE_i=X_{z+1}$
$(z\in\fkJ_{0,{\dkpch_1}-1})$, from which we have $[X_z,\trF_i]=(z)_{\tq_{ii}}
(1-\tq_{ii}^{z-1}\tq_{ij}\tq_{ji})\trL_\talpii X_{z-1}$
$(z\in\fkJ_{0,\dkpch_1})$.
Note that $N_{ij}\leq {\dkpch_1}-1$ and $X_{\trN_{ij}}=\trT_i(\trE_j)$.
Hence $X_{z_1}\tv^{\prime\prime}=0$  ($z_1\in\fkJ_{0,\trN_{ij}}$)
since $X_{\trN_{ij}}\tv^{\prime\prime}=0$ and $\trF_i\tv^{\prime\prime}=0$
(see $(1^\prime)$).
We have $X_{z_2}=0$ ($z_2\in\fkJ_{\trN_{ij}+1,\dkpch_1}$) since 
$X_{\trN_{ij}+1}=\trE_i\trT_i(\trE_j)-\tq_{ii}^{\trN_{ij}}\tq_{ij}\trT_i(\trE_j)\trE_i
=\trK_\talpii\trT_i(\trF_i\trE_j)
-\tq_{ii}^{\trN_{ij}}\tq_{ij}\trT_i(\trE_j\trF_i)\cdot\tq_{ii}\trK_\talpii =
\trK_\talpii\trT_i([\trF_i,\trE_j])=0$.
Hence $\trE_j\trE_i^r\tv^{\prime\prime}=0$  ($r\in\fkJ_{0,{\dkpch_1}-1}$)
since $X_{z_3}\tv^{\prime\prime}=0$  ($z_3\in\fkJ_{0,{\dkpch_1}-1}$).
By \eqref{eqn:vpEikmvpp}, $\trE_j\tv^\prime=0$, as desired.
This completes the proof. \hfill $\Box$

\subsection{Shapovalov determinats} \label{subsection:ShDet}
Recall the convention \eqref{eqn:abbrU}.
Define the $\bK$-linear map $\trShbctal:\trU\to\trUo$
by 
\begin{equation*}
\begin{array}{l}
\trShbctal(X)= \\ \quad 
\left\{\begin{array}{ll}
X & \quad\mbox{if $X\in\trUo$,} \\
0 & \quad\mbox{if $X\in\trUm_{-\tmu}\trUo\trUp_\tlambda$
with $\tlambda$, $\tmu\in\tfkAPip$ such that $(\tlambda,\tmu)\ne (0,0)$.}
\end{array}\right.
\end{array}
\end{equation*}
Let $\tlambda\in\tfkAPip$. By \eqref{eqn:bprtoftB}, we see that
\begin{equation} \label{eqn:trShXY}
\begin{array}{l}
\trShbctal(XY)-(\trtvttbhmtPi(X,\HopfS(Y))
\trK_\tlambda+\trtvttbhmtPi(X,Y)\trL_\tlambda) \\
\quad\displaystyle{\in\bigoplus_{{{\tmu,\tnu\in\tfkAPip\setminus\{0\},}\atop{\tmu+\tnu=\tlambda}}}
\bK\trK_\tmu\trL_{\tnu}} \\
\quad\quad(X\in\trUp_\tlambda,Y\in\trUm_{-\tlambda}).
\end{array}
\end{equation} Let $\mtlambda:=\dim\trUp_\tlambda$.
Let $\hatX:\fkJ_{1,\mtlambda}\to\trUp_\tlambda(\in\trUo)$ and $\hatY:\fkJ_{1,\mtlambda}\to\trUm_{-\tlambda}$
be maps.
Let
\begin{equation} \label{eqn:pretrShXYMt}
\dhXhY:=\det[\trtvttbhmtPi(\hatX(x),\hatY(y))]_{1\leq x,y\leq\mtlambda}\in\bK,
\end{equation} and
\begin{equation} \label{eqn:trShXYMt}
\mclSbctaltl[\hatX,\hatY]:=[\trShbctal(\hatX(x)\hatY(y))]_{1\leq x,y\leq\mtlambda}\in\Mat(\mtlambda,\trUo).
\end{equation} 
By \eqref{eqn:trShXY}, $\dhXhY$ is the coefficient of $\trL_{\mtlambda\tlambda}$ 
of $\det\mclSbctaltl[\hatX,\hatY]$.
Hence, by \eqref{eqn:dimUmlUpld}, we see that
\begin{equation} \label{eqn:dettrShXYMt}
\begin{array}{l}
\mbox{$\dhXhY\ne 0$ and $\det\mclSbctaltl[\hatX,\hatY]\ne 0$
if $\{\hatX(x)|x\in\fkJ_{1,\mtlambda}\}$ is a $\bK$-basis of $\trUp_\tlambda$} \\
\mbox{and $\{\hatY(y)|y\in\fkJ_{1,\mtlambda}\}$ is a $\bK$-basis of $\trUm_{-\tlambda}$.}
\end{array}
\end{equation}

Let $\trLambda\in\rmChtbhmchecktal$. 
Let $\mclM(\trLambda)$ mean $\MtbctaltrLam$.
For $\tlambda\in\tfkAPi$, let
$\mclM(\trLambda)_{-\tlambda}:=\trUm_{-\tlambda}\tv_\trLambda$.
Then $\mclM(\trLambda)=\oplus_{\tlambda\in\tfkAPip}\mclM(\trLambda)_{-\tlambda}$
as a $\bK$-linear space.
Let $\NtbctaltrLam$ be the largest $\trU$-submodule 
$\mclN^\prime$ of $\mclM(\trLambda)$ such that 
$\mclN^\prime\subset\oplus_{\tlambda\in\tfkAPip\setminus\{0\}}\mclM(\trLambda)_{-\tlambda}$.
Let $\mclN(\trLambda)$ mean $\NtbctaltrLam$.
For $\tlambda\in\tfkAPi$, let
$\mclN(\trLambda)_{-\tlambda}:=\mclN(\trLambda)\cap\mclM(\trLambda)_{-\tlambda}$.
Then $\mclN(\trLambda)=\oplus_{\tlambda\in\tfkAPip\setminus\{0\}}\mclN(\trLambda)_{-\tlambda}$
as a $\bK$-linear space.
Let 
$\LtbctaltrLam:=\mclM(\trLambda)/\mclN(\trLambda)$.
Then $\LtbctaltrLam$ is a non-zero irreducible $\trU$-module.
Let $\mclL(\trLambda)$ mean $\LtbctaltrLam$.
Let $v_\trLambda:=\tv_\trLambda+\mclN(\trLambda)\in\mclL(\trLambda)$.
For $\tlambda\in\tfkAPi$, let
$\mclL(\trLambda)_{-\tlambda}:=\trUm_{-\tlambda}v_\trLambda$.
Then $\mclL(\trLambda)=\oplus_{\tlambda\in\tfkAPip}\mclL
(\trLambda)_{-\tlambda}$
as a $\bK$-linear space.
We have  $\dim\mclL(\trLambda)_{-\tlambda}
=\dim\mclM(\trLambda)_{-\tlambda}-\dim\mclN(\trLambda)_{-\tlambda}$
($\tlambda\in\tfkAPi$).

By a standard argument (cf. \cite[Section~6]{HY10}), we have
\begin{lemma} \label{lemma:EssShapo} Recall the convention \eqref{eqn:abbrU}.
Let $\tlambda\in\tfkAPip$. Let $\mtlambda:=\dim\trUp_\tlambda$.
Let $\{X_x|x\in\fkJ_{1,\mtlambda}\}$ be a $\bK$-basis of $\trUp_\tlambda$.
Let $\{Y_y|y\in\fkJ_{1,\mtlambda}\}$ be a $\bK$-basis of $\trUm_{-\tlambda}$
{\rm{(}}see also \eqref{eqn:dimUmlUpl}{\rm{)}}.
Let $\trLambda\in\rmChtbhmchecktal$.
Then
\begin{equation}\label{eqn:EssShpEq}
\dim\LtbctaltrLam_{-\tlambda}=\MatRank([\trLambda(\trShbctal(X_xY_y))]_{1\leq x,y\leq\mtlambda}).
\end{equation}
Moreover
\begin{equation}\label{eqn:EssShpKr}
\begin{array}{l}
\{\,Y\in\trUm_{-\tlambda}\,|\,Y\tv_\trLambda\in\NtbctaltrLam_{-\tlambda}\,\} \\
\quad =\{\,Y\in\trUm_{-\tlambda}\,|\,\trLambda(\trShbctal(XY))=0\,\,\mbox{for all}\,\,X\in\trUp_\tlambda\,\}.
\end{array}
\end{equation}

\end{lemma}

Let $\trRpmapbctal$ be the set of maps
$f:\trRtbhmtPi\to\bZgeqo$ satisfying the condition that
$(f(\tal))_{\tbhm(\tal,\tal)}!\ne 0$
for all $\tal\in\trRtbhmtPi$. For $\tlambda\in\tfkAPip$,
let
\begin{equation}\label{eqn:DFtrRPmp}
\trRpmapbctal_\tlambda=\{\,f\in\trRpmapbctal\,|\,
\sum_{\tal\in\trRptbhmtPi}f(\tal)\tal=\tlambda\,\}.
\end{equation} From Theorem~\ref{theorem:LusPBW},
it follows that
\begin{equation}\label{eqn:preLPBW}
\CardtrRptbtP<\infty\,\,\Rightarrow\,\,
|\trRpmapbctal_\tlambda|=\mtlambda
\quad(\tlambda\in\tfkAPip),
\end{equation} where $\mtlambda:=\dim\trUp_\tlambda=\dim\trUm_{-\tlambda}$.

For $\tlambda\in\tfkAPip$, $\tal\in\trRptbhmtPi$ and $t\in\bZgeqo$, let
\begin{equation*}
\trRpmapbctal_\tlambda(\tal;t):=
\{\,f\in\trRpmapbctal_\tlambda\,|\,f(\tal)\geq t\,\}.
\end{equation*}

Let $\CardtrRp:=\CardtrRptbtP$. Assume $\CardtrRp<\infty$.
As for the statement of Theorem~\ref{theorem:LusPBW},
we see that 
$|\trRpmapbctal_\tlambda(\dotbeta_x;t)|=\dim(\trUm_\tlambda\cap\trUm \dottrF_x^t)$
for $\tlambda\in\tfkAPip$, $x\in\fkJ_{1,\CardtrRp}$ and $t\in\bN$.

By \eqref{eqn:trShXY} and \cite[Theorem~7.3]{HY10}, we have

\begin{theorem}{\rm{(\cite[Theorem~7.3]{HY10})}}\label{theorem:Shapo} 
Assume that $\CardtrRptbtP<\infty$.
Assume that $\tbhm(\tal,\tal)\ne 1$ for all $\tal\in\trRptbhmtPi$.
Let $\tlambda\in\tfkAPip$. 
Let $\mtlambda:=\dim\trUp_\tlambda$.
Let $\dhXhY$ and $\mclSbctaltl[\hatX,\hatY]\in\Mat(m_\tlambda,\trUo)$
be as in \eqref{eqn:pretrShXYMt}, \eqref{eqn:trShXYMt} respectively.
Let $v_\tal:=\rmMax\{\,t\in\bZgeqo\,|\,\trRpmapbctal_\tlambda(\tal;t)\ne\emptyset\,\}$
for $\tal\in\trRptbhmtPi$.
Then, as an element of $\trUo$, we have
\begin{equation}\label{eqn:ShapoEq}
\begin{array}{l}
\det\mclSbctaltl[\hatX,\hatY] \\
\quad\displaystyle{=\dhXhY\cdot
\prod_{\tal\in\trRptbhmtPi}\prod_{t_\tal=1}^{v_\tal}
(-\trhobctal(\tal)\tbhm(\tal,\tal)^{-t_\tal}
\trK_\tal+\trL_\tal)^{|\trRpmapbctal_\tlambda(\tal;t_\tal)|}.}
\end{array}
\end{equation}
{\rm{(}}As for $\dhXhY$, recall \eqref{eqn:dettrShXYMt}.{\rm{)}}
\end{theorem}

We also give a proof of Theorem~\ref{theorem:Shapo}
after Theorem~\ref{theorem:keythma} below.

\section{Applications of Sections~{\rm{\ref{section:App}}}
and {\rm{\ref{section:AppII}}} to Shapovalov factors}
\label{section:ShFac}

\subsection{Nichols locally closed spaces of bi-homomorphisms}\label{subsection:LCSbh}

The argument here has originally been given in \cite[Section~7]{HY10}.
\begin{equation}\label{eqn:assNlcs}
\begin{array}{l}
\mbox{In Section~\ref{section:ShFac}, assume that $\CardtrRptbtP<\infty$} \\
\mbox{and that $\tbhm(\tal,\tal)\ne 1$ for all $\tal\in\trRptbhmtPi$.}
\end{array}
\end{equation}

Let $\tbX$ be the set of all bi-homomorphisms $\tbhm^\prime:\trtfkAPi\times\trtfkAPi\to\bKt$.
Recall $\CardfkIN=\CardfkI$.

\begin{lemma}\label{lemma:locltbhm}
Let $\tbhm^\prime\in\tbX$.
Then $\trRp(\tbhm^\prime,\checktal)=\trRptbhmtPi$
if and only if
\begin{equation}\label{eqn:lcltcnd}
\trN^{\tbhm^\prime,\checktal^\prime}_{i,j}
=\trN^{\tbhm,\checktal^\prime}_{i,j}\quad
(\checktal^\prime\in\checktBtbhmchPictal,\,\,i,\,j\in\fkI).
\end{equation}
\end{lemma}
{\it{Proof.}} We can easily see that `only-if'-part follows from
the first equation of \eqref{eqn:defNij}.
We prove `if'-part. Assume that  
\eqref{eqn:lcltcnd} holds.
Let $\CardtrRp:=\CardtrRptbtP$.
Recall $\hatf$ from \eqref{eqn:defhatf}.
By \eqref{eqn:trRptbhmtPi} and \eqref{eqn:lcltcnd}, we have 
$\trRptbhmtPi\subset\trRp(\tbhm^\prime,\checktal)$.
By \eqref{eqn:prvprtd}, for $k\in\fkJ_{0,\CardtrRp}$
there exists a subset $X_k$ of $\trRp(\tbhm^\prime,\checktal)$
with ${\frac {k-|X_k|} 2}\in\fkJ_{0,{\frac k 2}}$
and $\trRp(\tbhm^\prime,\checktal_{\hatf,k})=
(-X_k)\cup(\trRp(\tbhm^\prime,\checktal)\setminus X_k)$.
Since $\checktal_{\hatf,\CardtrRp}(\fkI)=-\tPi$, we have
$\trRp(\tbhm^\prime,\checktal_{\hatf,\CardtrRp})=-\trRp(\tbhm^\prime,\checktal)$,
so $X_\CardtrRp=\trRp(\tbhm^\prime,\checktal)$, whence
$\trRp(\tbhm^\prime,\checktal)\leq\CardtrRp$.  
Hence $\trRptbhmtPi=\trRp(\tbhm^\prime,\checktal)$,
as desired.
\hfill $\Box$
\newline\par
For $d\in\bKtprim$ and $p=\sum_{i,j\in\fkI}c_{i,j}\talpii\otimes\talpij
\in\tfkAPi\otimes_\bZ\tfkAPi$ with $c_{i,j}\in\bZ$, 
let 
\begin{equation*}
V(d;p):=\{\,(\,z_{i,j}\,|\,i,j\in\fkI\,)\in\bKtll\,|\,
d=\prod_{i,j\in\fkI}z_{i,j}^{c_{i,j}}\,\},
\end{equation*}
and $D(d;p):=\bKtll\setminus V(d;p)$.
Then $V(d;p)$ (resp. $D(d;p)$) are a Nichols closed (resp. open)
subsets of $\bKtll$.

Define the Nichols closed subsets $\bV^\tbhm_t$ of $\bKtll$ ($t\in\fkJ_{1,2}$) by
\begin{equation*}
\begin{array}{l} 
\bV^\tbhm_1:=\displaystyle{\bigcap_{\checktal^\prime\in\checktBtbhmchPictal}\,\bigcap_{i,\,j\in\fkI,\,i\ne j}}
\bigl(
V(1;(N^{\tbhm,\checktal^\prime}_{i,j}+1)\checktal^\prime(i)\otimes\checktal^\prime(i))
\\
\quad\quad\quad\quad\quad\quad\quad\quad
\cup V(1;N^{\tbhm,\checktal^\prime}_{i,j}
\checktal^\prime(i)\otimes\checktal^\prime(i)
+\checktal^\prime(i)\otimes\checktal^\prime(j)
+\checktal^\prime(j)\otimes\checktal^\prime(i))\bigr), 
\end{array}
\end{equation*} and 
\begin{equation*}
\bV^\tbhm_2:=\bigcap_{{\tal\in\trRptbhmtPi,}\atop{\kpch(\tbhm(\tal,\tal))\geq 2}}V(1;\kpch(\tbhm(\tal,\tal))\tal\otimes\tal).
\end{equation*}
Define the Nichols open subsets $\bD^\tbhm_t$ ($t\in\fkJ_{1,3}$) of $\bKtll$ by
\begin{equation*}
\bD^\tbhm_1:=\cap_{\tal\in\trRptbhmtPi}D(1;\tal\otimes\tal),
\end{equation*} 

\begin{equation*}
\begin{array}{l}
\bD^\tbhm_2:=\displaystyle{\bigcap_{\checktal^\prime\in\checktBtbhmchPictal}\,\bigcap_{i,\,j\in\fkI,\,i\ne j}
\,\bigcap_{s=1}^{N^{\tbhm,\checktal^\prime}_{i,j}}
\bigl(
D(1;s\checktal^\prime(i)\otimes\checktal^\prime(i))}
\\
\quad\quad\quad\quad\quad\quad\quad\quad
\cap D(1;(s-1)\checktal^\prime(i)\otimes\checktal^\prime(i)
+\checktal^\prime(i)\otimes\checktal^\prime(j)
+\checktal^\prime(j)\otimes\checktal^\prime(i))\bigr), 
\end{array}
\end{equation*} and
\begin{equation*}
\bD^\tbhm_3:=\bigcap_{{\tal\in\trRptbhmtPi},\atop{\kpch(\tbhm(\tal,\tal))\geq 2}}\,
\bigcap_{s=1}^{\kpch(\tbhm(\tal,\tal))-1}D(1;s\tal\otimes\tal).
\end{equation*} 

Define the bijection $\trxi:\tbX\to\bKtll$
by $\trxi(\tbhm^\prime):=(\,\tbhm^\prime(\talpii,\talpij)\,|\,i,j\in\fkI\,)$.
Let $\tbXtbhmone:=\trxi^{-1}(\bV^\tbhm_1\cap\bD^\tbhm_1\cap\bD^\tbhm_2)$,
and $\tbXtbhmtwo:=\trxi^{-1}(\bV^\tbhm_1\cap\bV^\tbhm_2\cap
\bD^\tbhm_1\cap\bD^\tbhm_2\cap\bD^\tbhm_3)$.
By Lemma~\ref{lemma:locltbhm}, we have
\begin{equation}\label{eqn:spcofR}
\tbXtbhmone=\{\,\tbhm^\prime\in\tbX\,|\,\tbhm^\prime(\tal,\tal)\ne 1\,
(\tal\in\trRptbhmtPi),\,\,\trRp(\tbhm^\prime,\tPi)=\trRptbhmtPi\,\},
\end{equation} and 
\begin{equation}\label{eqn:spcofRht}
\begin{array}{lcl}
\tbXtbhmtwo & = & \{\,\tbhm^\prime\in\tbXtbhmone\,|\,\kpch(\tbhm^\prime(\tal,\tal))
=\kpch(\tbhm(\tal,\tal)) \\
& & \quad\quad\quad\quad\quad 
\mbox{for all
$\tal\in\trRptbhmtPi$ with $\kpch(\tbhm(\tal,\tal))\geq 2$}\,\}.
\end{array}
\end{equation}
Let
\begin{equation*}
\tbXtprim:=\{\,\tbhm^\prime\in\tbX\,|\,
\tbhm^\prime(\talpii,\talpij)\in\bKtprim\,
(i,j\in\fkI)\,\}.
\end{equation*} 
\begin{remark}
By Theorem~\ref{theorem:basicc}, for $t\in\fkJ_{1,2}$,
$\trxi(\tbXtprim\cap\tbX^\tbhm_t)$ is dense in $\trxi(\tbX^\tbhm_t)$
under the Zariski topology on $\bKtll$. 
\end{remark}

Let $\tlambda\in\tfkAPip$.
Define the Nichols open subset $\bD^{\tbhm;\tlambda}_4$  of $\bKtll$ by 
\begin{equation*} 
\bD^{\tbhm;\tlambda}_4:=\bigcap_{{\tal\in\trRptbhmtPi,}\atop{\kpch(\tbhm(\tal,\tal))=0}}\,
\bigcap_{{t\in\bN,}\atop{r(\tal,t)\geq 1}}D(1;t\tal\otimes\tal),
\end{equation*} where
$r(\tal,t):=|\trRpmapbctal_\tlambda(\tal;t)|$.
Let $\tbXtbhmtlambda:=\trxi^{-1}(\bD^{\tbhm;\tlambda}_4\cap\trxi(\tbXtbhmtwo))$.
By \eqref{eqn:preLPBW}, \eqref{eqn:spcofR}, \eqref{eqn:spcofRht}
and Theorem~\ref{theorem:LusPBW}, we have
\begin{equation}\label{eqn:spcofRb}
\tbXtbhmtlambda =\{\,\tbhm^\prime\in\tbXtbhmtwo\,|\,|\trRpmap^{\tbhm^\prime,\checktal}_\tlambda|
=|\trRpmapbctal_\tlambda|\,\}.
\end{equation} 

\begin{remark}
By Theorem~\ref{theorem:basicc}, 
$\trxi(\tbXtprim\cap\tbXtbhmtlambda)$ is dense in $\trxi(\tbXtbhmtlambda)$
under the Zariski topology on $\bKtll$. 
\end{remark}

\subsection{Nichols locally closed spaces of characters}\label{subsection:LCSch}
Let $\tbUo$ be the $\bK$-algebra defined by the
$\trK_\tlambda$, $\trL_\tlambda$ ($\tlambda\in\tfkAPi$)
and the relations
$\trK_0=\trL_0=1$,
$\trK_\tlambda \trK_\tmu=\tK_{\tlambda+\tmu}$,
$\trL_\tlambda \trL_\tmu=\tL_{\tlambda+\tmu}$,
$\trK_\tlambda \trL_\tmu=\tL_\tmu \trK_\tlambda$
($\tlambda$, $\tmu\in\tfkAPi$).
Then, for $\tbhm^\prime\in\tbX$ and $\checktal^\prime\in\checktBtbhmchPi$, we have the $\bK$-algebra isomorphism 
$\trzeta^{\tbhm^\prime,\checktal^\prime}:\tbUo\to\trUo(\tbhm^\prime,\checktal^\prime)$
defined by $\trzeta^{\tbhm^\prime,\checktal^\prime}(\trK_\tlambda\trL_\tmu):=\trK_\tlambda\trL_\tmu$
($\tlambda$, $\tmu\in\tfkAPi$).

Define the bijection $\treta:\ChtbUo\to\bKttwl\cong\bKtkC\times\bKtkC$ by 
\begin{equation}\label{eqn:dftreta}
\treta(\trLambda):=(\,(\,\trLambda(\tK_\talpii)\,|\,i\in\fkI\,),\,
(\,\trLambda(\tL_\talpij)\,|\,j\in\fkI\,)\,).
\end{equation}

For $d\in\bKtprim$, and $\lambda=\sum_{i\in\fkI}{\hat c}_i\talpii$,
$\mu=\sum_{i\in\fkI}
{\check c}_i\talpii
\in\tfkAPi$ with ${\hat c}_i$, ${\check c}_i\in\bZ$,
define the Nichols closed subset $V(d;\lambda,\mu)$ of $\bKttwl$ by
\begin{equation*}
\begin{array}{l}
V(d;\lambda,\mu) \\
\quad:=\{\,
(\,(\,{\hat z}_i\,|\,i\in\fkI\,),\,
(\,{\check z}_j\,|\,j\in\fkI\,)\,)
\in\bKttwl \,|\,
d=
(\prod_{{\hat i}\in\fkI}{\hat z}_i^{{\hat c}_i})
(\prod_{{\check i}\in\fkI}{\check z}_j^{{\check c}_j})\,\}.
\end{array}
\end{equation*} 

For $\tlambda\in\tfkAPi$ and $d\in\bKtprim$, let $\bG(\tlambda;d):=\{\,\trLambda\in\ChtbUo\,|\,\trLambda(\tK_\tlambda-d\tL_\tlambda)=0\,\}$,
so $\bG(\tlambda;d)=\treta^{-1}(V(d;\tlambda,-\tlambda))$.

\begin{lemma}\label{lemma:NEMlem}
Let $x\in\bN$, and let
$\tbeta_y\in\trRptbhmtPi$
and $d_y\in\bKtprim$ $(y\in\fkJ_{1,x})$. 
Assume that $(\tbeta_y,d_y)\ne(\tbeta_1,d_1)$
for $y\in\fkJ_{2,x}$.
Then
$\treta(\bG(\tbeta_1;d_1)
\cap(\bigcap_{y=2}^x(\ChtbUo\setminus\bG(\tbeta_y;d_y))))$
is a dense subset of $\treta(\bG(\tbeta_1;d_1))$
under the {\rm{(}}relative{\rm{)}} Zariski topology.
\end{lemma}
{\it{Proof.}} By Lemmas~\ref{lemma:tRainfty}
and \ref{lemma:trRinfty},
there exists $\checktal^\prime\in\checktBtbhmchPi$
and $j\in\fkI$ such that
$\tbeta_1=\checktal^\prime(j)$.
By \eqref{eqn:simplpr},
$\tal\in(\oplus_{i\in\fkI}\bZ\checktal^\prime(i))\setminus\bZ\checktal^\prime(j)$
for $\tal\in\trRptbhmtPi\setminus\{\tbeta_1\}$.
Then we can easily see the statement
by \eqref{eqn:rprecrs} (for $r=1$), 
and Lemma~\ref{lemma:fncvr}.
\hfill $\Box$
\newline\par
For $d\in\bKtprim$, $p=\sum_{i=1}^l\sum_{j=1}^lc_{i,j}\talpii\otimes\talpij
\in\tfkAPi\otimes_\bZ\tfkAPi$ with $c_{i,j}\in\bZ$, 
and $\lambda=\sum_{i\in\fkI}{\hat c}_i\talpii$,
$\mu=\sum_{j\in\fkI}
{\check c}_j\talpij
\in\tfkAPi$ with ${\hat c}_i$, ${\check c}_j\in\bZ$,
define the Nichols closed subset $V(d;p,\lambda,\mu)$ of $\bKtlltwl\cong\bKtll\times\bKttwl$ by
\begin{equation*}
\begin{array}{l}
V(d;p,\lambda,\mu) \\
\quad:=\{\,(\,(\,z_{i,j}\,|\,i,j\in\fkI\,),\,
(\,(\,{\hat z}_{{\hat i}}\,|\,{\hat i}\in\fkI\,),\,
(\,{\check z}_{{\check i}}\,|\,{\check i}\in\fkI\,)\,)\,)
\in\bKtlltwl\\
\quad\quad\quad\,|\,
d=(\prod_{i,j\in\fkI}z_{i,j}^{c_{i,j}})
(\prod_{{\hat i}\in\fkI}{\hat z}_{{\hat i}}^{{\hat c}_{{\hat i}}})
(\prod_{{\check i}\in\fkI}{\check z}_{\check i}^{{\check c}_{{\check i}}})\,\}.
\end{array}
\end{equation*} 
For $t\in\bZ$ and $\tlambda=\sum_{i\in\fkI}m_i\talpii\in\tfkAPi$ with $m_i\in\bZ$, let
\begin{equation*}
\bY(\tlambda;t):=(\trxi^{-1}\times\treta^{-1})(V(1;(\sum_{i\in\fkI}m_i\talpii\otimes\talpii)-t\tlambda\otimes\tlambda,\tlambda,-\tlambda)).
\end{equation*}
We can easily see that
\begin{equation}\label{eqn:prYtbt}
\bY(\tlambda;t) 
=\{\,(\tbhm^\prime,\trLambda)\in\tbX\times\ChtbUo\,|\, 
\trLambda(\trho^{\tbhm^\prime,\checktal}
(\tlambda)\tK_\tlambda-\tbhm^\prime(\tlambda,\tlambda)^t\tL_\tlambda)=0\,\}.
\end{equation} 
For $\tbhm^\prime\in\tbX$, $\tlambda\in\tfkAPi$
and $t\in\bZ$, let $\bY^{\tbhm^\prime}(\tlambda;t)
:=\{\,\trLambda\in\ChtbUo\,|\,(\tbhm^\prime,\trLambda)
\in\bY(\tlambda;t)\,\}$;
note that $\bY^{\tbhm^\prime}(\tlambda;t)=\bG(\tlambda;{\frac {\tbhm^\prime(\tlambda,\tlambda)^t} {\trho^{\tbhm^\prime,\checktal}
(\tlambda)}})$ if $\tbhm^\prime\in\tbXtprim$ and $\tlambda\in\trRp(\tbhm^\prime,\checktal)$.

\subsection{Rank estimation of Shapovalov matrices}
\label{subsection:prNdegShap}

For finite sequences $\bari=(i_1,\ldots,i_x)$, $\barj=(j_1,\ldots,j_x)\in\fkI^x$ in
$\fkI$
with $x\in\bN$ such that 
$\sum_{y=1}^x\checktal(i_y)=\sum_{y=1}^x\checktal(j_y)$,
define the map $\tromega_{{\bar i},{\bar j}}:\tbX\times\ChtbUo\to\bKt$ by
\begin{equation*}
\tromega_{\bari,\barj}(\tbhm^\prime,\trLambda):=
(\trLambda\circ(\trzeta^{\tbhm^\prime,\checktal})^{-1}\circ\trSh^{\tbhm^\prime,\checktal})(\trE_{i_1}\cdots\trE_{i_x}\trF_{j_1}\cdots\trF_{j_x}).
\end{equation*}
Then $\tromega_{\bari,\barj}\circ(\trxi^{-1}\times\treta^{-1})$ is a regular map 
on the affine variety $\bKtlltwl$.  

For $\tlambda=\sum_{i\in\fkI}x_i\talpii\in\tfkAPip\setminus\{0\}$ with $x_i\in\bZgeqo$,
let $\fkI^{(\tlambda)}$ be the set of sequences $(i_1,\ldots,i_x)$,
with 
$\sum_{y=1}^x\checktal(i_y)=\tlambda$,
where we let $x:=\sum_{i\in\fkI}x_i\in\bN$. 

\begin{theorem}\label{theorem:keythma} Recall \eqref{eqn:assNlcs}.
Let $\tlambda\in\tfkAPip\setminus\{0\}$. 
Let $\mtlambda:=\dim\trUptbhmtPi_\tlambda$.
Let $\tal\in\trRptbhmtPi$
and $t\in\bN$.
Let $r:=|\trRpmapbctal_\tlambda(\tal;t)|(\in\bZgeqo)$,
and assume $r\geq 1$.
Let $\trLambda\in\bY^\tbhm(\tal;t)$.
Let $\bari^{(x)}$, $\barj^{(x)}\in\fkI^{(\tlambda)}$
$(x\in\fkJ_{1,\mtlambda-r+1})$. 
Then
\begin{equation*}
\det([\tromega_{\bari^{(x)},\barj^{(y)}}(\tbhm,\trLambda)]_{1\leq x,y\leq\mtlambda-r+1})=0.
\end{equation*}
In particular, 
letting $X$ be the element of $\Mat(\mtlambda,\bK)$ obtained from 
$\mclSbctaltl[\hatX,\hatY]$ of Theorem~{\rm{\ref{theorem:Shapo}}} by sending its
components by $\trLambda\circ(\trzeta^{\tbhm,\checktal})^{-1}$,
\begin{equation}\label{eqn:keythmaneq}
\MatRank(X)\leq\mtlambda-r.
\end{equation}

\end{theorem}
{\it{Proof.}} Define the map $f:\tbX\times\ChtbUo\to\bKt$
by 
\begin{equation*}
f(\tbhm^\prime,\trLambda^\prime):=
\det([\tromega_{\bari^{(x)},\barj^{(y)}}(\tbhm^\prime,\trLambda^\prime)]_{1\leq x,y\leq\mtlambda-r+1}).
\end{equation*}
Then $f\circ(\trxi^{-1}\times\treta^{-1})$ is a regular map on $\bKtlltwl$.

Let 
$\tbhm_1\in\tbXtbhmtlambda\cap\tbXtprim$. Using 
\eqref{eqn:spcofRb}, \eqref{eqn:prYtbt},
Theorem~\ref{theorem:SVth}.
and Lemma~\ref{lemma:NEMlem},
we have $f(\tbhm_1,\trLambda_1)=0$
for $\trLambda_1\in\bY^{\tbhm_1}(\tal;t)$.
In particular, $f(\tbhm_1,\trLambda_2)=0$
for $\trLambda_2\in\bY^{\tbhm_1}(\tal;t)\cap\treta^{-1}(\bKkkgqi)$.
By Theorem~\ref{theorem:basicc}, we have
$f(\tbhm_2,\trLambda_3)=0$
for $\tbhm_2\in\tbXtbhmtlambda$ and $\trLambda_3\in\bY^{\tbhm_2}(\tal;t)$.
This completes the proof, since $\tbhm\in\tbXtbhmtlambda$.
\hfill $\Box$.
\newline\par
{\it{Proof of Theorem~{\rm{\ref{theorem:Shapo}}}.}}
By Theorem~\ref{theorem:LusPBW}, we have
\begin{equation*}
\mtlambda=\sum_{\tal\in\trRptbhmtPi}\sum_{t_\tal=1}^{v_\tal}|\trRpmapbctal_\tlambda(\tal;t_\tal)|.
\end{equation*} 
Then \eqref{eqn:ShapoEq} follows from \eqref{eqn:trShXY}, \eqref{eqn:dettrShXYMt}, \eqref{eqn:keythmaneq},
Lemma~\ref{lemma:aplmt} and \eqref{eqn:irrKL} below. \hfill $\Box$

\subsection{Non-degeneracy of Shapovalov factors}
\label{subsection:NdegShap}
Recall $\trU=\trUtbhmtPi$, $\trUo=\trUotbhmtPi$, $\trUp=\trUptbhmtPi$ and $\trUm=\trUmtbhmtPi$ from \eqref{eqn:abbrU}.  
Note that $\trUo$ is a unique factorization domain
whose invertible elements are $\trK_\tlambda\trL_\tmu$
($\tlambda$, $\tmu\in\tfkAPi$).
By Lemma~\ref{lemma:BtB} and \ref{lemma:trRinfty}, for $\tal\in\trRptbhmtPi$,
there exists $\tPi^\prime\in\tBtbhmchPi$ such that $\tal\in\tPi^\prime$. Then we see the following.

\begin{equation}\label{eqn:irrKL}
\begin{array}{l}
\mbox{For $\tal\in\trRptbhmtPi$ and $z\in\bKt$, $z\trK_\tal+\trL_\tal$ is an irreducible element
} \\
\mbox{of $\trUo$ (which is shown by \eqref{eqn:assNlcs} and Lemma~\ref{lemma:trRinfty}),} \\ 
\mbox{and $z^\prime\trK_{\tal^\prime}+\trL_{\tal^\prime}\notin\trUo\cdot(z\trK_\tal+\trL_\tal)$} \\
\mbox{for $\tal^\prime\in\trRptbhmtPi$ and $z^\prime\in\bKt$ with $(\tal^\prime,z^\prime)\ne (\tal,z)$.} 
\end{array}
\end{equation} 
\begin{equation}\label{eqn:assShap-2}
\begin{array}{l}
\mbox{In Subsection~\ref{subsection:NdegShap}, we fix $\tlambda\in\tfkAPip$ and $\tal\in\trRptbhmtPi$,} \\
\mbox{assume 
$0\ne|\trRpmapbctal_\tlambda|
(=\dim\trUp_\tlambda)$, let $\mtlambda:=|\trRpmapbctal_\tlambda|$,} \\
\mbox{let $t\in\bN$, $r:=|\trRpmapbctal_\tlambda(\tal;t)|$, $f:=-\trhobctal(\tal)\tbhm(\tal,\tal)^{-t}
\trK_\tal+\trL_\tal\in\trUo$,} \\
\mbox{let $\dhXhY\in\bK$ and $\mclSbctaltl[\hatX,\hatY]\in\Mat(\mtlambda,\trUo)$ be as in Theorem~\ref{theorem:Shapo},} \\
\mbox{and assume $\dhXhY\ne 0$.}
\end{array}
\end{equation}

\begin{theorem}\label{theorem:nodgthm}
There exist $h\in\trUoSetmsf$ and
$P_z\in\Mat(\mtlambda,\trUo)$
with $\det (P_z)\in\trUoSetmsf$
$(z\in\fkJ_{1,2})$
such that 
\begin{equation}\label{eqn:ndgthEq}
P_1\mclSbctaltl[\hatX,\hatY]P_2=h\cdot (\sum_{x=1}^{\mtlambda-r}E_{x,x}+f\cdot \sum_{y=\mtlambda-r+1}^\mtlambda E_{y,y}).
\end{equation} 
\end{theorem}
{\it{Proof.}}
The claim of Theorem~\ref{theorem:nodgthm} can be obtained 
as an immediate consequence of 
\eqref{eqn:ShapoEq},
\eqref{eqn:irrKL},
Lemma~\ref{lemma:maintool} and Theorem~\ref{theorem:keythma}. \hfill $\Box$

\begin{remark}\label{remark:rcHY10}
By \eqref{eqn:ShapoEq}, we immediately have Equation
\eqref{eqn:eqlmd} for $G=\mclSbctaltl[\hatX,\hatY]$. Theorem~\ref{theorem:nodgthm},
claims that $m^\prime=\mtlambda$, and $c_x=0$ ($x\in\fkJ_{1, \mtlambda-r}$), 
$c_y=1$ ($y\in\fkJ_{\mtlambda-r+1,\mtlambda}$)
for that equation.
\end{remark}

\begin{theorem}\label{theorem:ndgthCor}
Assume $r\geq 1$. Let $P_z$ $(z\in\fkJ_{1,2})$ and $h$ be as in Theorem~{\rm{\ref{theorem:nodgthm}}}.
Then there exists $g\in(\trUoSetmsf)\cap(\trUo\det(P_1P_2 h))$ 
satisfying the condition that
for $\trLambda\in\rmCh(\trUo)$
with $\trLambda(f)=0$ and 
$\trLambda(g)\ne 0$, there exists
$v\in\MtbctaltrLam_{-t\tal}\setminus\{0\}$
such that 
\begin{equation}\label{eqn:ndCEqVan}
\trE_i\cdot v=0\quad\quad (i\in\fkI),
\end{equation} and 
\begin{equation}\label{eqn:ndgthCrEq}
\dim \trUm_{-\tlambda+t\tal}\cdot v=r
\quad\mbox{and}\quad\trUm_{-\lambda+t\tal}\cdot v
=\NtbctaltrLam_{-\tlambda}.
\end{equation}
\end{theorem}
{\it{Proof.}} For each $\tmu\in\tfkAPip$, let $a_\tmu$ be 
the element of $\trUo$ obtained from the one of \eqref{eqn:ShapoEq}
by replacing $\tlambda$ with $\tmu$, and assume $a_\tmu\ne 0$.

Define the $\bK$-algebra automorphism $\trGamma:\trUo\to\trUo$
by 
\begin{equation*}
\trGamma(\trK_\tmu\trL_\tnu):={\frac {\tbhm(\tal,\tnu)^t} {\tbhm(\tmu,\tal)^t}}\trK_\tmu\trL_\tnu
\quad\quad(\tmu,\,\tnu\in\tfkAPi).
\end{equation*}

Let 
\begin{equation*}
g_1:=\prod_{{\tmu\in\tfkAPip,}\atop{t\tal-\tmu\in\tfkAPip\setminus\{0\}}} a_\tmu\in\trUoSetmsf.
\end{equation*}
Let $\trLambda\in\rmCh(\trUo)$ be such that
$\trLambda(f)=0$ and $\trLambda(g_1)\ne 0$.
By \eqref{eqn:EssShpEq} and \eqref{eqn:ShapoEq},
there exists $v\in\MtbctaltrLam_{-t\tal}\setminus\{0\}$
such that $\trE_i\cdot v=0$ ($i\in\fkI$). Note that
\begin{equation}\label{eqn:UvNLm}
\trU\cdot v=\trUm\cdot v\subset 
\NtbctaltrLam.
\end{equation}
Define the $\trU$-module homomorphism 
$p:\Mtbctal(\trLambda\circ\trGamma)\to\MtbctaltrLam$ 
by $p(X\cdot \tv_{\trLambda\circ\trGamma}):=X\cdot v$
($X\in\trU$). 
Let 
\begin{equation}\label{eqn:defgtwo}
g_2:=h\cdot\det(P_1P_2)\in\trUoSetmsf.
\end{equation}
From now on until the end of this proof, we also assume
$\trLambda(g_2)\ne 0$. By \eqref{eqn:EssShpEq} and \eqref{eqn:ndgthEq},
we have
\begin{equation}\label{eqn:dmNLMLb}
\dim\NtbctaltrLam_{-\tlambda}=r.
\end{equation}

Let $\tq_\tal:=\tbhm(\tal,\tal)$ and $\hckpchita:=\kpch(\tq_\tal)$.

{\it{Case-$1$. Assume $\hckpchita=0$.}} 
We have
\begin{equation}\label{eqn:Gmalm}
\forall\tmu\in\tfkAPip,\,\,\trGamma(a_\tmu)\in\trUoSetmsf
\end{equation} since
$\trGamma(-\trhobctal(\tal)\tq_\tal^{-s}\trK_\tal+\trL_\tal)=
\tq_\tal^t\cdot (-\trhobctal(\tal)\tq_\tal^{-s-2t}\trK_\tal+\trL_\tal)$
for all $s\in\bN$ (note that $-s-2t\ne -t$).
Let $g_3:=a_{\tlambda-t\tal}$. Then $\trGamma(g_3)\in\trUoSetmsf$
by \eqref{eqn:Gmalm}.
Assume $(\trLambda\circ\trGamma)(g_3)\ne 0$.
By \eqref{eqn:EssShpEq} and \eqref{eqn:ShapoEq},
we have $\ker p\cap\Mtbctal(\trLambda\circ\trGamma)_{-\tlambda+t\tal}=\{0\}$.
Hence, since $r=|\trRpmapbctal_{\tlambda-t\tal}|=\dim\trUm_{-\tlambda+t\tal}$,
by \eqref{eqn:preLPBW}, \eqref{eqn:UvNLm}
and \eqref{eqn:dmNLMLb}, we have
\begin{equation*}
\NtbctaltrLam_{-\tlambda}=\trUm_{-\tlambda+t\tal}\cdot v.
\end{equation*}
Hence, by \eqref{eqn:dmNLMLb}, we may put $g:=g_1g_2\cdot\trGamma(g_3)$.

{\it{Case-$2$. Assume $\hckpchita\geq 2$.}}
Let $f^\prime:=-\trhobctal(\tal)\tq_\tal^{-\hckpchita+t}\trK_\tal+\trL_\tal$.
Define $g_4\in\trUo\setminus \trUo f^\prime$ 
in the same way as that for $g_2$ (see \eqref{eqn:defgtwo})
with $\tlambda-t\tal$ and $\hckpchita-t$ in place
of $\tlambda$ and $t$ respectively.
Since $\trGamma(f^\prime)=\tq_\tal^tf$, 
we have $\trGamma(g_4)\in\trUoSetmsf$.
Assume $(\trLambda\circ \trGamma)(g_4)\ne 0$.
By  \eqref{eqn:EssShpEq}, \eqref{eqn:preLPBW} and \eqref{eqn:ndgthEq},
we have
\begin{equation}\label{eqn:dmLLGmbpta}
\begin{array}{l}
\dim\Ltbctal(\trLambda\circ\trGamma)_{-\tlambda+t\tal}
=|\trRpmapbctal_{\tlambda-t\tal}|-|\trRpmapbctal_{\tlambda-t\tal}(\tal;\hckpchita-t)| \\
=|\{\,x\in\trRpmapbctal_{\tlambda-t\tal}\,|\,x(\tal)\in\fkJ_{0,\hckpchita-t-1}\,\}|=r.
\end{array}
\end{equation}
Hence $\dim p(\Mtbctal(\trLambda\circ\trGamma)_{-\tlambda+t\tal})\geq r$.
By \eqref{eqn:UvNLm} and \eqref{eqn:dmNLMLb}, we may put 
$g:=g_1g_2\cdot\trGamma(g_4)$.
This completes the proof. 
\hfill $\Box$

\section{Skew centers}\label{section:SkewCenters}
\subsection{Definition of the skew centers}\label{subsection:DefSkewCenters}

Recall \eqref{eqn:abbrU}.
Define the $\bK$-linear map
\begin{equation*}
\trpilqtbtal:\trU\to\rmSpan_\bK(\trUm\trUo)
\end{equation*}
by $\trpilqtbtal(YZX):=\Hopfe(X)YZ$
($X\in\trUp$, $Z\in\trUo$, $Y\in\trUm$).

\begin{lemma}\label{lemma:nondgng}
For $X\in\trU$, if $\trpilqtbtal(XY)=0$ for all $Y\in\trUm$,
then $X=0$.
\end{lemma}
{\it{Proof.}} For $\tlambda\in\tfkAPip$,
let $m_\tlambda:=\dim\trUp_\tlambda$,
and
let  $\{X_{\tlambda,x}\,|\,x\in\fkJ_{1,m_\tlambda}\}$
and $\{Y_{-\tlambda,y}\,|\,y\in\fkJ_{1,m_\tlambda}\}$
be $\bK$-base of $\trUp_\tlambda$
and $\trUm_{-\tlambda}$ respectively.
Let $X$ be as in the statement, and
express it as
\begin{equation*}
X=\sum_{\tlambda,\tmu\in\tfkAPip}\sum_{y=1}^{\mtlambda}\sum_{x=1}^{\mtmu}
Y_{-\tlambda,y}Z^{(\tlambda,\tmu)}_{y,x}X_{\tmu,x},
\end{equation*} where $Z^{(\tlambda,\tmu)}_{y,x}\in\trUo$.
Let $\tmu\in\tfkAPip$ be such that
$Z^{(\tlambda,\tmu)}_{y,x}=0$
for all $\tlambda\in\tfkAPip$ and all $\tmu\in\tfkAPip$
with $\tmu-\tlambda^\prime\in\tfkAPip\setminus\{0\}$.
Then for $y^\prime\in\fkJ_{1,m_\tmu}$, we have
\begin{equation*}
0=\trpilqtbtal(XY_{-\tmu,y^\prime})=\sum_{\tlambda\in\tfkAPip}\sum_{y=1}^{\mtlambda}
\sum_{x=1}^{m_\tmu}
Y_{-\tlambda,y}Z^{(\tlambda,\tmu)}_{y,x}\trShbctal(X_{\tmu,x}Y_{-\tmu,y^\prime}),
\end{equation*} so
$\sum_{x=1}^{m_{\tmu}}Z^{(\tlambda,\tmu)}_{y,x}\trShbctal(X_{\tmu,x}Y_{-\tmu,y^\prime})=0$,
see Theorem~\ref{theorem:DefofGQG}~(1)~{\rm{(GQG4-2)}}.
By \eqref{eqn:trShXYMt}, $Z^{(\tlambda,\tmu)}_{y,x}=0$ for all $\tlambda\in\tfkAPip$ and all
$y\in\fkJ_{1,m_\tlambda}$.
Thus we have $X=0$, as desired. \hfill $\Box$
\newline\par
Let $\newpara:\trtfkAPi\to\bKt$ be a $\bZ$-module homomorphism. 
Note 
\begin{equation*}
\trU_0=\oplus_{\tlambda\in\tfkAPip}
\rmSpan_\bK(\trUm_{-\tlambda}\trUo\trUp_\tlambda),
\end{equation*} where recall that $\trU_0$ means $\trUtbhmtPi_0$,
see \eqref{eqn:abbrU} and \eqref{eqn:abbrUgr}.
Define the $\bK$-subspace $\prtrZnewparabctal$ of $\trU_0$
by
\begin{equation}\label{eqn:defprtrZomega}
\prtrZnewparabctal:=\{\,Z\in\trU_0\,|\,\forall \tlambda\in\tfkAPi,\,\forall X\in\trU_\tlambda,\,
ZX
=\newpara(\tlambda)XZ\,\}.
\end{equation}

Let $\trHCmapbctaln:\prtrZnewparabctal\to\trUo$
be the $\bK$-linear homomorphism defined by
\begin{equation}\label{eqn:deftrHComega}
\trHCmapbctaln:=(\trShbctal)_{|\prtrZnewparabctal}.
\end{equation}
We call $\trHCmapbctaln$ the {\it{$\newpara$-Harish-Chandra map
associated with $(\tbhm,\checktal)$}}.

\begin{lemma}\label{lemma:HCinj}
$\trHCmapbctaln$ is injective.
\end{lemma}
{\it{Proof.}}
Let $Z\in\ker\trHCmapbctaln$, and 
express it as $Z=\sum_{\tlambda\in\tfkAPip}Z_\tlambda$
with $Z_\tlambda\in\rmSpan_\bK(\trUm_{-\tlambda}\trUo\trUp_\tlambda)$.
Note $Z_0=0$.
Then for all $Y=\sum_{\tlambda\in\tfkAPip}Y_{-\tlambda}\in\trUm$
with $Y_{-\tlambda}\in\trUm_{-\tlambda}$, we have 
\begin{equation*}
\trpilqtbtal(ZY)=\sum_{\tlambda\in\tfkAPip}\trpilqtbtal(\newpara(-\tlambda)Y_{-\tlambda}Z)
=\sum_{\tlambda\in\tfkAPip}\newpara(-\tlambda)Y_{-\tlambda}Z_0=0.
\end{equation*}
By Lemma~\ref{lemma:nondgng},
$Z=0$, as desired. \hfill $\Box$

\begin{lemma}\label{lemma:scZg}
Let $X\in\trU_0$ and $Z\in\trUo$.
Assume that
\begin{equation}\label{eqn:ass-scZg}
\forall \tlambda\in\tfkAPip,\,\forall Y_{-\tlambda}\in\trUm_{-\tlambda},\quad
\trpilqtbtal(XY_{-\tlambda})=\newpara(-\tlambda)Y_{-\tlambda}Z.
\end{equation}
Then $X\in\prtrZnewparabctal$ and $\trHCmapbctaln(X)=Z$.
\end{lemma}
{\it{Proof.}} Let $X\in\trU_0$ and $Z\in\trUo$ be as in 
\eqref{eqn:ass-scZg}.
For all $i\in\fkI$,
all $\tlambda\in\tfkAPip$
and all $Y_{-\tlambda}\in\trUm_{-\tlambda}$, we have
\begin{equation*}
\begin{array}{l}
\trpilqtbtal((XE_i-\newpara(\talpii)E_iX)Y_{-\tlambda}) \\
\quad =\trpilqtbtal(XE_iY_{-\tlambda})-\newpara(\talpii-\tlambda)\trpilqtbtal(E_iY_{-\tlambda}Z) \\
\quad =\trpilqtbtal(X[E_i,Y_{-\tlambda}])-\newpara(\talpii-\tlambda)\trpilqtbtal(E_iY_{-\tlambda}Z) \\
\quad =\newpara(\talpii-\tlambda)
\trpilqtbtal([E_i,Y_{-\tlambda}]Z)-\newpara(\talpii-\tlambda)\trpilqtbtal(E_iY_{-\tlambda}Z) \\
\quad =-\newpara(\talpii-\tlambda)\trpilqtbtal(Y_{-\tlambda}E_iZ) \\
\quad =0,
\end{array}
\end{equation*} and
\begin{equation*}
\begin{array}{l}
\trpilqtbtal((XF_i-\newpara(-\talpii)F_iX)Y_{-\tlambda}) \\
\quad=\newpara(-\tlambda-\talpii)F_iY_{-\tlambda}Z
-\newpara(-\tlambda-\talpii)F_iY_{-\tlambda}Z \\
\quad=0.
\end{array}
\end{equation*} 
By Lemma~\ref{lemma:nondgng},
$X\in\prtrZnewparabctal$.
Moreover $\trHCmapbctaln(X)=\trpilqtbtal(X\cdot 1)=\newpara(0)Z=Z$, as desired.
\hfill $\Box$

\subsection{Formulation of Harish-Chandra-type theorem}\label{subsection:SmEq}
Let $i\in\fkI$ and $\ckpchi:=\kpch(\tbhm(\talpii,\talpii))$. Let $\newpara:\tfkAPi\to\bKt$ be a $\bZ$-module homomorphism. We can easily see that
\begin{equation}\label{eqn:prtrZni}
\begin{array}{l}
\mbox{for $X=\sum_{\tlambda\in\tfkAPip}X_\tlambda\in\prtrZnewparabctal$
with $X_\lambda\in\rmSpan_\bK(\trUm_{-\tlambda}\trUo\trUp_\tlambda)$,} \\
\mbox{letting $Y:=\sum_{k=0}^\infty X_{k\talpii}$, we have} \\
\mbox{$Y\in\trUtbhmtPiith\cap\trU_0$
and $Y\trE_i=\newpara(\talpii)\trE_iY$, $Y\trF_i=\newpara(-\talpii)\trF_iY$.}
\end{array}
\end{equation}
Define the $\bK$-algebra isomorphism 
\begin{equation*}
\trjmthtbtalni:\trUo(\tbhm,\tPiangi)\to\trUotbhmtPi
\end{equation*}
by
\begin{equation} \label{eqn:Shiftop}
\trjmthtbtalni(\trK_\tlambda\trL_\tmu):=
\newpara(\talpii)^{\ckpchi-1}\cdot
{\frac {\tbhm(\talpii,\tmu)^{\ckpchi-1}} {\tbhm(\tlambda,\talpii)^{\ckpchi-1}}}
\trK_\tlambda\trL_\tmu
\quad\quad
(\tlambda,\,\tmu\in\tfkAPi),
\end{equation} where $\ckpchi:=\kpch(\tbhm(\talpii,\talpii))$.

\begin{lemma} \label{lemma:sfteq} Assume  $\CardtrRptbtP<\infty$.
Then
\begin{equation} \label{eqn:eqsfteq}
(\trHCmapbctaln\circ\trTtbtali)(X)
=(\trjmth^{\tbhm,\tPiangi}_{\newpara;i}\circ\trHCmapbctalnangi)(X)\quad
\quad(X\in\prtrZnewparabctalnangi),
\end{equation} where note 
$\trTtbtali(\prtrZnewparabctalnangi)
=\prtrZnewparabctal$.
\end{lemma}
{\it{Proof.}} 
For $\tlambda\in\tfkAPipni$,
let 
\begin{equation*}
\acuUtlambda
:=\rmSpan_\bK(\trUmtbtalni_{-\tlambda}\trUotbtalni\trUptbtalni_\tlambda).
\end{equation*} Then $\trUtbtalni_0=\oplus_{\tlambda\in\tfkAPipni}\acuUtlambda$.
Let $X=\sum_{\tlambda\in\tfkAPipni}X_\tlambda\in\prtrZnewparabctalnangi$
with $X_\tlambda\in\acuUtlambda$.
Let $Y:=\sum_{k=0}^\infty X_{-k\talpii}$. Similarly to \eqref{eqn:prtrZni},
we have $Y\in\trU(\tbhm,\tPiangi;i)\cap\trUtbtalni_0$
and $Y\trE_i=\newpara(-\talpii)\trE_iY$, $Y\trF_i=\newpara(\talpii)\trF_iY$,
where note that $\trE_i$, $\trF_i\in\trUtbtalni$.
Then by \cite[Lemma~3.1]{BY15} and the first equation of \eqref{eqn:LuIsoeqn},
we have 
$(\trShbctal\circ\trTtbtali)(Y)=(\trjmth^{\tbhm,\tPiangi}_{\newpara;i}\circ\trShbctalni)(Y)$.
Let $Z:=X-Y$. Clearly $\trShbctalni(Z)=0$.
By \eqref{eqn:extlgdd}, we have
$(\trShbctal\circ\trTtbtali)(Z)=0$.
Hence we have \eqref{eqn:eqsfteq}.
\hfill $\Box$
\newline\par
Assume  $\CardtrRptbtP<\infty$.
Let $\newpara:\tfkAPi\to\bKt$ be a $\bZ$-module homomorphism.
For each $\tbeta\in\trRtbhmtPi(=\trRptbhmtPi\cup(-\trRptbhmtPi))$, 
let $\tfkBbctanptb$
be the linear $\bK$-subspace of $\trUotbhmtPi$ formed by
the elements 
\begin{equation*}
\sum_{(\tlambda,\tmu)\in\tfkAPi^2}a_{(\tlambda,\tmu)}
\trK_\tlambda\trL_\tmu
\end{equation*} with $a_{(\tlambda,\tmu)}\in\bK$
satisfying the following equations $(e1)_\tbeta$-$(e4)_\tbeta$.
In $(e1)_\tbeta$-$(e4)_\tbeta$, let $\tq_\tbeta:=\tbhm(\tbeta,\tbeta)$,
$\hckpchitb:=\kpch(\tq_\tbeta)$ and 
$\newparatlmbtea:=\newpara(\tbeta)\cdot{\frac {\tbhm(\tbeta,\tmu)} {\tbhm(\tlambda,\beta)}}$.
\newline\par
$(e1)_\tbeta$ For $(\tlambda,\tmu)\in\tfkAPi^2$ and $t\in\bZ\setminus\{0\}$,
if $\tq_\tbeta\ne 1$,
$\hckpchitb=0$ 
and $\newparatlmbtea=\tq_\tbeta^t$,
then the equation
$a_{(\tlambda+t\tbeta,\tmu-t\tbeta)}
=\trhobctaltb^t\cdot a_{(\tlambda,\tmu)}$ holds.
\newline\par
$(e2)_\tbeta$ For $(\tlambda,\tmu)\in\tfkAPi^2$, if
$\hckpchitb=0$
and $\newparatlmbtea\ne \tq_\tbeta^t$
for all $t\in\bZ$, then the equation $a_{(\tlambda,\tmu)}=0$ holds.
\newline\par
$(e3)_\tbeta$ 
For $(\tlambda,\tmu)\in\tfkAPi^2$
and $t\in\fkJ_{1,\hckpchitb-1}$,
if $\tq_\tbeta\ne 1$,
$\hckpchitb\geq 2$ and
$\newparatlmbtea=\tq_\tbeta^t$,
then the equation 
\begin{equation*}
\begin{array}{l}
\displaystyle{\sum_{x=-\infty}^\pinfty a_{(\tlambda+(\hckpchitb x+t)\tbeta,\tmu-(\hckpchitb x+t)\tbeta)}
\trhobctaltb^{-(\hckpchitb x+t)}} \\
\quad 
=\displaystyle{\sum_{y=-\infty}^\pinfty a_{(\tlambda+\hckpchitb y\tbeta,\tmu-\hckpchitb y\tbeta)}
\trhobctaltb^{-\hckpchitb y}}
\end{array}
\end{equation*}
holds. 
\newline\par
$(e4)_\tbeta$ For $(\tlambda,\tmu)\in\tfkAPi^2$,
if $\hckpchitb\geq 2$
and $\newparatlmbtea\ne\tq_\tbeta^m$
for all
$m\in\fkJ_{0,\hckpchitb-1}$, then
the $\hckpchitb-1$ equations
\begin{equation*}
\begin{array}{l}
\displaystyle{\sum_{x=-\infty}^\pinfty a_{(\tlambda+(\hckpchitb x+t)\tbeta,\tmu-(\hckpchitb x+t)\tbeta)}
\trhobctaltb^{-(\hckpchitb x+t)}} \\
\quad \displaystyle{=\sum_{y=-\infty}^\pinfty 
a_{(\tlambda+\hckpchitb y\tbeta,\tmu-\hckpchitb y\tbeta)}
\trhobctaltb^{-\hckpchitb y}} \quad(t\in\fkJ_{1,\hckpchitb-1})
\end{array}
\end{equation*} hold.
\newline\newline
Let 
\begin{equation}\label{eqn:deftfkBbctanp}
\tfkBbctanp:=
\bigcap_{\tbeta\in\trRptbhmtPi}\tfkBbctanptb\,(=\bigcap_{\tbeta\in\trRtbhmtPi}\tfkBbctanptb),
\end{equation} where note that $\tfkBbctanp(-\tbeta)=\tfkBbctanptb$.
By \eqref{eqn:tbtbrho} and \eqref{eqn:Shiftop}, we can directly see that 
\begin{equation}\label{eqn:shifttbi}
\trjmthtbtalni(\tfkBbctanpbctaltb)=\tfkBbctanptb\quad (i\in\fkI,\,\tbeta\in\trRtbhmtPi(=\trRtbtal)).
\end{equation}

Now we have the main statement in this subsection.
\begin{proposition}\label{proposition:spsibsc} Assume $\CardtrRptbtP<\infty$.
Then
\begin{equation}\label{eqn:spsibsc-1}
\rmIm\trHCmapbctaln
\subset\tfkBbctanp.
\end{equation}
\end{proposition}
{\it{Proof.}} By \eqref{eqn:prtrZni} and \cite[Lemma~2.6]{BY15}, we have 
$\rmIm\trHCmapbctaln\subset\cap_{i\in\fkI}\tfkBbctanp(\talpii)$.
Tnen \eqref{eqn:spsibsc-1} follows from
\eqref{eqn:trRptbhmtPi}, \eqref{eqn:shifttbi} and Lemma~\ref{lemma:sfteq}.
\hfill $\Box$

\subsection{Same value on a submodule}

\begin{lemma}\label{lemma:smscP}
Let $\tbeta\in\trRptbhmtPi$. Let $\tq_\tbeta:=\tbhm(\tbeta,\tbeta)$ and $\hckpchitb:=\kpch(\tq_\tbeta)$.
Let $\trLambda\in\rmCh(\trUotbhmtPi)$ and $t\in\bN$ 
be such that
\begin{equation}\label{eqn:smscPpreq}
\trLambda(\trK_\tbeta\trL_{-\tbeta})=
{\frac {\tq_\tbeta^t} {\trhobctaltb}}.
\end{equation} Assume that $t\in\fkJ_{1,\hckpchitb-1}$ if $\hckpchitb\geq 2$.
Define $\Lambda^\prime\in\rmCh(\trUotbhmtPi)$ by
\begin{equation}\label{eqn:smscPppreq}
\Lambda^\prime(\trK_\tlambda\trL_\tmu):=
{\frac {\tbhm(\tbeta,\tmu)^t} {\tbhm(\tlambda,\tbeta)^t}}\trLambda(\trK_\tlambda\trL_\tmu)\quad
(\tlambda,\,\tmu\in\tfkAPi). 
\end{equation} Then
\begin{equation}\label{eqn:smscPeq}
\Lambda^\prime(P)=\newpara(\tbeta)^{-t}\Lambda(P)\quad
(P\in\tfkBbctanptb).
\end{equation}
\end{lemma}
{\it{Proof.}} Let $P=\sum_{\tlambda, \tmu\in\tfkAPi}a_{(\tlambda,\tmu)}\trK_\tlambda\trL_\tmu\in\tfkBbctanptb$ with $a_{(\tlambda,\tmu)}\in\bK$.
For $(\tlambda,\tmu)\in\tfkAPi^2$.
let $P_{\tlambda,\tmu}:=\sum_{x=-\infty}^\infty
a_{(\tlambda+x\tbeta,\tmu-x\tbeta)}
\trK_{\tlambda+x\tbeta}\trL_{\tmu-x\tbeta}$,
so $P_{\tlambda,\tmu}\in\tfkBbctanptb$.
To show \eqref{eqn:smscPeq}, 
it suffices to show that for all $(\tlambda,\tmu)\in\tfkAPi^2$,
\begin{equation*}
(*)_{\tlambda,\tmu}\quad\quad\quad\Lambda^\prime(P_{\tlambda,\tmu})=\newpara(\tbeta)^{-t}\Lambda(P_{\tlambda,\tmu})
\end{equation*} hold.

By \eqref{eqn:smscPpreq}
and \eqref{eqn:smscPppreq}, we have
\begin{equation}\label{eqn:smscPpf-1}
\Lambda^\prime(\trK_\tbeta\trL_{-\tbeta})=
{\frac {\tq_\tbeta^{-t}}
{\trhobctaltb}}.
\end{equation}

Fix $(\tlambda,\tmu)\in\tfkAPi^2$ until the end of this proof.
Let $\newparatlmbtea:=\newpara(\tbeta)\cdot{\frac {\tbhm(\tbeta,\tmu)} {\tbhm(\tlambda,\beta)}}$.

{\it{Case-{\rm{0}}}}. If $\tq_\tbeta=\newparatlmbtea=1$, $(*)_{\tlambda,\tmu}$ is clear.

{\it{Case-{\rm{1}}}}. Assume that
$\hckpchitb=0$
and 
\begin{equation}\label{eqn:smscPas-1}
\newparatlmbtea=\tq_\tbeta^s
\quad\mbox{for some $s\in\bZ$.} 
\end{equation} By {\it{Case-{\rm{0}}}}, we may assume $\tq_\tbeta\ne 1$.
Let
\begin{equation*}
P^\prime_{\tlambda,\tmu}:=
\left\{\begin{array}{ll}
a_{(\tlambda,\tmu)}\trK_\tlambda\trL_\tmu & \mbox{if $s=0$}, \\
a_{(\tlambda,\tmu)}\trK_\tlambda\trL_\tmu
+a_{(\tlambda+s\tbeta,\tmu-s\tbeta)}
\trK_{\tlambda+s\tbeta}\trL_{\tmu-s\tbeta} & \mbox{if $s\ne 0$}. 
\end{array}\right.
\end{equation*}
Since  $\newpara_{\tlambda+x\tbeta,\tmu-x\tbeta;\tbeta}
=\tq_\tbeta^{s-2x}$,
to show $(*)_{\tlambda,\tmu}$,
it suffices to show that \eqref{eqn:smscPeq}
for $P^\prime_{\tlambda,\tmu}$ holds.
If $s=0$, it is clear. Assume $s\ne 0$. Then
\begin{equation*}
\begin{array}{l}
\Lambda^\prime(P^\prime_{\tlambda,\tmu}) \\
\quad = \Lambda^\prime(\trK_\tlambda\trL_\tmu)\cdot
(a_{(\tlambda,\tmu)}+a_{(\tlambda+s\tbeta,\tmu-s\tbeta)}
\Lambda^\prime(\trK_\tbeta\trL_{-\tbeta})^s) \\
\quad = \Lambda^\prime(\trK_\tlambda\trL_\tmu)\cdot
a_{(\tlambda,\tmu)}\cdot
(1+\trhobctaltb^s
\Lambda^\prime(\trK_\tbeta\trL_{-\tbeta})^s)
\quad\mbox{(by $(e1)_\tbeta$)}  \\
\quad = {\frac {\tbhm(\tbeta,\tmu)^t} {\tbhm(\tlambda,\tbeta)^t}}\Lambda(\trK_\tlambda\trL_\tmu)\cdot a_{(\tlambda,\tmu)}\cdot
(1+\trhobctaltb^s\cdot\tq_\tbeta^{-st}\trhobctaltb^{-s}) \\
\quad\quad\quad\mbox{(by 
\eqref{eqn:smscPppreq} and \eqref{eqn:smscPpf-1})} \\
\quad = \newpara(\tbeta)^{-t}\tq_\tbeta^{st}\cdot\Lambda(\trK_\tlambda\trL_\tmu)\cdot a_{(\tlambda,\tmu)}\cdot
(1+\tq_\tbeta^{-st}) \quad\mbox{(by \eqref{eqn:smscPas-1})} \\
\quad = \newpara(\tbeta)^{-t}\Lambda(\trK_\tlambda\trL_\tmu)\cdot a_{(\tlambda,\tmu)}\cdot
(\tq_\tbeta^{st}+1) \\
\quad = \newpara(\tbeta)^{-t}\Lambda(\trK_\tlambda\trL_\tmu)\cdot a_{(\tlambda,\tmu)}\cdot
(\trhobctaltb^s\Lambda(\trK_\tbeta\trL_{-\tbeta})^s+1) 
\quad\mbox{(by \eqref{eqn:smscPpreq})} \\
\quad = \newpara(\tbeta)^{-t}\Lambda(\trK_\tlambda\trL_\tmu)\cdot
(a_{(\tlambda+s\tbeta,\tmu-s\tbeta)}\Lambda(\trK_\tbeta\trL_{-\tbeta})^s+a_{(\tlambda,\tmu)})\quad\mbox{(by $(e1)_\tbeta$)} \\
\quad = \newpara(\tbeta)^{-t}\Lambda(P^\prime_{\tlambda,\tmu}).
\end{array}
\end{equation*} 
Hence $(*)_{\tlambda,\tmu}$ holds,
as desired.

{\it{Case-{\rm{2}}}}. Assume that
$\hckpchitb=0$ 
and $\newparatlmbtea\ne \tq_\tbeta^s$
for all $s\in\bZ$.
Then $(e2)_\tbeta$ implies $P_{\tlambda,\tmu}=0$,
since  $\newpara_{\tlambda+x\beta,\tmu-x\tbeta;\tbeta}
=\tq_\tbeta^{-2x}\cdot\newparatlmbtea$.
Hence $(*)_{\tlambda,\tmu}$ holds,
as desired.

{\it{Case-{\rm{3}}}}. 
Assume that
$\hckpchitb\geq 2$ and
\begin{equation}\label{eqn:smscPas-3}
\newparatlmbtea=\tq_\tbeta^s
\quad\mbox{for some $s\in\fkJ_{0,\hckpchitb-1}$}.
\end{equation} By {\it{Case-{\rm{0}}}}, we may assume $\tq_\tbeta\ne 1$.
Let
\begin{equation*}
P^{\prime\prime}_{\tlambda,\tmu}:=
\left\{\begin{array}{ll}
\sum_{x=-\infty}^\infty a_{(\tlambda+\hckpchitb x\tbeta,\tmu-\hckpchitb x\tbeta)}\trK_{\tlambda+\hckpchitb x\tbeta}\trL_{\tmu-\hckpchitb x\tbeta} & \mbox{if $s=0$}, \\
\sum_{x=-\infty}^\infty (a_{(\tlambda+\hckpchitb x\tbeta,\tmu-\hckpchitb x\tbeta)}\trK_{\tlambda+\hckpchitb x\tbeta}\trL_{\tmu-\hckpchitb x\tbeta}
\\ \quad
+a_{(\tlambda+(\hckpchitb x+s)\tbeta,\tmu-(\hckpchitb x+s)\tbeta)}\trK_{\tlambda+(\hckpchitb x+s)\tbeta}\trL_{\tmu-(\hckpchitb x+s)\tbeta}) & \mbox{if $s\ne 0$}. 
\end{array}\right.
\end{equation*}
To show $(*)_{\tlambda,\tmu}$,
it suffices to show that
\eqref{eqn:smscPeq} for $P^{\prime\prime}_{\tlambda,\tmu}$ holds,
and it is clear if $s=0$.
Assume $s\in\fkJ_{1,\hckpchitb-1}$.
Then
\begin{equation*}
\begin{array}{l}
\Lambda^\prime(P^{\prime\prime}_{\tlambda,\tmu}) \\
\quad = \Lambda^\prime(\trK_\tlambda\trL_\tmu)\cdot
\sum_{x=-\infty}^\infty 
(a_{(\tlambda+\hckpchitb x\tbeta,\tmu-\hckpchitb x\tbeta)}\Lambda^\prime(\trK_\tbeta\trL_{-\tbeta})^{\hckpchitb x} \\
\quad\quad +a_{(\tlambda+(\hckpchitb x+s)\tbeta,\tmu-(\hckpchitb x+s)\tbeta)}
\Lambda^\prime(\trK_\tbeta\trL_{-\tbeta})^{\hckpchitb x+s}) \\
\quad = {\frac {\tbhm(\tbeta,\tmu)^t} {\tbhm(\tlambda,\tbeta)^t}}\Lambda(\trK_\tlambda\trL_\tmu)\cdot
\sum_{x=-\infty}^\infty 
(a_{(\tlambda+\hckpchitb x\tbeta,\tmu-\hckpchitb x\tbeta)}\trhobctaltb^{-\hckpchitb x} \\
\quad\quad +a_{(\tlambda+(\hckpchitb x+s)\tbeta,\tmu-(\hckpchitb x+s)\tbeta)}
\tq_\tbeta^{-st}\trhobctaltb^{-(\hckpchitb x+s)}) \\
\quad\quad\quad\mbox{(by \eqref{eqn:smscPppreq} and \eqref{eqn:smscPpf-1})} \\
\quad = \newpara(\tbeta)^{-t}\tq_\tbeta^{st}\cdot\Lambda(\trK_\tlambda\trL_\tmu)\cdot(1+\tq_\tbeta^{-st})\cdot
\sum_{x=-\infty}^\infty a_{(\tlambda+\hckpchitb x\tbeta,\tmu-\hckpchitb x\tbeta)}\trhobctaltb^{-\hckpchitb x} \\
\quad\quad\quad\mbox{(by \eqref{eqn:smscPas-3} and $(e3)_\tbeta$)} \\
\quad = \newpara(\tbeta)^{-t}\Lambda(\trK_\tlambda\trL_\tmu)(1+\tq_\tbeta^{st})
\sum_{x=-\infty}^\infty a_{(\tlambda+\hckpchitb x\tbeta,\tmu-\hckpchitb x\tbeta)}\trhobctaltb^{-\hckpchitb x} \\
\quad = \newpara(\tbeta)^{-t}\Lambda(\trK_\tlambda\trL_\tmu)
\sum_{x=-\infty}^\infty (a_{(\tlambda+\hckpchitb x\tbeta,\tmu-\hckpchitb x\tbeta)}\trhobctaltb^{-\hckpchitb x} \\
\quad\quad +\tq_\tbeta^{st}\cdot a_{(\tlambda+(\hckpchitb x+s)\tbeta,\tmu-(\hckpchitb x+s)\tbeta)}\trhobctaltb^{-(\hckpchitb x+s)})\quad\mbox{(by $(e3)_\tbeta$)} \\
\quad = \newpara(\tbeta)^{-t}\Lambda(\trK_\tlambda\trL_\tmu)
\sum_{x=-\infty}^\infty (a_{(\tlambda+\hckpchitb x\tbeta,\tmu-\hckpchitb x\tbeta)}\Lambda(\trK_\tbeta\trL_{-\tbeta})^{\hckpchitb x}+ \\
\quad\quad +a_{(\tlambda+(\hckpchitb x+s)\tbeta,\tmu-(\hckpchitb x+s)\tbeta)}\Lambda(\trK_\tbeta\trL_{-\tbeta})^{(\hckpchitb x+s)}) 
\quad\mbox{(by \eqref{eqn:smscPpreq})} \\
\quad = \newpara(\tbeta)^{-t}\Lambda(P^{\prime\prime}_{\tlambda,\tmu}).
\end{array} 
\end{equation*}
Hence \eqref{eqn:smscPeq}
for $P^{\prime\prime}_{\tlambda,\tmu}$ holds,
as desired.

{\it{Case-{\rm{4}}}}. 
Assume that
$\hckpchitb\geq 2$ and
\begin{equation*}
\newparatlmbtea\ne \tq_\tbeta^s
\quad\mbox{for all $s\in\fkJ_{0,\hckpchitb-1}$}.
\end{equation*} Then
\begin{equation*}
\begin{array}{l}
\Lambda(P_{\tlambda,\tmu}) \\
\quad = \Lambda(\trK_\tlambda\trL_\tmu)\cdot
\sum_{x=-\infty}^\infty 
a_{(\tlambda+x\tbeta,\tmu-x\tbeta)}\Lambda(\trK_\tbeta\trL_{-\tbeta})^x \\
\quad = \Lambda(\trK_\tlambda\trL_\tmu)\cdot
\sum_{x=-\infty}^\infty 
a_{(\tlambda+x\tbeta,\tmu-x\tbeta)}\tq_\tbeta^{xt}\trhobctaltb^{-x}\quad\mbox{(by \eqref{eqn:smscPpreq})} \\
\quad = \Lambda(\trK_\tlambda\trL_\tmu)\cdot
\sum_{z=0}^{\hckpchitb-1}\tq_\tbeta^{zt}\sum_{y=-\infty}^\infty 
a_{(\tlambda+(\hckpchitb y+z)\tbeta,\tmu-(\hckpchitb y+z)\tbeta)}\trhobctaltb^{-(\hckpchitb y+z)} \\
\quad = \Lambda(\trK_\tlambda\trL_\tmu)\cdot
(\hckpchitb)_{\tq_\tbeta^t}\cdot
\sum_{y=-\infty}^\infty 
a_{(\tlambda+\hckpchitb y\tbeta,\tmu-\hckpchitb y\tbeta)}
\trhobctaltb^{-\hckpchitb y}\quad\mbox{(by $(e4)_\tbeta$)} \\
\quad = 0 \quad\mbox{(since 
$(\hckpchitb)_{\tq_\tbeta^t}=0$)}.
\end{array}
\end{equation*} Similarly, by \eqref{eqn:smscPpf-1},
we also have $\Lambda^\prime(P_{\tlambda,\tmu})=0$.
Hence  $(*)_{\tlambda,\tmu}$ holds,
as desired. This completes the proof.
\hfill $\Box$

\section{Kac argument} \label{section:KKarg}
\subsection{Matrix equations}\label{subsection:ME}
Note \eqref{eqn:abbrU}. In this subsection, keep the notation and assumption of 
Theorem~\ref{theorem:Shapo}; however let
\begin{equation}\label{eqn:defremclS}
\remclS:=\mclSbctaltl[\hatX,\hatY]\in\Mat(m_\tlambda,\trUo).
\end{equation}
Assume that \eqref{eqn:dettrShXYMt} is the case, so $\det\remclS\ne 0$.
Assume $\tlambda\ne 0$.
Let $\tfkAPiptlambdast:=\{\tmu\in\tfkAPip|\tlambda-\tmu\in\tfkAPip,\tmu\ne\tlambda\}$.
For $\tmu\in\tfkAPiptlambdast$ with $\trRpmapbctal_\tmu\ne\emptyset$,
let $\mtmu:=|\trRpmapbctal_\tmu|$, 
and let $\{\,X_{\tmu,x}\,|\,x\in\fkJ_{1,\mtmu}\,\}$
and $\{\,Y_{-\tmu,y}\,|\,y\in\fkJ_{1,\mtmu}\,\}$
be base of $\trUp_\tmu$ and $\trUm_{-\tmu}$
respectively (recall \eqref{eqn:DFtrRPmp} and 
\eqref{eqn:preLPBW}).
Assume $X_{0,0}=Y_{0,0}=1$.

Let $Z_{\tmu,y,x}\in\trUo$ for $\tmu\in\tfkAPiptlambdast$
and $x$, $y\in\fkJ_{1,\mtmu}$.
Set
\begin{equation*}
V_{<\tlambda}:=\sum_{\tmu\in\tfkAPiptlambdast}\sum_{x=1}^{\mtmu}\sum_{y=1}^{\mtmu}
Y_{-\tmu,x}Z_{\tmu,x,y}X_{\tmu,y}.
\end{equation*}

\begin{equation}\label{eqn:assVlqb}
\begin{array}{l}
\mbox{Let $P\in\tfkBbctanp$, and assume that}\\ 
\mbox{$\trpilqtbtal(V_{<\tlambda}Y^\prime)=\newpara(-\tmu)Y^\prime P$
for all $\tmu\in\tfkAPiptlambdast$ and all
$Y^\prime\in\trUm_{-\tmu}$.}
\end{array}
\end{equation}
Then $Z_{0,0,0}=P$.

Recall $\hatX(x)$ and $\hatY(y)$ from Theorem~\ref{theorem:Shapo}.
Define $G_{x,y}\in\trUo$ with $x$, $y\in\fkJ_{1,m}$
by
\begin{equation}\label{eqn:lwreq}
\trpilqtbtal(V_{<\tlambda}\hatY(y))=\newpara(-\tlambda)\hatY(y)P
-\sum_{x=1}^\mtlambda \hatY(x)G_{x,y}
\quad\quad(y\in\fkJ_{1,\mtlambda}).
\end{equation}
Let $Z_{x,y}\in\trUo$ for $x$, $y\in\fkJ_{1,\mtlambda}$.
Set
\begin{equation*}
V_{\leq\tlambda}:=V_{<\tlambda}+\sum_{x=1}^\mtlambda\sum_{y=1}^\mtlambda
\hatY(y)Z_{y,x}\hatX(x).
\end{equation*} 
For $y\in\fkJ_{1,\mtlambda}$, by \eqref{eqn:lwreq}, we have
\begin{equation}\label{eqn:preMatEq}
\begin{array}{l}
\trpilqtbtal(V_{\leq\tlambda}\hatY(y)) \\
\quad
=\trpilqtbtal(V_{<\tlambda}\hatY(y))+\sum_{x^\prime=1}^\mtlambda\sum_{y^\prime=1}^\mtlambda
\trpilqtbtal(\hatY(y^\prime)Z_{y^\prime,x^\prime}\hatX(x^\prime)\hatY(y))
\\
\quad = \newpara(-\tlambda)\hatY(y)P
-\sum_{x=1}^\mtlambda\hatY(x)G_{x,y} \\
\quad\quad +\sum_{x^\prime=1}^\mtlambda\sum_{y^\prime=1}^\mtlambda
\hatY(y^\prime)Z_{y^\prime,x^\prime}
\trShbctal(\hatX(x^\prime)\hatY(y)).
\end{array}
\end{equation}
Define 
$\mclG\in\Mat(\mtlambda,\trUo)$
by
$\mclG:=[G_{x,y}]_{1\leq x,y\leq\mtlambda}$.
Define 
$\mclZ\in\Mat(\mtlambda,\trUo)$
by 
$\mclZ:=[Z_{x,y}]_{1\leq x,y\leq\mtlambda}$.
By \eqref{eqn:preMatEq}, we have
\begin{equation}\label{eqn:CeqZG}
\mbox{$\trpilqtbtal(V_{\leq\tlambda}\hatY(y))=\newpara(-\tlambda)\hatY(y)P$
for all $y\in\fkJ_{1,\mtlambda}$
if and only if $\mclG=\mclZ \remclS$},
\end{equation} where recall $\remclS$ from \eqref{eqn:defremclS}.

For $\mclZ^\prime=[Z^\prime_{x,y}]_{1\leq x,y\in\leq\mtlambda}$
$\in
\Mat(m,\trUo)$ 
and $\trLambda\in
\rmCh(\trUo)$,
let 
\begin{equation*}
\trLambda(\mclZ^\prime):=[\trLambda(Z^\prime_{x,y})]_{1\leq x,y\in\leq\mtlambda}
\in\Mat(\mtlambda,\bK).
\end{equation*}

In the proof of Lemma~\ref{lemma:EXTlem} below,
we use a well-known argument originally give in \cite{Kac84}
(see also \cite[The proof of Theorem~13.1.1]{Musson12}).

\begin{lemma}\label{lemma:EXTlem}
There exist $Z_{x,y}\in\trUo$ 
$(x,\,y\in\fkJ_{1,m})$
such that \eqref{eqn:CeqZG} is the case.
\end{lemma}
{\it{Proof.}}
Let $\tal\in\trRptbhmtPi$ and $t\in\bN$
be such that $\trRpmapbctal_\tlambda(\tal;t)\ne\emptyset$.
Let $f\in\trUo$ and $g\in\trUoSetmsf$ be the ones of
Theorems~\ref{theorem:ndgthCor}
for $\tal$ and $t$.
Let $\trLambda\in\rmCh(\trUo)$ be such that
$\trLambda(f)=0$ and $\trLambda(g)\ne 0$.
Let $v\in\MtbctaltrLam_{-t\tal}\setminus\{0\}$ be
as in Theorem~\ref{theorem:ndgthCor}.
For all $Y^\prime\in\trUm_{-\tlambda+t\tal}$,
we have
\begin{equation}\label{eqn:VlqbEQ}
\begin{array}{l}
V_{<\tlambda}Y^\prime\cdot v \\
\quad = \newpara(-\tlambda+t\tal)Y^\prime P\cdot v
\quad(\mbox{by \eqref{eqn:assVlqb}})\\
\quad =\newpara(-\tlambda+t\tal)\newpara(\tal)^{-t}\trLambda(P)Y^\prime\cdot v
\quad(\mbox{by Lemma~\ref{lemma:smscP}}) \\
\quad =\newpara(-\tlambda)\trLambda(P)Y^\prime\cdot v.
\end{array}
\end{equation} By \eqref{eqn:ndgthCrEq} and \eqref{eqn:VlqbEQ},
we have 
\begin{equation}\label{eqn:VlqbEQd}
V_{<\tlambda}\cdot v^\prime=\newpara(-\tlambda)\trLambda(P)v^\prime
\quad(v^\prime\in\mclN(\trLambda)_{-\tlambda}).
\end{equation} Let $z_\mtlambda={^t}[z_1,\ldots,z_\mtlambda]\in\bK^\mtlambda$
and $Y^z:=\sum_{y=1}^\mtlambda z_y\hatY(y)\in \trUm_{-\tlambda}$.
By \eqref{eqn:lwreq}, we have
\begin{equation}\label{eqn:lwreqd}
V_{<\tlambda}Y^z\cdot\tv_\trLambda 
= \newpara(-\tlambda)\trLambda(P)Y^z\cdot\tv_\trLambda
-\sum_{x,y=1}^\mtlambda z_y\trLambda(G_{x,y})
\hatY(x)
\cdot\tv_\trLambda
\end{equation} By \eqref{eqn:EssShpKr} and \eqref{eqn:ndgthCrEq}, 
\begin{equation}\label{eqn:Yztv}
\ker\trLambda(\remclS)=\{z\in\bK^m|Y^z\tv_\trLambda\in\mclN(\trLambda)_{-\tbeta}\}\,\,
\mbox{and}\,\,\dim\ker\trLambda(\remclS)=r,
\end{equation} where recall $\remclS$ from \eqref{eqn:defremclS}.
By \eqref{eqn:VlqbEQd}, \eqref{eqn:lwreqd} and \eqref{eqn:Yztv}, we have
\begin{equation}\label{eqn:SLCL}
\ker \trLambda(\remclS)\subset \ker \trLambda(\mclG)
\end{equation} (for any for $\trLambda\in\rmCh(\trUo)$
with $\trLambda(f)=0$ and $\trLambda(g)\ne 0$).

Recall that $\trUo$ is a unique factorization domain.
As in \eqref{eqn:dftreta},
we can identify $\rmCh(\trUo)$
with the affine space $(\bKt)^{2\CardfkIN}$,
where recall $\CardfkIN=\CardfkI$.
Let $r:=|\trRpmapbctal_\tlambda(\tal;t)|$,
as in Theorem~\ref{theorem:nodgthm}.
By \eqref{eqn:ShapoEq}, $f^r$ divides $\det\remclS$
and $\rmgcd \{f^r,{\frac {\det\remclS} {f^r}}\}=1$.
Then,
by \eqref{eqn:ShapoEq}, 
\eqref{eqn:irrKL}, 
\eqref{eqn:Yztv}, \eqref{eqn:SLCL}
and Lemma~\ref{lemma:maintool}, we see that
\begin{equation}\label{eqn:CSKLZ}
\mclG\remclS^{-1}=
{\frac {f^r} {\det\remclS}}
\mclZ^\prime \,\,\mbox{for 
some $\mclZ^\prime\in\Mat(\mtlambda,\trUo)$}.
\end{equation}

Recall
from \eqref{eqn:ShapoEq} that
\begin{equation*}
\rmgcd\{\,{\frac {\det\remclS} {(-\trho(\tbeta)
\tq_\tbeta^{-t}\trK_\tbeta+\trL_\tbeta)^{r(\tbeta;t)}}}\,|\,\tbeta\in\trRp,\,
t\in\bN,\,\trRpmapbctal_\tlambda(\tbeta;t)\ne\emptyset\,\}=1,
\end{equation*} where $r(\tbeta;t):=|\trRpmapbctal_\tlambda(\tbeta;t)|$.
Hence, by  \eqref{eqn:CSKLZ}, we can see $\mclG\remclS^{-1} 
\in\Mat(\mtlambda,\trUo)$. This completes the proof.
\hfill $\Box$

\subsection{Estimation of degree of $\mclZ$}\label{subsection:est}
Keep the notation and assumption of Subsection~\ref{subsection:ME}.
For $t\in\bZgeqo$, define the
$\bK$-linear subspace $\trU^{0,(t)}$ of $\trUo$ as follows.
Let $\trU^{0,(0)}:=\bK 1_\trU$.
For $t\in\bN$, let $\trU^{0,(t)}:=\sum_{i\in\fkI}
(\trK_\talpii\trU^{0,(t-1)}+\trL_\talpii\trU^{0,(t-1)})$,
and $\trU^{0,(-t)}:=\{0\}$. 
Note that $\trU^{0,(t)}\trU^{0,(t^\prime)}\subset\trU^{0,(t+t^\prime)}$
for $t$, $t^\prime\in\bZgeqo$.
Let $\trU^{0,(\infty)}:=\oplus_{t=0}^\infty\trU^{0,(t)}$.
Then $\trU^{0,(\infty)}=\oplus_{\tmu,\tnu\in\tfkAPip}\bK\trK_\tmu\trL_{\tnu}$,
that is to say,
\begin{equation}\label{eqn:idnU}
\begin{array}{l}
\mbox{$\trU^{0,(\infty)}$
can be identified with the polynomial $\bK$-algebra} \\
\mbox{in the $2\CardfkIN$-variables
$\trK_\talpii$, $\trL_\talpii$
($i\in\fkI$).}
\end{array}
\end{equation}
Let $l\in\bZgeqo$ be such that
$\trK_\tlambda\in\trU^{0,(l)}$. 
By \eqref{eqn:trShXY} and \eqref{eqn:dettrShXYMt}, we have 
\begin{equation}\label{eqn:predegdetS}
\remclS\in\Mat(\mtlambda,\trU^{0,(l)})\setminus\{0\},
\end{equation}
and
\begin{equation}\label{eqn:degdetS}
\det\remclS\in\trU^{0,(\mtlambda\cdot l)}\setminus\{0\}.
\end{equation}
\begin{lemma}\label{lemma:estlm}
Let $P$ and $\mclZ$ be as in Lemma~{\rm{\ref{lemma:EXTlem}}}.
Let $\tmu$, $\tnu\in\tfkAPip$ be such that 
$\trK_\tmu\trL_{\tnu}P\in\oplus_{t=0}^k\trU^{0,(t)}$ for some $k\in\bZgeqo$.
Assume that \eqref{eqn:CeqZG} is the case.
Assume $l>k$. Then $\mclZ=0$.
\end{lemma}
{\it{Proof.}} By the definition of $\tfkBbctanp$, we may 
assume $\trK_\tmu\trL_{\tnu}P\in\trU^{0,(k)}$.
We show that
\begin{equation}\label{eqn:estlma}
\trK_\tmu\trL_{\tnu} Z_{x,y}\in\trU^{0,(k-l)}
\quad(x,\,y\in\fkJ_{1,\mtlambda}).
\end{equation} We use an induction on $l$. By \eqref{eqn:lwreq}
and \eqref{eqn:estlma} (for lower $l$'s),
we may assume that 
\begin{equation}\label{lemma:estlmD}
\trK_\tmu\trL_{\tnu}G_{x,y}\in\trU^{0,(k)}
\quad(x,\,y\in\fkJ_{1,\mtlambda}).
\end{equation}

Let $\retmclS=[\retmclS_{x,y}]_{1\leq x,y\leq \mtlambda}$
be the cofactor matrix of $\remclS$,
so $\remclS^{-1}={\frac 1 {\det\remclS}}\retmclS$.
Recall $\mclG=[G_{x,y}]_{1\leq x,y\leq m_\tlambda}$ from \eqref{eqn:CeqZG}.
Then $\mclG\retmclS=\det\remclS\cdot\mclZ$.
By \eqref{lemma:estlmD}, $\trK_\tmu\trL_{\tnu}\mclG\in\Mat(\mtlambda,\trU^{0,(k)})$.
By \eqref{eqn:predegdetS}, $\retmclS\in\Mat(\mtlambda,\trU^{0,((\mtlambda-1)l)})$.
Hence $\trK_\tmu\trL_{\tnu}\mclG\retmclS\in\Mat(\mtlambda,\trU^{0,(k+(\mtlambda-1)l)})$.
Note that $\mclZ\in\Mat(\mtlambda,\trUo)$
and $\trUo$ can be identified with the Laurent polynomial $\bK$-algebra
in the $2\CardfkIN$-variables
$\trK_\talpii$, $\trL_\talpii$
($i\in\fkI$).
By \eqref{eqn:ShapoEq}, \eqref{eqn:irrKL} and \eqref{eqn:idnU}, 
we have $\trK_\tmu\trL_{\tnu}\mclZ\in\Mat(\mtlambda,\trU^{0,(k-l)})$,
which means that \eqref{eqn:estlma} holds.
This completes the proof. \hfill $\Box$

\subsection{Harish-Chandra-type theorem}\label{subsection:sbmainTh}
\begin{lemma}\label{lemma:premainTh}
For any $P\in\tfkBbctanp$,
there exists a unique $V\in\prtrZnewparabctal$ such
that $\trHCmapbctaln(V)=P$.
\end{lemma}
{\it{Proof.}} By \eqref{eqn:assVlqb}, \eqref{eqn:CeqZG}, 
and Lemmas~\ref{lemma:EXTlem},
and \ref{lemma:estlm},
for any $P\in\tfkBbctanp$,
there exists $V\in\trU_0$
such that the same equations as those of \eqref{eqn:ass-scZg}
with $P$ and $V$ in place of 
$X$ and $Z$ respectively hold.
Then this lemma follows from 
Lemmas~\ref{lemma:HCinj} and 
\ref{lemma:scZg}.
\hfill $\Box$
\newline\par
By Proposition~\ref{proposition:spsibsc} and Lemma~\ref{lemma:premainTh},
we have our main theorem below.

\begin{theorem}\label{theorem:mainTh}
Assume that $\CardtrRptbtP<\infty$.
Assume that $\tbhm(\tal,\tal)\ne 1$ for all $\tal\in\trRptbhmtPi$.
Then 
\begin{equation*}
\rmIm\trHCmapbctaln=\tfkBbctanp.
\end{equation*}
{\rm{(}}Recall 
that $\trHCmapbctaln$ is injective, see Lemma~{\rm{\ref{lemma:HCinj}}}.{\rm{)}}
\end{theorem} 

\section{Symmetric case}\label{section:symcase}
\subsection{Subspaces of skew graded centers}\label{subsection:CmSgc}
In Subsection~\ref{subsection:CmSgc},
assume $\tbhm$ to be any bi-homomorphism such that
$\tbhm(\tlambda,\tmu)=\tbhm(\tmu,\tlambda)$ ($\tlambda$, $\tmu\in\tfkAPi$).

Define the $\bK$-subalgebra $\trUodag$ of $\trUo=\trUotbhmtPi$ by 
$\trUodag:=\oplus_{\tlambda\in\tfkAPi}\bK\trK_\tlambda\trL_{-\tlambda}$.
Define the $\bK$-subspace $\trUdag_0$ of $\trU_0=\trUtbhmtPi_0$ by
\begin{equation*}
\trUdag_0:=\bigoplus_{\tlambda\in\tfkAPip}\rmSpan_\bK(
\trUm_{-\tlambda}\trK_{-\tlambda}\trUodag\trUp_\tlambda).
\end{equation*}
Then for a $\bZ$-module homomorphism 
$\newpara:\trtfkAPi\to\bKt$, we have 
\begin{equation*}
\prtrZnewparabctal=
\bigoplus_{\tlambda\in\tfkAPi}(\prtrZnewparabctal
\cap\trK_\tlambda\trUdag_0).
\end{equation*}

Let $\chkI$ be the two-sided ideal of
the $\bK$-algebra $\trU=\trUtbhmtPi$ generated by all the elements
$\trK_\tlambda\trL_\tlambda-1$ ($\tlambda\in\trtfkAPi$).
Let $\chktrU:=\trU/\chkI$ (the quotient $\bK$-algebra),
and let $\chkpi:\trU\to\chktrU$ be the canonical map.
Let $\chktrK_\tlambda:=\chkpi(\trK_\tlambda)$,
($\tal\in\trtfkAPi$),
and $\chktrE_i:=\chkpi(\trE_i)$, $\chktrF_i:=\chkpi(\trF_i)$
($i\in\fkI$).
Let $\chktrUo:=\chkpi(\trUo)$.
Let $\chktrU_\tlambda:=\chkpi(\trU_\tlambda)$
($\tlambda\in\trtfkAPi$).
We can easily see that $\ker(\chkpi_{|\trUp})=\ker(\chkpi_{|\trUm})=\{0\}$,
$\chktrUo=\oplus_{\tlambda\in\trtfkAPi}\bK\chktrK_\tlambda$, 
and that the $\bK$-linear homomorphism 
$\chkpi(\trUm)\otimes\chktrUo\otimes\chkpi(\trUp)\to\chktrU$ defined by
$\chktrXm\otimes\chktrY\otimes\chktrXp\mapsto\chktrXm\chktrY\chktrXp$ is bijective.

Let $\newpara:\trtfkAPi\to\bKt$ be a $\bZ$-module homomorphism.
Let $\chkprtrZ_\newpara(\tbhm,\checktal)$ be the $\bK$-subspace of $\chktrU^y$
defined 
in the same way as that for $\prtrZnewparabctal$
with $\chktrU_\tmu$
in place of $\trU_\tmu$ ($\tmu\in\tfkAPi$),
where note that 
$\chktrU=\oplus_{\tlambda \in\trtfkAPi}\chktrU_\tlambda$.
Let $\tfkAPi^\prime:=\{\,\sum_{i\in\fkI}x_i\talpii\in\tfkAPi\,|\,x_i\in\fkJ_{0,1}\,(i\in\fkI)\,\}$.
We can see
\begin{equation*}
\chkprtrZ_\newpara(\tbhm,\checktal)=\bigoplus_{\tlambda\in\tfkAPi^\prime}(\chkprtrZ_\newpara(\tbhm,\checktal)
\cap\chktrK_\tlambda\chkpi(\trUdag_0)).
\end{equation*}

It is easy to see:
\begin{lemma}\label{lemma:preSymMain}
For any $\tlambda\in\tfkAPi^\prime$, 
\begin{equation*}
\chkpi_{|\prtrZnewparabctal
\cap\trK_\tlambda\trUdag_0}:\prtrZnewparabctal
\cap\trK_\tlambda\trUdag_0\to\chkprtrZ_\newpara(\tbhm,\checktal)
\cap\chktrK_\tlambda\chkpi(\trUdag_0)
\end{equation*}
is the $\bK$-linear isomorphism. In particular,
$\chkpi(\prtrZnewparabctal)=\chkprtrZ_\newpara(\tbhm,\checktal)$.
Furthermore, we have the injective $\bK$-linear map
${\check{\trHCmap}}^{\tbhm,\checktal}_\newpara:\chkprtrZ_\newpara(\tbhm,\checktal)\to\chktrUo$ such that
${\check{\trHCmap}}^{\tbhm,\checktal}_\newpara(\chkpi(X)):=\chkpi(\trHCmapbctaln(X))$
{\rm{(}}$X\in\prtrZnewparabctal${\rm{)}}.
\end{lemma}

{\it{Proof of Theorem~{\rm{\ref{theorem:SYMMAIN}}}.}}
By Theorem~\ref{theorem:mainTh} and Lemma~\ref{lemma:preSymMain}, 
we have Theorem~\ref{theorem:SYMMAIN}.
\hfill $\Box$

\subsection{Skew centers for quantum superalgebras associated to
basic classical Lie superalgebras} \label{subsection:SymExamples}
Assume that $\CardtrRptbtP<\infty$.
Assume that $\tbhm(\tal,\tal)\ne 1$ for all $\tal\in\trRptbhmtPi$.
Let $\tq\in\bKt$ be such that $\tq\ne 1$ and $\kpch(\tq)=0$. 
Assume that there exist a $\bZ$-module homomorphism
$\eparity:\tfkAPi\to\bZ$
and a $\bZ$-module bihomomorphism
$\langle\,,\,\rangle:\tfkAPi\times\tfkAPi\to\bZ$
such that $\tbhm(\tlambda,\tmu)=(-1)^{\eparity(\tlambda)\eparity(\tmu)}\tq^{\langle\tlambda,\tmu\rangle}$.
Let $\trRiso(\tbhm,\checktal):=\{\tal\in\trRtbhmtPi|\langle\tal,\tal\rangle= 0\}$ and
$\trRniso(\tbhm,\checktal):=\trRtbhmtPi\setminus\trRiso(\tbhm,\checktal)$.
Note 
$\trRiso(\tbhm,\checktal)=\{\tal\in\trRtbhmtPi|\tbhm(\tal,\tal)=-1\}$.
For $\tal\in\trRniso(\tbhm,\checktal)$, 
define the $\bZ$-module isomorphism $\espri_\tal
:\tfkAPi\to\tfkAPi$ by $\espri_\tal(\mu):=\tmu-
{\frac {2\langle\tmu,\tal\rangle} {\langle\tal,\tal\rangle}}\tal$.
Let $\eWprime$ be the subgroup of 
$\rmAutbZ(\tfkAPi)$ generated by 
$\espri_\tal$ with $\tal\in\trRniso(\tbhm,\checktal)$.
For $w^\prime\in\eWprime$,  define the $\bK$-linear
isomorphism $f_{w^\prime}:\chktrUo\to\chktrUo$ 
by $f_{w^\prime}(\chktrK_\tlambda):=\chktrK_{w^\prime(\tlambda)}$,
so $f_{w^\prime_1}\circ f_{w^\prime_2}=f_{w^\prime_1w^\prime_2}$,
i.e., this is a $\eWprime$-action on $\chktrUo$.
Assume that there exists a $\bZ$-module homomorphism $\trhoprime:\tfkAPi\to\bKt$ 
such that $\trhoprime(\lambda)^2=\trhobctal(\lambda)$
for all $\lambda\in\tfkAPi$.
Define the $\bK$-linear homomorphism $\phiprime:
\chktrUo\to\chktrUo$
by $\phiprime(\chktrK_\tlambda):=\trhoprime(-\tlambda)\chktrK_\tlambda$
($\tlambda\in\tfkAPi$).
Let $\newpara:\trtfkAPi\to\bKt$ be a $\bZ$-module homomorphism
with $\newpara(\trtfkAPi)\subset\{-1,1\}$.
Let $X=\sum_{\tlambda\in\tfkAPi}a_\tlambda
\chktrK_\tlambda\in\chktrUo$ with $a_\tlambda\in\bK$.
By Theorem~{\ref{theorem:SYMMAIN}},
we see that
$X\in\rmIm{\check{\trHCmap}}^{\tbhm,\checktal}_\newpara$ if and only if the following 
(${\check{e}}^\prime 1$)-(${\check{e}}^\prime 3$) hold.
\newline\par
(${\check{e}}^\prime 1$) Let  $\tlambda\in\tfkAPi$. Then
$a_\tlambda=0$ if there 
exists $\tbeta\in\trRniso(\tbhm,\tPi)$ such that 
$(-1)^{\eparity(\tbeta)\eparity(\tlambda+t\tbeta)}\tq^{\langle\tbeta,\tlambda+t\tbeta\rangle}
\ne\newpara(\tbeta)$ for all $t\in\bZ$.
\par
(${\check{e}}^\prime 2$) $f_{w^\prime}(\phiprime(X))=\phiprime(X)$
for all $w^\prime\in\eWprime$.
\par
(${\check{e}}^\prime 3$) Let $\tlambda\in\tfkAPi$ and $\tbeta\in\trRiso(\tbhm,\tPi)$.
Then
$\sum_{t=-\infty}^\pinfty (-1)^t\trhobctal(-t\tbeta)a_{\tlambda+2t\tbeta}
=0$ if
$(-1)^{\eparity(\tbeta)\eparity(\tlambda)}\tq^{\langle\tbeta,\tlambda\rangle}
\ne\newpara(\tlambda)$.
\newline\newline
Define the $\bZ$-module homomorphism $\newpara^\prime:\trtfkAPi\to\bKt$ 
by $\newpara^\prime(\tlambda):=\newpara(\tlambda)(-1)^{\eparity(\tlambda)}$
($\tlambda\in\tfkAPi$).
Let $\chktrU^\addsig$ be the $\bK$-algebra satisfying the condition that
$\chktrU$ is a $\bK$-subalgebra of $\chktrU^\addsig$
and
that there exists $\addsig\in\chktrU^y$ such that
$\addsig^2=1$, $\addsig\chktrK_\tlambda\addsig=\chktrK_\tlambda$
($\tlambda\in\tfkAPi$),
$\addsig\chktrE_i\addsig=(-1)^{\eparity(\talpii)}\chktrE_i$,
$\addsig\chktrF_i\addsig=(-1)^{\eparity(\talpii)}\chktrF_i$
($i\in\fkI$)
and $\chktrU^\addsig=\chktrU\oplus\chktrU\addsig$
as a $\bK$-linear space.
Let $\chkprtrZ_\newpara(\tbhm,\checktal)^\addsig$ be the $\bK$-subspace of $\chktrU^\addsig$
defined 
in the same way as that for $\prtrZnewparabctal$
with $\chktrU_\tmu\oplus(\chktrU_\tmu)\addsig$
in place of $\trU_\tmu$ ($\tmu\in\tfkAPi$).
Then $\chkprtrZ_\newpara(\tbhm,\checktal)^\addsig
=\chkprtrZ_\newpara(\tbhm,\checktal)\oplus\chkprtrZ_{\newpara^\prime}(\tbhm,\checktal)\addsig$.
Define the $\bK$-linear map $({\check{\trHCmap}}^{\tbhm,\checktal}_\newpara)^\addsig:
\chkprtrZ_\newpara(\tbhm,\checktal)^\addsig\to\chktrUo
\oplus(\chktrUo)\addsig$ by $({\check{\trHCmap}}^{\tbhm,\checktal}_\newpara)^\addsig(X+X^\prime\addsig):=
{\check{\trHCmap}}^{\tbhm,\checktal}_\newpara(X)
+{\check{\trHCmap}}^{\tbhm,\checktal}_{\newpara^\prime}(X^\prime)\addsig$
($X\in\chkprtrZ_\newpara(\tbhm,\checktal)$, $X^\prime\in\chkprtrZ_{\newpara^\prime}(\tbhm,\checktal)$).

Assume that the characteristic of $\bK$ is zero.
Let ${\bar{\trK}}_\tlambda:=\chktrK_\tlambda\addsig^{\eparity(\tlambda)}$
($\tlambda\in\tfkAPi$),
and ${\bar{\trE}}_i:=\chktrE_i$, ${\bar{\trF}}_i:=\chktrF_i\addsig^{\eparity(\talpii)}$
($i\in\fkI$).
Let ${\bar{\trU}}$ be the $\bK$-subalgebra of $\chktrU^\addsig$
generated by these elements, so 
$\chktrU^\addsig={\bar{\trU}}\oplus{\bar{\trU}}\addsig$
as a $\bK$-linear space.
Let ${\bar{\trU}}^0:=\oplus_{\tlambda\in\tfkAPi}\bK{\bar{\trK}}_\tlambda$.
For $t\in\bZ$, let $\tfkAPi[t]:=\{\tlambda\in\tfkAPi|(-1)^{\eparity(\tlambda)}
=(-1)^t\}$. Let ${\bar{\trU}}^0[t]:=\oplus_{\tlambda\in\tfkAPi[t]}\bK{\bar{\trK}}_\tlambda$,
so ${\bar{\trU}}^0={\bar{\trU}}^0[0]\oplus{\bar{\trU}}^0[1]$.
For $t\in\fkJ_{0,1}$,
let ${\bar{\trU}}_0[t]:=(\oplus_{\tlambda\in\tfkAPi[0]}\rmSpan_\bK(
\chkpi(\trUm_{-\tlambda}){\bar{\trU}}^0[t]\chkpi(\trUp_\tlambda)))
\oplus(\oplus_{\tlambda\in\tfkAPi[1]}\rmSpan_\bK(
\chkpi(\trUm_{-\tlambda}){\bar{\trU}}^0[t+1]\chkpi(\trUp_\tlambda)\addsig))$.
Then 
${\bar{\trU}}\cap\chkprtrZ_\newpara(\tbhm,\checktal)^\addsig=
({\bar{\trU}}_0[0]\cap\chkprtrZ_\newpara(\tbhm,\checktal))
\oplus({\bar{\trU}}_0[1]\cap\chkprtrZ_{\newpara^\prime}(\tbhm,\checktal)\addsig)$.

Let $X=\sum_{\tlambda\in\tfkAPi}a_\tlambda
{\bar{\trK}}_\tlambda$ with $a_\tlambda\in\bK$.
We see that $X\in{\bar{\trU}}\cap\chkprtrZ_\newpara(\tbhm,\checktal)^\addsig$
if and only if it satisfies the conditions 
(${\check{e}}^\prime 1^\prime$)-(${\check{e}}^\prime 3^\prime$)
obtained from (${\check{e}}^\prime 1$)-(${\check{e}}^\prime 3$)
by replacing $\newpara$ with $\newpara^\prime$.
By \cite{Hec09}
(see also \cite{AYY15}),
${\bar{\trU}}$ is isomorphic to a quantum superalgebra $U_q({\mathfrak{s}})$
of a finite dimensional contragredient Lie superalgebra ${\mathfrak{s}}$.
Let ${\bar {R}}$ be a set of roots of ${\mathfrak{s}}$.
Then
${\bar {R}}$ can be identified with $\trRtbhmtPi
\cup\{2\tbeta|\tbeta\in\trRtbhmtPi,\,\langle\tbeta,\tbeta\rangle\ne 0,\,
\eparity(\tbeta)\in 2\bZ-1\}$.
Let ${\bar {R}}^+:={\bar {R}}\cap\tfkAPip$.
Let $b:=\sum_{\tbeta\in{\bar {R}}^+}(-1)^{\eparity(\tbeta)}\tbeta$.
Then $\trhobctal(\tlambda)=(-1)^{\eparity(\tlambda)}\tq^{\langle b,\tlambda\rangle}$
($\lambda\in\tfkAPi$) (see also \cite[Corollary~8.5.4]{Musson12}).
Let ${\bar {R}}^\prime:=\{\tbeta\in{\bar {R}}|\eparity(\tbeta)\in 2\bZ\}$.
Then ${\bar {R}}^\prime$
can be identified with the root system of 
the complex reductive Lie algebra realized as the even part of ${\mathfrak{s}}$.
Let $\{\tal^\prime_i|i\in\fkI^\prime\}$ be a base of ${\bar {R}}^\prime$,
where $\fkI^\prime$ is a non-empty subset of $\fkI$.
For $i\in\fkI$, 
$\tal^{\prime\prime}_i\in\trRtbhmtPi$ be such that
$r\tal^{\prime\prime}_i=\tal^\prime_i$ for some $r\in\fkJ_{1,2}$.
Assume $X\in{\bar{\trU}}\cap\chkprtrZ_\newpara(\tbhm,\checktal)^\addsig$. 
By (${\check{e}}^\prime 1^\prime$), 
for $\tlambda\in\tfkAPi$, if $a_\tlambda\ne 0$, then
${\frac {\langle\tal^{\prime\prime}_i,\tlambda\rangle} 
{\langle\tal^{\prime\prime}_i,\tal^{\prime\prime}_i\rangle}}\in\bZ$
for all $i\in\fkI^\prime$.

\begin{remark} \label{remark:gooddirection}
It may be essential to study for obtaining a result for $\trU$ similar to
that of \cite{SV2011} for the basic classical Lie superalgebras. 
If such a study is achieved, we will have a $\bK$-basis of $\prtrZnewparabctal$
in terms of the Grothendieck ring of $\trU$.
\end{remark}

{\bf{Acknowledgment.}} This research was in part supported by
Japan's Grand-in-Aid for Scientific Research (C), 16K05095,
and the one 25400040.

\vspace{1cm}
{Punita Batra, 
Harish-Chandra Research Institute,
Chhatnag Road,
Jhunsi,
Allahabad 211 019,
India}, \newline
{E-mail: batra@hri.res.in}
\vspace{0.5cm}

{Hiroyuki Yamane, Department of Mathematics,
Faculty of Science, University of Toyama,
3190 Gofuku, Toyama-shi, Toyama 930-8555, Japan}, \newline
{E-mail: hiroyuki@sci.u-toyama.ac.jp}
\end{document}